\documentclass[12pt]{article}
\usepackage{amsfonts}
\usepackage{latexsym,amssymb,amsfonts}
\usepackage{mathrsfs}
\usepackage{graphicx}
\usepackage{psfrag}
\usepackage{amsmath, amsbsy}
\usepackage{amsopn, amstext}
\usepackage{amsmath,multicol,amsopn}
\usepackage{latexsym,amssymb}

\def\sqr#1#2{{\vcenter{\vbox{\hrule height.#2pt
              \hbox{\vrule width.#2pt height#1pt \kern#1pt \vrule width.#2pt}
              \hrule height.#2pt}}}}
\def\signed #1{{\unskip\nobreak\hfil\penalty50
              \hskip2em\hbox{}\nobreak\hfil#1
              \parfillskip=0pt \finalhyphendemerits=0 \par}}
\def\endpf{\signed {$\sqr69$}}

\def\dbR{{\mathop{\rm l\negthinspace R}}}
\def\dbC{{\mathop{\rm l\negthinspace\negthinspace\negthinspace C}}}

\def\3n{\negthinspace \negthinspace \negthinspace }
\def\2n{\negthinspace \negthinspace }
\def\1n{\negthinspace }

\def\dbC{{\mathbb{C}}}

\def\dbE{{\mathbb{E}}}
\def\dbF{{\mathbb{F}}}

\def\dbH{{\mathbb{H}}}

\def\dbN{{\mathbb{N}}}

\def\dbP{{\mathbb{P}}}

\def\dbR{{\mathbb{R}}}

\def\rank{{\rm rank}}

\def\ds{\displaystyle}

\def\={\buildrel \triangle \over =}

\def\resp{{\it resp. }}
\def\a{\alpha}
\def\b{\beta}

\def\d{\delta}
\def\e{\varepsilon}

\def\l{\lambda}

\def\n{\nabla}
\def\si{\sigma}
\def\t{\times}
\def\f{\varphi}
\def\th{\theta}
\def\om{\omega}

\def\ns{\noalign{\ss} }
\def\nm{\noalign{\ms} }

\def\D{\Delta}

\def\L{\Lambda}
\def\Si{\Sigma}

\def\Om{\Omega}

\def\ol{\overline}
\def\dbm{\mathbf{m}}

\def\cA{{\cal A}}
\def\cB{{\cal B}}
\def\cC{{\cal C}}

\def\cF{{\cal F}}

\def\cH{{\cal H}}
\def\cI{{\cal I}}
\def\cJ{{\cal J}}

\def\cL{{\cal L}}
\def\cM{{\cal M}}

\def\cQ{{\cal Q}}

\def\cU{{\cal U}}

\def\cl{{\cal l}}
\def\mE{{\mathbb{E}}}

\def\no{\noindent}

\def\ss{\smallskip}
\def\ms{\medskip}

\def\q{\quad}
\def\qq{\qquad}
\def\hb{\hbox}

\def\lan{\mathop{\langle}}
\def\ran{\mathop{\rangle}}

\def\max{\mathop{\rm max}}
\def\min{\mathop{\rm min}}
\def\exp{\mathop{\rm exp}}
\def\sup{\mathop{\rm sup}}

\def\h{\widehat}
\def\wt{\widetilde}
\def\cd{\cdot}
\def\cds{\cdots}

\def\inf{\hbox{\rm inf$\,$}}

\def\ae{\hbox{\rm a.e.{ }}}
\def\as{\hbox{\rm a.s.{ }}}

\def\span{\hbox{\rm span$\,$}}
\def\tr{\hbox{\rm tr$\,$}}
\def\cl{\overline}

\def\Re{{\mathop{\rm Re}\,}}

\def\|{\Big |}
\def\({\Big (}
\def\){\Big )}
\def\[{\Big[}
\def\]{\Big]}
\def\be{\begin{equation}}
\def\bel{\begin{equation}\label}
\def\ee{\end{equation}}
\def\bt{\begin{theorem}}
\def\bcd{\begin{condition}}
\def\ecd{\end{condition}}
\def\et{\end{theorem}}
\def\bc{\begin{corollary}}
\def\ec{\end{corollary}}
\def\bde{\begin{definition}}
\def\ede{\end{definition}}
\def\bl{\begin{lemma}}
\def\el{\end{lemma}}
\def\bp{\begin{proposition}}
\def\ep{\end{proposition}}
\def\br{\begin{remark}}
\def\er{\end{remark}}
\def\ba{\begin{array}}
\def\ea{\end{array}}
\def\ed{\end{document}}
\def\ns{\noalign{\ms}}
\def\ds{\displaystyle}

\def\square#1{\vbox{\hrule\hbox{\vrule height#1%
     \kern#1\vrule}\hrule}}
\def\rectangle#1#2{\vbox{\hrule\hbox{\vrule height#1%
     \kern#2\vrule}\hrule}}

\font\tenbb=msbm10 \font\sevenbb=msbm7
\font\fivebb=msbm5

\newfam\bbfam
\scriptscriptfont\bbfam=\fivebb
\textfont\bbfam=\tenbb
\scriptfont\bbfam=\sevenbb

\newtheorem{lemma}{Lemma}[section]
\newtheorem{remark}{Remark}[section]
\newtheorem{example}{Example}[section]
\newtheorem{theorem}{Theorem}[section]
\newtheorem{corollary}{Corollary}[section]

\newtheorem{definition}{Definition}[section]
\newtheorem{proposition}{Proposition}[section]
\newtheorem{condition}{Condition}[section]

\makeatletter
   
   \@addtoreset{equation}{section}
\makeatother

\begin{document}
\title{\bf A Mini-Course on Stochastic Control\thanks{This is a lecture notes of a short introduction to stochastic control. It was written for the LIASFMA (Sino-French International Associated  Laboratory for Applied Mathematics) Autumn School ``Control and Inverse Problems of Partial Differential Equations" at Zhejiang University, Hangzhou, China from October 17 to October 22, 2016. The second named author thanks Professors Jean-Michel Coron and Tatsien Li for their kind invitation, and Professor Gang Bao and his team for their hospitality during the teaching of this course.} \ms}

\author{Qi L\"{u}\thanks{School of Mathematics, Sichuan University, Chengdu 610064, Sichuan Province, China. The research of this
author is partially supported by NSF of China
under grants 11471231, the Fundamental Research
Funds for the Central Universities in China
under grant 2015SCU04A02 and Grant MTM2014-52347
of the MICINN, Spain. {\small\it E-mail:}
{\small\tt lu@scu.edu.cn}.} ~~~ and~~~ Xu
Zhang\thanks{School of Mathematics, Sichuan
University, Chengdu 610064, Sichuan Province,
China. The research of this author is partially
supported by the NSF of China under grants
11231007, the PCSIRT
under grant IRT$\!\_$15R53 and the Chang Jiang Scholars Program
from the Chinese Education Ministry. {\small\it
E-mail:} {\small\tt zhang$\_$xu@scu.edu.cn}.}}

\date{}

\maketitle

\begin{abstract}\no
 This course is addressed to giving a short introduction to control theory of stochastic systems, governed by stochastic differential equations in both finite and infinite dimensions. We will mainly explain the new phenomenon and difficulties in the study of controllability and optimal control problems for these sort of equations. In particular, we will show by some examples that both the formulation of stochastic control problems and the tools to solve them may differ considerably from their deterministic counterpart.

\end{abstract}


\tableofcontents


\section{Introduction}


It is well-known that control theory was founded
by N. Wiener in 1948. After that, this theory
was greatly extended to various complicated
settings and widely used in sciences and
technologies.

Clearly, ``control" means a suitable manner for people to change the dynamics of a system under consideration. There are two fundamental issues in control theory. One is ``feasibility", or in the language of control theory, controllability, which means that one may find at least one way to achieve a goal. Another is ``optimality", or optimal control, which indicates that, one hopes to find the best way, in some sense, to achieve the goal.

Roughly speaking, control theory can be divided into two parts. The first part is control theory for deterministic systems, and the second part is that for stochastic systems.
Of course, these two parts are not completely separated but rather they are inextricably linked each other.

Control theory for deterministic systems  can be
again divided into two parts. The first part is
control theory for finite dimensional systems,
mainly governed by ordinary  differential
equations, and the second part is that for
(deterministic) distributed parameter systems,
mainly described  by differential equations in
infinite dimensional spaces, typically by
partial differential equations. Control theory
for finite dimensional systems is by now
relatively mature. There exist a huge list of
works on control theory for distributed
parameter systems but it is still quite active.

Likewise, control theory for stochastic  systems
can be divided into two parts. The first part is
control theory for stochastic finite dimensional
systems, governed by stochastic (ordinary)
differential equations, and the second part is
that for stochastic distributed parameter
systems, described by stochastic differential
equations in infinite dimensions, typically by
stochastic partial differential equations.

One can find a huge list of publications on control theory for stochastic finite dimensional systems and its applications, say, in mathematical finance. Nevertheless, most of the existing works in this respect are mainly addressed/related to the optimal control problems. As we shall see later in this course, so far controllability theory for stochastic finite dimensional systems is NOT well-developed.

Control theory for stochastic distributed parameter systems, is, in
our opinion, still at its very beginning stage. This is actually a rather new branch of mathematical control theory, which is indeed the main concern of this course (See \cite{LZ3, LZ4} for more material).

One of the most essential difficulties in the study of control theory
for stochastic distributed parameter systems is that, compared to the deterministic setting, people know very little about stochastic evolution equation (and in particular, about stochastic partial differential equations) although significant progresses have been made there, especially in recent years. On the other hand, as we shall show in this course, both the formulation of stochastic control problems in infinite dimensions and
the tools to solve them may differ considerably from their
deterministic/finite-dimensional counterparts. Because of this, one has
to develop new mathematical tools to solve some problems in this field.

The rest of this course is organized as follows. In Section \ref{s2}, we collect some preliminary results (without proofs) from probability theory and stochastic
analysis. In Sections \ref{sec ex finite} and \ref{s3}, we analyze respectively the controllability and optimal controls for stochastic differential equations in finite dimensions; while in Sections \ref{s5}--\ref{ch-Pon}, we consider the same problems but for stochastic evolution equations in infinite dimensions.

\section{Some preliminary results from probability theory and stochastic
analysis}\label{s2}

For the proofs of the results presented in this section, we refer to \cite{Prato, LZ3}. In what follows, we shall denote by $\cC$ a generic positive constant, which may change from one place
to another.

\subsection{Probability,
random variables and expectation}
Fix a nonempty set $\Omega$ and a  $\si$-field $\cF$ on $\Omega$.
Let $(\Omega,\cF, \dbP)$ be a complete probability space, i.e. a complete measure space  for
which $\dbP(\Omega)=1$.
Any point $\omega\in\Omega$ is called a sample, any
$A\in\cF$ is called an event and $\dbP(A)$ represents
the probability of event $A$. If an event
$A\in\cF$ is such that $\dbP(A)=1$, then we may
alternatively say that $A$ holds, $\dbP$-a.s.,
or simply $A$ holds a.s.

Let $H$ be a Hilbert space. Each $H$-valued, strongly measurable
function $f:\ (\Omega,\cF)\to (H,\cB(H))$ is called an ($H$-valued)
random variable. Clearly, $f^{-1}(\cB(H))$ is a
sub-$\si$-field of $\cF$, which is called the $\si$-field generated
by $f$, denoted by $\si(f)$. Further,  if $f$ is Bochner integrable w.r.t. the measure
$\dbP$, i.e. the integral
 $$
 \dbE f\equiv\int_\Omega f(\omega)d\, \dbP(\omega)
 $$
exists, then we say that $f$ has a mean. We also call $\dbE f$ the
(mathematical) expectation of $f$.
For a given index set $\Lambda$ and a
family of $H$-valued,  random variables $\{f_\l\}_{\l\in
\Lambda}$ (defined on $(\Omega,\cF)$), we denote by
$\si(f_\l;\l\in \Lambda)$ the $\si$-field
generated by $\cup_{\l\in \Lambda}\si(f_\l)$.

For any $p\in [1,\infty)$, denote by
$L_{\mathcal F}^p(\Omega;H)\equiv
L^p(\Omega,{\mathcal F},\dbP;H)$ the set of all
random variables $f$ such  that $|f|_H^p$ has
means. It is a Banach space with the norm $
|f|_{L_{\mathcal F}^p(\Omega)}=\left(\int_\Omega
|f|_H^pd\dbP\right)^{1/p}$. In particular,
$L_{\mathcal F}^2(\Omega;H)$ is a Hilbert space.
We simply denote $L_\cF^p(\Omega;\dbR)$ by
$L_\cF^p(\Omega)$. For any $f\in L_{{\mathcal
F}}^2(\Omega)$, we define the variance of $f$ by
 $$
 \hb{Var}\; f=\dbE(f-\mE f)^2.
 $$

Let $A,B\in\cF$. We say that $A$ and $B$ are
independent if $\dbP(A\cap B)=\dbP(A)\dbP(B)$.
Let $\cJ_1$ and $\cJ_2$ be two subsets of $\cF$.
We say that $\cJ_1$ and $\cJ_2$ are independent
if $\dbP(A\cap B)=\dbP(A)\dbP(B)$ for any
$A\in\cJ_1$ and $B\in\cJ_2$. Let $f,\ g:\
(\Omega,\cF)\to(H,\cB(H))$ be two random
variables. We say that $f$ and $g$ (\resp $f$
and $\cJ_1$) are independent if $\si(f)$ and
$\si(g)$ (\resp $\si(f)$ and $\cJ_1$) are
independent.

Let $X:\ (\Omega,\cF)\to (\dbR,\cB(\dbR))$ be a
random variable. We call
 $
 F(x)\equiv \dbP\{X\le x\}
 $
the distribution function of $X$. If for some
function $p(\cd)$, one has
 $$
 F(x)=\int_{-\infty}^{x}p(\xi)d\xi,
 $$
then the function $p(\cd)$ is called the density
of $X$. If $p(\cdot)$ is of the following form:
 $$
 p(x)=(2\pi
 \mu)^{-1/2}\exp\left\{-\frac{1}{2\mu}(x-\l)^2\right\},
 $$
where $\l\in\dbR$, $\mu\in \dbR^+$, then $X$ is
called a normally distributed random variable (or
$X$ is a normal distribution). Clearly,
$\l$ and $\mu$ are the mean and variance of $f$,
respectively.

Assume that $\cJ\subset{\mathcal F}$ is a given sub-$\si$-field and $f\in L_\cF^1(\Omega;H)$. Define a function on $\cJ$
by
 $$
 \nu(B)=\int_Bfd\dbP, \qq\forall\;B\in\cJ.
 $$
It is easy to see that $\nu$ is an ($H$-valued)
vector measure of bounded variation on
$(\Omega,\cJ)$, and $\nu(B)=0$ whenever $\dbP(B)=0$. Hence, there is a (unique) function in
$L_\cJ^1(\Omega;H)$,
denoted by $\mE(f\;|\;\cJ)$, such that
  \bel{1c1e14}
  \int_B\dbE(f\,|\,\cJ)d\dbP=\int_Bfd\dbP, \qq\forall\;B\in\cJ.
  \ee
This function is called the conditional
expectation of $f$ given $\si$-field $\cJ$.

\begin{example} Let $B_1,B_2\subset{\mathcal F}$ such that
$B_1\cup B_2=\Omega$,  $B_1\cap B_2=\emptyset$, and $\dbP(B_k)>0$
for all $k=1,2$. Let $f\in L_\cF^1(\Omega;H)$ and $\cJ=\{\emptyset,\Omega,B_1,B_2\}$. Then
 $$
 \dbE(f\,|\,\cJ)(\omega)=\sum_{k=1}^2 \frac{1}{
 \dbP(B_k)}\int_{B_k}fd\dbP \chi_{B_k}(\omega).
 $$
\end{example}

We collect some basic properties of  conditional
expectation as follows.

\begin{theorem}
Let $\cJ$ be a sub-$\si$-field of $\cF$ and $f\in L_\cF^1(\Omega;H)$. It
holds that:

{\rm 1)} The map $\mE(\cd\;|\;\cJ):\
L_{\cF}^1(\Omega;H)\to L_{\cJ}^1(\Omega;H)$ is
linear and continuous;

{\rm 2)} $\mE(a\;|\;\cJ)=a$,
$\dbP|_{\cJ}\mbox{-}\as\!\!$, $\forall\;a \in H$;

{\rm 3)} If $\a \in L_{\cJ}^1(\Omega)$ satisfies
$\a f\in L_{\cF}^1(\Omega;H)$, then
$$
\mE(\a f\;|\;\cJ)=\a \mE(f\;|\;\cJ),\qq
\dbP|_{\cJ}\mbox{-}\as
$$
In particular, $\mE(\a\;|\;\cJ)=\a$,
$\dbP|_{\cJ}-\as\!\!$;

{\rm 4)} If $f$ is independent of $\cJ$, then
$$
\mE(f\;|\;\cJ)= \mE f,\qq
\dbP|_{\cJ}\mbox{-}\as\!\!;
$$

{\rm 5)} Let $\cJ'$ be a sub-$\si$-field of
$\cJ$. Then
$$
\mE( \mE(f\;|\;\cJ)\;|\;\cJ')=  \mE(
\mE(f\;|\;\cJ')\;|\;\cJ)=  \mE(f\;|\;\cJ'), \qq
\dbP|_{\cJ'}\mbox{-}\as\!\!;
$$

{\rm 6)} (Jensen's
inequality) Let $\phi:\ H\to \dbR$ be a convex
function such that $\phi(f)\in L^1_\cF(\Omega)$.
Then
$$
\phi( \mE(f\;|\;\cJ))\leq \mE(\phi(f)\;|\;\cJ),
\qq \dbP|_{\cJ}\mbox{-}\as
$$
In particular, for any $p\geq 1$,
$$
\big| \mE(f\;|\;\cJ)\big|_H^p\leq
\mE(|f|_H^p\;|\;\cJ), \qq \dbP|_{\cJ}\mbox{-}\as
$$
provided that $ \mE|f|_H^p$ exists.
\end{theorem}

\subsection{Stochastic processes}

Let $\cI=[0,T]$ with $T>0$. A family of $H$-valued random variables
$\{X(t)\}_{t\in \cI}$ is called a stochastic
process. For any $\omega\in\Omega$, the map $t\mapsto X(t,\omega)$ is called a
sample path (of $X$). We will interchangeably use $\{X(t)\}_{t\in
\cI}$, $X(\cd)$ or even $X$ to denote a (stochastic) process.

An
($H$-valued) process $ X(\cdot)$ is said to be
continuous
({\it resp.}, c\'adl\`ag, i.e., right-continuous with left
limits) if there is a $\dbP$-null
set $N \in\cF$, such that for any
$\omega\in\Omega\setminus N$, the sample path
$X(\cdot,\omega)$ is continuous (\resp\!\!,
c\'adl\`ag) in $H$. In a similar way, one can define
right-continuous stochastic processes,
etc. Two ($H$-valued) processes $\ X(\cdot)\ $ and $\
\cl X(\cdot)\ $ are said to be \index{stochastic
equivalence} stochastically equivalent if
$\dbP(\{X(t)=\cl X(t)\})=1$  for any $t\in \cI$.
In this case, one is said to be a modification
of the other.

We call a family of sub-$\si$-fields
$\{\cF_t\}_{t\in \cI}$ in $\cF$ a
\index{filtration} filtration if
$\cF_{t_1}\subset \cF_{t_2}$ for all $
t_1,t_2\in \cI \mbox{ with } t_1\leq t_2$. For
any $t\in \cI$, we put
$$
\cF_{t+}\=\bigcap_{s\in (t,+\infty)\cap
\cI}\cF_s,\qq \cF_{t-}\=\bigcup_{s\in [0,t)\cap
\cI}\cF_s.
$$
If $\cF_{t+}=\cF_t$ (\resp $\cF_{t-}=\cF_t$),
then $\{\cF_t\}_{t\in \cI}$ is said to be right
(\resp left) continuous. In the sequel, for
simplicity, we write $\mathbf{F}=\{\cF_t\}_{t\in
\cI}$ unless we want to emphasize what $\cF_t$
or $I$ exactly is. We call
$(\Omega,\cF,\mathbf{F}, \dbP)$ a filtered
probability space.

We say that
$(\Omega,\cF,\mathbf{F}, \dbP)$ satisfies
 the usual condition if
$(\Omega,\cF,\dbP)$ is complete, $\cF_0$
contains all $\dbP$-null sets in $\cF$, and
$\mathbf{F}$ is right continuous. We shall keep
these assumptions in what follows unless stated
otherwise.

\begin{definition}
Let $X(\cd)$ be an $H$-valued process.

{\rm 1)} \index{measurable stochastic process}
$X(\cd)$ is said to be measurable if the map
$(t,\omega)\mapsto X(t,\omega)$ is strongly
$(\cB(\cI)\t \cF)/\cB(H)$-measurable;

{\rm 2)} \index{adapted stochastic process}
$X(\cd)$ is said to be $\mathbf{F}$-adapted if
it is measurable, and for each $t\in \cI$, the
map $\omega\mapsto X(t,\omega)$ is strongly
$\cF_t/\cB(H)$-measurable;

{\rm 3)}  $X(\cd)$ is said to be
$\mathbf{F}$-progressively measurable if for
each $t\in \cI$, the map $(s,\omega)\mapsto
X(s,\omega)$ from $[0,t]\times\Omega$ to $H$ is
strongly $(\cB([0,t])\t
\cF_t)/\cB(H)$-measurable.
\end{definition}

A set $A\in \cI\times \Omega$ is called
progressively measurable w.r.t. $\mathbf{F}$ if
the process $\chi_{A}(\cdot)$ is progressive.
The class of all progressively measurable sets
is a $\si$-field, called the
progressive
$\si$-field w.r.t. $\mathbf{F}$, denoted by
$\dbF$. One can show that, an ($H$-valued) process $\f:[0,T]\times\Omega\to
H$ is $\mathbf{F}$-progressively measurable if
and only if it is strongly $\dbF$-measurable.

It is clear that if $X(\cd)$ is
$\mathbf{F}$-progressively measurable, it must
be $\mathbf{F}$-adapted. Conversely, it can be
proved that, for any $\mathbf{F}$-adapted
process $X(\cd)$, there is an
$\mathbf{F}$-progressively measurable process
$\wt X(\cd)$ which is stochastically equivalent
to $X(\cd)$. For
this reason, in the
sequel, by saying that a process $X(\cd)$ is
$\mathbf{F}$-adapted, we mean that it is
$\mathbf{F}$-progressively measurable.

For any
$p,q\in[1,\infty)$, write
$$
\begin{array}{ll}
\ds L^p_\dbF(\Omega;L^q(0,T;H)) \=\! \Big\{ \f:
(0,T)\!\times\! \Omega \!\to\! H\,\Big|\,
\f(\cd)\hb{
is $\dbF$-adapted and }\dbE\(\int_0^T\!|\f(t)|_H^qdt\)^{\frac{p}{q}}\!<\!\infty\Big\},\\
\ns\ds L^q_\dbF(0,T;L^p(\Omega;H)) \=\! \Big\{
\f :(0,T)\! \times\! \Omega\!\to\! H  \Big|
\f(\cd)\hb{ is $\dbF$-adapted and
}\int_0^T\!\(\dbE|\f(t)|_H^p\)^{\frac{q}
{p}}dt<\!\infty\Big\}.
\end{array}
$$
Similarly, we may also define (for $1\le
p,q<\infty$)
$$
\left\{
\begin{array}{ll}
\ds L^\infty_\dbF(\Omega;L^q(0,T;H)),\q
L^p_\dbF(\Omega;L^\infty(0,T;H)),\q
L^\infty_\dbF
(\Omega;L^\infty(0,T;H)),\\
\ns\ds L^\infty_\dbF(0,T;L^p(\Omega;H)),\q
L^q_\dbF(0,T;L^\infty(\Omega;H)),\q
L^\infty_\dbF(0,T;L^\infty(\Omega;H)).
\end{array}
\right.
$$
All these spaces are Banach spaces
(with the canonical norms). In the sequel, we
shall simply denote
$L^p_\dbF(\Omega;L^p(0,T;H))\equiv
L^p_\dbF(0,T;L^p(\Omega;H))$ by
$L^p_\dbF(0,T;H)$; and further simply denote
$L^p_\dbF(0,T;\dbR)$ by $L^p_\dbF(0,T)$.

For any $p\in [1,\infty)$, set
$$
\begin{array}{ll}\ds
L^{p}_{\dbF}(\Omega;C([0,T];H))\=\Big\{\f:[0,T]\times\Omega\to
 H \,\Big|\, \f(\cd)\hb{
is continuous, }\\
\ns\ds\qq\qq\qq\qq \qq\qq
\mathbf{F}\mbox{-adapted and
}\dbE\big(|\f(\cdot)|_{C([0,T];H)}^p\big)<\infty\Big\}
\end{array}
$$
and
$$
\begin{array}{ll}\ds
C_{\dbF}([0,T];L^{p}(\Omega;H))\=\Big\{\f:[0,T]\times\Omega\to
H\,\Big|\,\f(\cd)\hb{
is $\mathbf{F}$-adapted }\\
\ns\ds\qq\qq\qq\qq \qq \mbox{ and }\f(\cd):[0,T]
\to L^p_{\cF_T}(\Omega;H)\mbox{ is
continuous}\Big\}.
\end{array}
$$
One can show that both
$L^{p}_{\dbF}(\Omega;C([0,T];H))$ and
$C_{\dbF}([0,T];L^{p}(\Omega;H))$ are Banach
spaces with canonical norms
$|\f(\cd)|_{L^{p}_{\dbF}(\Omega;C([0,T];H))}\!=\!\big(\dbE(|\f(\cdot)|_{C([0,T];H)}^p)\!\big)^{1/p}$
and
$|\f(\cd)|_{C_{\dbF}([0,T];L^{p}(\Omega;H))}=\max_{t\in
[0,T]}\big(\dbE(|\f(t)|_H^p)\big)^{1/p}$,
respectively.
Also, we denote by
$D_{\dbF}([0,T];L^p(\Omega;H))$ the Banach space
of all processes such that $X(t)$ is c\`adl\`ag
in $L^p_{\cF_T}(\Omega;H)$, w.r.t. $t\in [0,T]$,
such that
$|\mathbb{E}|X(\cdot)|^p_H|^{1/p}_{L^{\infty}(0,T)})
< \infty$, with the canonical norm.

We need to introduce two important classes of
stochastic processes, i.e., Brownian motion and
martingale.

\begin{definition}
A continuous $\mathbf{F}$-adapted process $W(\cd)$ is called a
$1$-dimensional Brownian motion
(over $\cI$), if for all $s,t\in\cI$ with $0\le s<t<T$, $W(t)-W(s)$ is independent of
${\mathcal F}_s$, and normally distributed with mean $0$ and
variance $t-s$. In addition, if $\dbP(W(0)=0)=1$, then $W(\cd)$ is
called a $1$-dimensional standard Brownian motion.
\end{definition}

In the seuqel, we fix a $1$-dimensional standard Brownian motion on
$(\Omega,\cF,\mathbf{F}, \dbP)$. Write
 \begin{equation}\label{a2.2}
 \cF_t^W\=\si(W(s);\; s\in [0,t])\subset\cF_t,\qq\forall\;t\in\cI.
 \end{equation}
Generally, the filtration
$\{\cF_t^W\}_{t\in\cI}$ is left-continuous, but
not necessarily right-continuous. Nevertheless,
the augmentation $\{\hat \cF_t^W\}_{t\in\cI}$ of
$\{ \cF_t^W\}_{t\in\cI}$ by adding all
$\dbP$-null sets is continuous, and $W(\cd)$ is
still a Brownian motion on the (augmented)
filtered probability space
$(\Omega,\cF,\{\hat\cF_t^W\}_{t\in\cI}, \dbP)$.
In the sequel, by saying that $\mathbf{F}$ is
the natural filtration generated by $W(\cd)$, we
mean that $\mathbf{F}$ is generated as in
(\ref{a2.2}) with the above augmentation, and
hence in this case $\mathbf{F}$ is continuous.

\begin{definition}
An $H$-valued, $\mathbf{F}$-adapted process
$X=\{X(t)\}_{t\in \cI}$ is called an
$\mathbf{F}$-martingale, if
$X(t)$ is Bochner integrable for each
$t\in \cI$, and
$E(X(t)\;|\;\cF_s)=X(s)\ \as\!\!$, for
any $t,s\in \cI$ with $s< t$.
\end{definition}

Clearly, for any $f\in L_{\mathcal F}^1(\Omega;H)$, the process
$\{\dbE(f\;|\;{\mathcal F}_t)\}_{t\in \cI}$ is an
${\mathbb{F}}$-martingale.

Write
$$\ba{ll}
\cM^2[0,T]\!=\bigm\{\!X\!\in\!
L^2_{{\mathbb{F}}}(0,T;H)\!\bigm|\!X\hb{ is a
right-continuous, ${\mathbb{F}}$-martingale with
}\!X(0)\!=\!0,\dbP\mbox{-a.s.}\!\bigm\},\\
 \ns
\cM^2_c[0,T]=\bigm\{X\in\cM^2[0,T]\bigm| X \hb{ is continuous}\bigm\}.\ea
 $$
Define
$$|X|_{\cM^2[0,T]}=\sqrt{\dbE|X(T)|_H^2},\qq\forall\; X\in\cM^2[0,T].
$$
Then, $(\cM^2[0,T],|\cd|_{\cM^2[0,T]})$ is a
Hilbert space, and $\cM^2_c[0,T]$ is a closed subspace of
$\cM^2[0,T]$.

\subsection{It\^o's integral and its properties}

We now define the It\^o integral
 \bel{a2.3}
 \int_0^TX(t)dW(t)
 \ee
of an $H$-valued, $\mathbf{F}$-adapted stochastic process $X(\cd)$ (satisfying suitable conditions) w.r.t. a Brownian motion
$W(t)$. Note that one cannot define (\ref{a2.3}) to be a
Lebesgue-Stieltjes type integral by regarding $\omega$ as a parameter. Indeed, the map
$t\ni[0,T]\mapsto W(t, \cd)$ is nowhere
differentiable, $\dbP$-a.s.

Denote by $\cL_0$  the class of simple processes
$f\in L^2_{{\mathbb{F}}}(0,T;H)$ of the forms:
 \bel{6.2}
 \displaystyle f(t,\omega)=\sum_{j=0}^n
 f_j(\omega)\chi_{[t_j,t_{j+1})}(t), \qq (t,\omega)\in[0,T]\t\Omega,
 \ee
where $0=t_0<t_1<\cdots<t_{n+1}= T$, $f_j$ is
${\mathcal F}_{t_j}$-measurable with
$\sup\big\{|f_j(\omega)|_H\;\big|\;
j\in\{0,\cdots,n\},\,\omega\in\Omega\big\}<\infty$.
One can show that $\cL_0$ is dense in
$L_{{\mathbb{F}}}^2(0,T;H)$.

We now define the It\^o integral (\ref{a2.3}) as a mapping
$
 f\in L_{{\mathbb{F}}}^2(0,T)\mapsto I(f)\in \cM_c^2[0,T]$.
First, assume that $f\in\cL_0$ takes the form of (\ref{6.2}). Then
we set
 \bel{6.7}
 I(f)(t,\omega)=\sum_{j=0}^n f_j(\omega)[W(t\wedge
 t_{j+1},\omega)-W(t\wedge t_j,\omega)].
 \ee
It is easy to show that $I(f )\in \cM_c^2[0,T]$ and the following It\^o isometry holds:
 \bel{6oo7}|I(f
)|_{\cM^2[0,T]}=|f
|_{L_{{\mathbb{F}}}^2(0,T;H)}.
 \ee
Generally, for $f \in
L_{{\mathbb{F}}}^2(0,T;H)$, one can find a
sequence of $\{f _k\}\subset \cL_0$ such that
$|f _k-f |_{L_{{\mathbb{F}}}^2(0,T;H)}\to 0$ as
$k\to\infty$. Since $|I(f _k)-I(f
_j)|_{\cM^2[0,T]}=|f _k-f
_j|_{L_{{\mathbb{F}}}^2(0,T;H)}$, one deduces
that $\{I(f_k)\}_{k=1}^\infty$ is a Cauchy
sequence in $\cM^2[0,T]$ and therefore, it
converges to a unique element $X\in \cM^2[0,T]$.
Clearly, $X$ is determined uniquely by $f$ and
is independent of the particular choice of $\{f
_k\}_{k=1}^\infty$. This process is called the
It\^o integral of $f \in
L_{{\mathbb{F}}}^2(0,T;H)$ w.r.t. the Brownian
Motion $W(\cd)$. We shall denote it by
 $
 \int_0^tf (s)dW(s)$ or simply
 $\int_0^tf  dW.
 $

\begin{theorem}
Let $f , g \in L_\dbF^2(0,T;H)$, $a,b\in
 L_{\cF_s}^2(\Omega)$, $T\ge t>s\ge0$. Then

{\rm 1)}
 $
\int_s^t(af +b g )dW=a\int_s^tf dW+b\int_s^t g
dW, \dbP\mbox{-}\as\!\!;
 $

{\rm 2)}
$\dbE\big(\int_s^t f  dW\;\big|\;
\cF_s\big)=0, \dbP\mbox{-}\as\!\!;
 $

{\rm 3)}
 $\dbE\big(\lan\int_s^tf dW,\;\int_s^t g
dW{\ran}_H\,\big|\;
\cF_s\big)=\dbE\big(\int_s^t
\lan f (r,\cd),g
(r,\cd)
{\ran}_Hdr \;\big|\;\cF_s\big),\; \dbP\mbox{-}\as\!\!$;

{\rm 4)} The stochastic process
$ \{\int_0^t f(s)dW(s)\}_{t\in [0,T]} $
is a martingale.
\end{theorem}

For any $p\in (0,\infty)$, denote by
$L_\dbF^{p,loc}(0,T;H)$ the
set of $\mathbf{F}$-adapted stochastic processes
$f(\cd)$ satisfying only
$\int_0^T|f(t)|_{H}^pdt
<\infty$, $\dbP$-a.s.
One can define the It\^o integral $\int_0^t\Phi  dW$ for $\Phi\in L_\dbF^{2,loc}(0,T;H)$ (See \cite{LZ3, YZ} for more details).

\begin{definition} An $H$-valued, $\mathbf{F}$-adapted
process $X(\cd)$ is called an It\^o process if there exist two $H$-valued stochastic processes $\phi(\cd)\in L_\dbF^{1,loc}(0,T;H)$
and $\Phi(\cd)\in L_\dbF^{2,loc}(0,T;H)$ such that
\bel{5.888}
X(t)=X(0)+\int_0^t\phi(s)ds+\int_0^t\Phi(s)dW(s),\quad \dbP\hbox{-a.s.}, \;\;\forall\; t\in[0,T].
\ee
\end{definition}

The following fundamental result is known as {\it It\^o's formula}.

\begin{theorem} Let
$X(\cd)$ be given by (\ref{5.888}). Let $F:
[0,T]\t H\to \dbR$ be a function such that
its partial derivatives $F_t$, $F_x$ and
$F_{xx}$ are uniformly continuous on any
bounded subset of $[0,T]\t H$. Then,
\begin{equation}\label{7e2}
\begin{array}{ll}
F(t,X(t))- F(0,X(0))\\
\ns\ds=\int_0^t
F_x(s,X(s)) \Phi(s)dW(s)+\int_0^t\Big[F_t(s,X(s))+\big\langle
F_x(s,X(s)),\phi(s)\big\rangle_H\\
\ns \ds\q+\frac{1}{2}\langle
F_{xx}(s,X(s))\Phi(s),\Phi(s)\rangle_{H}\Big]ds,\quad \dbP\hbox{-a.s.}, \;\;\forall\; t\in[0,T].
\end{array}
\end{equation}
\end{theorem}

The following deep result, known as the {\it Burkholder-Davis-Gundy
inequality}, links It\^o's integral to the Lebesgue/Bochner integral.

\begin{theorem}\label{BDG} For any $p> 0$, there exists a constant
$\cC_p>0$ such that for any $T>0$ and $f\in
L_\dbF^{p,loc}(0,T;H)$,
\begin{equation}\label{BDGQ-eq2}
\ba{ll}
\ds\frac{1}{\cC_p}
\dbE\(\int_0^T|f(s)|_H^2
ds\)^{\frac{p}{2}}\leq
\dbE\(\sup_{t\in [0,T]}\Big|\int_0^t
f(s)dW(s)\Big|_H^p\)\leq \cC_p
\dbE\(\int_0^T|f(s)|_H^2
ds\)^{\frac{p}{2}}.
\ea
\end{equation}
\end{theorem}


\subsection{Stochastic evolution equations}\label{Ch-SEE}

In what follows, we shall always assume that $H$ is a separable Hilbert space, and
$A$ is an unbounded linear operator (with domain
$D(A)$ on $H$), which is the infinitesimal
generator of a $C_0$-semigroup $\{S(t)\}_{t\geq
0}$. Denote by $A^*$ the dual operator of $A$.
Clearly, $D(A)$ is a Hilbert space with the
usual graph norm, and $A^*$ is the infinitesimal
generator of $\{S^*(t)\}_{t\geq 0}$, the dual
$C_0$-semigroup of $\{S(t)\}_{t\geq 0}$.

Let us consider the following stochastic evolution equation:
\begin{equation}\label{c1-system1}
\left\{
\begin{array}{ll}\ds
dX(t) = \big[A X(t) + F(t,X(t))\big]dt +
\wt F(t,X(t))dW(t) &\mbox{ in } (0,T],\\
\ns\ds X(0)=X_0.
\end{array}
\right.
\end{equation}
Here $X_0\in L^p_{\cF_0}(\Omega;H)$ (for some $p\geq 2$),
and $F(\cd,\cd)$ and $\wt F(\cd,\cd)$
are measurable functions from
$[0,T]\times\Omega\times H$ to $H$, satisfying
the following conditions:
\begin{condition}\label{c1-Lip-Con1}
\begin{equation}
\left\{
\begin{array}{ll}\ds
|F(t,y)-F(t,z)|_{H}\leq \cC |y-z|_H,\q\forall \,y,z\in H,\
\ae t\in [0,T],\,\dbP\mbox{-}\as\!\!,\\
\ns\ds |\wt F(t,y)-\wt F(t,z)|_{H}\leq \cC
|y-z|_H,\q\forall \,y,z\in H,\
\ae t\in [0,T],\,\dbP\mbox{-}\as\!\!,\\
\ns\ds F(\cd,0)\in
L^p_\dbF(\Omega;L^1(0,T;H)),\q \wt F(\cd,0)\in
L^p_\dbF(\Omega;L^2(0,T;H)).
\end{array}
\right.
\end{equation}
\end{condition}

First, we give the notion of strong
solution to the equation \eqref{c1-system1}.

\begin{definition}\label{def-strong-sol}

An $H$-valued stochastic process $X(\cd)\in C_\dbF([0,T];L^p(\Omega;H))$ is
called a strong solution to \eqref{c1-system1}
if $X(t,\omega)\in D(A)$ for  $\ae (t,\omega)\in
  [0,T]\times\Omega$, $AX(\cd)\in
  L_\dbF^{1,loc}(0,T;H)$, and for all $t\in [0,T]$,
  $$
   X(t) = X_0 + \int_0^t \big[A X(s)
   + F(s,X(s)) \big] ds  + \int_0^t \wt F(s,X(s))dW(s),\;\dbP\mbox{-}\as
  $$
\end{definition}

Generally speaking, one needs very strong conditions to
guarantee the existence of a strong solution. Thus,
people introduce two types of ``weak" solutions.

\begin{definition}\label{c1-def-weak}
An $H$-valued stochastic process
$X(\cd)\in
C_\dbF([0,T];L^p(\Omega;H))$ is called a \index{weak solution} weak solution to
\eqref{c1-system1} if for any $t\in [0,T]$ and
$\xi\in D(A^*)$,
 $$
\begin{array}{ll}\ds
\big\langle X(t), \xi \big\rangle_H &\ds=
\big\langle X_0, \xi \big\rangle_H+ \int_0^t
\big(\big\langle X(s), A^*\xi \big\rangle_H +
\big\langle F(s,X(s)), \xi \big\rangle_H \big)
ds\\
\ns&\ds \q  + \int_0^t \big\langle \wt
F(s,X(s)), \xi
\big\rangle_HdW(s),\qquad\dbP\mbox{-}\as
\end{array}
$$
\end{definition}

\begin{definition}\label{c1-def-mild}
An $H$-valued stochastic process
$X(\cd)\in
C_\dbF([0,T];L^p(\Omega;H))$ is called a mild solution to
\eqref{c1-system1} if for any $t\in [0,T]$,
  $$
    X(t) = S(t)X_0 + \int_0^t
  S(t-s)F(s,X(s))ds + \int_0^t
  S(t-s)\wt
  F(s,X(s))dW(s),\q\dbP\mbox{-a.s.}
  $$
\end{definition}

It is easiest to show the well-posedness of
\eqref{c1-system1} in the framework of mild
solution among the above three kinds of solutions.
Indeed, we have the following result.

\begin{theorem}\label{ch-1-well-mild}
Let $p\geq 2$. Then, there is a unique mild
solution $X(\cd)\in C_{\dbF}([0,T];L^p(\Omega;H))$ to \eqref{c1-system1}. Moreover,
\begin{equation}\label{ch-1-well-mild-eq1}
\begin{array}{ll}\ds
|X(\cd)|_{C_\dbF([0,T];L^p(\Omega;H))} \leq
\cC\big(|X_0|_{L^p_{\cF_0}(\Omega;H)} +
|F(\cd,0)|_{L^p_\dbF(\Omega;L^1(0,T;H))} + |\wt
F(\cd,0)|_{L^p_\dbF(\Omega;L^2(0,T;H))}\big).
\end{array}
\end{equation}
\end{theorem}

If  $p>2$ or $\{S(t)\}_{t\geq 0}$ is a
contraction semigroup, then one can get a better
regularity for the mild solution with respect to
time, i.e.,
 $X(\cd,\omega)\in C([0,T];H)$, $\dbP$-a.s.
Here we only consider the latter case.

\begin{theorem}\label{ch-1-well-mild1}
If $A$ generates a contraction semigroup and
$p\ge 1$, then \eqref{c1-system1} admits a unique mild
solution $X(\cd)\in
L^p_\dbF(\Omega;C([0,T];H))$. Moreover,
\begin{equation}\label{ch-1-well-mild1-eq1}
\begin{array}{ll}\ds
|X(\cd)|_{L^p_\dbF(\Omega;C([0,T];H))}\leq
\cC\big(|X_0|_{L^p_{\cF_0}(\Omega;H)} +
|F(\cd,0)|_{L^p_\dbF(\Omega;L^1(0,T;H))}+ |\wt
F(\cd,0)|_{L^p_\dbF(\Omega;L^2(0,T;H))}\big).
\end{array}
\end{equation}
\end{theorem}

The following result indicates the space
smoothing effect of mild solutions to a class of
stochastic evolutions equations, say the
stochastic parabolic equation.

\begin{theorem}\label{ch-1-well-mild2}
Let $p\geq 1$. Assume that $A$ is a self-adjoint, negative
definite (unbounded linear) operator on $H$.
Then, the
equation \eqref{c1-system1} admits a unique mild
solution
$ X(\cd)\in L^p_\dbF(\Omega;C([0,T];H))\cap
L^p_\dbF(\Omega;L^2(0,T;D((-A)^{\frac{1}{2}}))).
$
Moreover,
\begin{equation}\label{ch-1-well-mild2-eq1}
\begin{array}{ll}\ds
|X(\cd)|_{L^p_\dbF(\Omega;C([0,T];H))} + |X(\cd)|_{L^p_\dbF(\Omega;L^2(0,T;D((-A)^{\frac{1}{2}})))} \\
\ns\ds \leq
\cC\big(|X_0|_{L^p_{\cF_0}(\Omega;H)}+
|F(\cd,0)|_{L^p_\dbF(\Omega;L^1(0,T;H))}+
|\wt
F(\cd,0)|_{L^p_\dbF(\Omega;L^2(0,T;H))}\big).
\end{array}
\end{equation}
\end{theorem}

Next result gives the relationship between mild
and weak solutions to \eqref{c1-system1}.

\begin{theorem}\label{ch-1-well-rel1}
Any  weak solution to \eqref{c1-system1} is also a
mild solution and vice versa.
\end{theorem}

Usually, the mild solution does not have enough
regularity. For example, when
establishing the pointwise identity for Carleman
estimate, we need the functions to be second order
differentiable in the sense of weak derivative
with respect to the spatial  variable. Nevertheless,
these problems can be solved by the following
strategy:

\begin{enumerate}
  \item Introduce some approximating equations with
strong solutions such that the limit of these
strong solutions is the mild or weak solution of
the original equation.
  \item  Obtain the desired
 properties for these strong solutions.
  \item Utilize the density argument to establish the desired properties for the
mild/weak solutions.
\end{enumerate}

There are many methods to implement the above
three steps in the setting of deterministic
partial differential equations. Roughly
speaking, any of these methods, which does not
destroy the adaptedness of the solution, can be
applied to stochastic partial differential
equations.  Here we only present one
approach.  Introduce an approximating system of
\eqref{c1-system1} as follows:
\begin{equation}\label{c1-system3}
\left\{
\begin{array}{ll}\ds
dX^\l(t) = AX^\l(t)dt + R(\l) F(t,X^\l(t))dt +
R(\l)\wt F(t,X^\l(t))dW(t) &\mbox{ in } (0,T],\\
\ns\ds X^\l(0)=R(\l)X_0\in D(A).
\end{array}
\right.
\end{equation}
Here   $\l\in \rho(A)$, the resolvent set
of $A$, and $R(\l)\=\l(\l I-A)^{-1}$
with $I$ being the identity operator on $H$.

\begin{theorem}\label{ch-2-app1}
For
each $X_0\in L^p_{\cF_0}(\Omega;H)$ with $p\geq 2$ and
$\l\in\rho(A)$, the equation \eqref{c1-system3}
admits a unique strong solution $X^\l(\cd)\in
C_\dbF([0,T];L^p(\Omega;H))$.
Moreover, as $\l\to\infty$, the solution
$X^{\l}(\cd)$ converges to $X(\cd)$ in
$C_\dbF([0,T];L^p(\Omega;H))$, where $X(\cd)$
solves \eqref{c1-system1} in the sense of the
mild solution.
\end{theorem}
%


\subsection{Backward stochastic  evolution
equations}


Backward stochastic differential
equations and more generally, backward stochastic  evolution
equations are by-products in the study of stochastic control theory, both of which have independent interest and been applied in other places.

Let us consider the following $H$-valued, backward stochastic  evolution
equation
\begin{equation}\label{c1-system4}
\left\{
\begin{array}{ll}\ds
dy(t) = -\big[Ay(t) + F(t,y(t),Y(t))\big]dt -
Y(t) dW(t) &\mbox{
in } [0,T),\\
\ns\ds y(T)=\xi.
\end{array}
\right.
\end{equation}
Here $\xi\in L^p_{\cF_T}(\Omega;H)$ (for some $p\geq 1$),
$F:[0,T]\times\Omega\times H\times
H\to H$ is a  measurable functionm satisfying that
\begin{equation}\label{c3-Lip}
\left\{
\begin{array}{ll}\ds
F(\cd,0,0)\in L^p_\dbF(\Om;L^1(0,T;H)),
\\
\ns\ds |F(t,y_1,z_1)-F(t,y_2,z_2)|_{H} \leq
\cC(|y_1-y_2|_H + |z_1-z_2|_{H}),
\\ \ns\ds \qq\qq\qq\qq\forall\;
y_1,y_2,z_1,z_2\in {H}, \ \ \ae t\in
[0,T],\ \ \dbP\mbox{-}\as
\end{array}
\right.
\end{equation}

Similarly to the case of stochastic  evolution
equations, one introduces below notions of strong, weak and mild
solutions to the equation \eqref{c1-system4}.

\begin{definition}\label{def-strong-sol-b}
A stochastic process $(y(\cd),Y(\cd))\in L^p_\dbF(\Omega;C([0,T];H))\times L^p_\dbF(\Om;L^2(0,T;H))$
is called a strong solution to
\eqref{c1-system4} if $y(t)\in D(A)$ for a.e. $(t,\omega)\in
  [0,T]\times\Omega$, $Ay(\cd)\in
  L^{1,loc}_\dbF(0,T;H)$, and for all $t\in [0,T]$,
  $$
   y(t)  \ds= \xi + \int_t^T \big[A y(s)
   +F(s,y(s),Y(s)) \big] ds   + \int_t^T
   Y(s) dW(s),\q\dbP\mbox{-}\as
  $$
\end{definition}

\begin{definition}\label{def-weak-sol-b}
A stochastic process $(y(\cd),Y(\cd))\in L^p_\dbF(\Omega;C([0,T];H))\times L^p_\dbF(\Om;L^2(0,T;H))$
is called a weak solution to \eqref{c1-system4}
if for any $t\in [0,T]$ and $\eta\in
D(A^*)$,
  $$
  \begin{array}{ll}\ds
   {\lan y(t),\eta\ran}_H&\ds= {\lan \xi,\eta\ran}_H  +
   \int_t^T
   {\lan y(s),A^*\eta\ran}_H ds   \\
   \ns&\ds \q - \int_t^T
  {\lan F(s,y(s),Y(s)),\eta\ran}_H ds- \int_t^T {\lan Y(s),\eta\ran}_HdW(s),\quad\dbP\mbox{-}\as
   \end{array}
  $$
  \end{definition}

\begin{definition}\label{def-mild-sol-b}
A stochastic process $(y(\cd),Y(\cd))\in L^p_\dbF(\Omega;C([0,T];H))\times L^p_\dbF(\Om;L^2(0,T;H))$
is called a mild solution to \eqref{c1-system4}
if for any $t\in [0,T]$,
  $$
  \begin{array}{ll}\ds
   y(t) =  S(T-t)\xi +\int_t^T S(s-t)F(s,y(s),Y(s))ds + \int_t^T
   S(s-t) Y(s) dW(s),\quad\dbP\mbox{-}\as
   \end{array}
  $$
\end{definition}

Similar to Theorem \ref{ch-1-well-mild} (but here one needs that the filtration $\mathbf{F}$ is natural), one can get the
well-posedness of \eqref{c1-system4} in the
sense of mild solution.
\begin{theorem}\label{ch-1-well-bmild}
Assume that $\mathbf{F}$ is the
natural filtration generated by $W(\cd)$. Then, for any $p\ge 1$ and $\xi\in L^p_{\cF_T}(\Omega;H)$, the
equation \eqref{c1-system4} admits a unique mild
solution $(y(\cd),Y(\cd))\in L^p_\dbF(\Omega;C([0,T];H))$ $\times L^p_\dbF(\Om;L^2(0,T;H))$ satisfying that
\begin{equation}\label{ch-1-well-bmild-eq1}
|(y,Y)|_{L^p_\dbF(\Omega;C([0,T];H))\times
L^p_\dbF(\Om;L^2(0,T;H))}   \leq
\cC\big(|\xi|_{L^p_{\cF_T}(\Omega;H)} +
|F(\cd,0,0)|_{L^p_\dbF(0,T;H)} \big).
\end{equation}
\end{theorem}

Also, similar to Theorem \ref{ch-1-well-rel1}, we have the following relationship between the weak
and mild solutions to \eqref{c1-system4}.

\begin{theorem}\label{ch-1-well-brel1}
A stochastic process $(y,Y)$ is a weak solution to \eqref{c1-system4}
if and only if it is a mild solution to the same equation.
\end{theorem}

Similarly to Theorem \ref{ch-1-well-mild2}, the following result describes the the smoothing effect of mild solutions to a class of backward
stochastic evolution equations.

\begin{theorem}\label{ch-1-bwell-mild2}
Let $\mathbf{F}$ be the
natural filtration generated by $W(\cd)$, $F(\cd,
0,0)\in L^1_\dbF(0,T;$ $L^2(\Omega;H))$, and $A$ be
a self-adjoint, negative definite (unbounded
linear) operator on $H$. Then, for any $ \xi\in
L^2_{\cF_T}(\Omega;H)$, the equation
\eqref{c1-system4} admits a unique mild solution
$(y(\cdot),Y(\cdot))\in
\big(L^2_\dbF(\Omega;C([0,T];$ $H))\cap
L^2_\dbF(0,T;D((-A)^{\frac{1}{2}}))\big)\times
L^2_\dbF(0,T;H)$.
Moreover,
\begin{equation}\label{ch-1-well-bmild2-eq1}
\begin{array}{ll}\ds
|y(\cd)|_{L^2_\dbF(\Omega;C([0,T];H))} + |y(\cd)|_{L^2_\dbF(0,T;D((-A)^{\frac{1}{2}}))} + |Y(\cd)|_{L^2_\dbF(0,T;H)} \\
\ns\ds \leq
\cC\big(|\xi|_{L^2_{\cF_T}(\Omega;H)} +
|F(\cd,0,0)|_{L^1_\dbF(0,T;L^2(\Omega;H))}\big).
\end{array}
\end{equation}
\end{theorem}

Similarly to \eqref{c1-system3}, we introduce an
approximating equation of \eqref{c1-system4} as
follows:
\begin{equation}\label{c1-system5}
\left\{
\begin{array}{ll}\ds
dy^\l(t) = -\big[Ay^\l(t)+ R(\l)
F(t,y^\l(t),Y(t))\big]dt -
R(\l)Y^\l(t)dW(t) &\mbox{ in } (0,T],\\
\ns\ds y^\l(T)=R(\l)\xi \in D(A).
\end{array}
\right.
\end{equation}
Similarly to Theorem \ref{ch-2-app1}, we have the following result.

\begin{theorem}\label{app th1}
Assume that $\mathbf{F}$ is the
natural filtration generated by $W(\cd)$, and $F(\cd,
0,0)\in L^1_\dbF(0,T;$ $L^2(\Omega;H))$.
Then, for each $\xi\in L^2_{\cF_T}(\Omega;H)$
and $\l\in\rho(A)$, the equation \eqref{c1-system5} admits a
unique strong solution $(y^\l(\cd),Y^\l(\cd))\in
L^2_\dbF(\Omega;C([0,T];D(A)))\times
L^2_\dbF(0,T;H)$. Moreover, as
$\l\to\infty$, $(y^\l(\cd),Y^\l(\cd))$ converges
to $(y(\cd),Y(\cd))$ (in $
L^2_\dbF(\Omega;C([0,T];H))\times
L^2_\dbF(0,T;H)$), the mild solution to
\eqref{c1-system4}.
\end{theorem}

Note that, in Theorems \ref{ch-1-well-bmild} and \ref{ch-1-bwell-mild2}--\ref{app th1}, we need the filtration $\mathbf{F}$ to be natural. For the general filtration, as we shall see later, we need to employ the stochastic transposition method (developed in \cite{LZ, LZ1, LZ2}) to show the well-posedness of the equation \eqref{c1-system4}.


\section{Controllability of stochastic (ordinary) differential equations}\label{sec ex finite}

In this section, we assume  $\mathbf{F}$ the
natural filtration generated by $W(\cd)$.

We begin with the following controlled system governed by a deterministic linear
ordinary differential equation:
 \bel{ols1}
 \left\{
 \ba{ll}
 \ds\frac{dy(t)}{dt} =Ay(t)+Bu(t),\qq t>0,\\
  \ns
 y(0)=y_0.
 \ea\right.
 \ee
In (\ref{ols1}), $A\in \mathbb{R}^{n\times n}$, $B\in
\mathbb{R}^{n\times m}$ ($n, m\in\mathbb{N}$), $y(\cd)$ is the {\it
state variable}, $u(\cd)$ is the {\it control variable},
$\mathbb{R}^n$ and $\mathbb{R}^m$ are respectively the {\it state and
control spaces}.

\begin{definition}\label{20191002def1}
The system (\ref{ols1}) is called
exactly controllable at time $T$ if for any
$y_0,y_T\in\mathbb{R}^n$, there
is a control $u(\cd)\in L^1(0,T;\mathbb{R}^m)$ such that the
solution $y(\cd)$ to (\ref{ols1}) satisfies $ y(T)=y_T$.
\end{definition}

One has the following result:

\begin{theorem}\label{20191002th1}
The system (\ref{ols1}) is exactly
controllable at time $T$ if and only if the  Kalman rank condition holds
 $$ \hbox{\rm
rank$\,$}[B,AB, \cdots, A^{n-1} B] = n.
 $$
 \end{theorem}

Write
 $
  G_T=\int_0^Te^{At}BB^\top e^{A^\top t}dt$. Further, one can show the following result:

\begin{theorem}\label{20161002t1}
If the system (\ref{ols1}) is exactly
controllable at time $T$, then $\det G_T\not=0$. Moreover, for any $y_0,y_T\in \dbR^n$, the control
  $$
  u^*(t)=-B^\top e^{A^\top (T-t)}G_T^{-1}(e^{AT}y_0-y_T)
  $$
 transfers $y_0$ to $y_T$ at time $T$.
 \end{theorem}

\begin{remark}From Theorem \ref{20161002t1}, it is easy to see that, if (\ref{ols1}) is exactly
controllable at time $T$ (by means of $L^1$-(in time) controls), then the same controllability can be achieved by
 using analytic-(in time) controls. Actually the same can be said for the case that the control class $L^1(0,T;\mathbb{R}^m)$ in Definition \ref{20191002def1} is replaced by $L^p(0,T;\mathbb{R}^m)$ for any $p\in [1,\infty]$. However, we shall see a completely
 different phenomenon evev in the simplest stochastic situation.
\end{remark}

Now, let us consider the
following controlled system governed by a stochastic linear
ordinary differential equation:
\begin{equation}\label{sto ode sy2}
\left\{
\begin{array}{lll}\ds
dy = (Ay+Bu)dt + (Cy + Du)dW(t) &\mbox{ in } [0,T],\\
\ns\ds y(0)=y_0,
\end{array}
\right.
\end{equation}
where $C\in \dbR^{n\t n}$ and $D\in \dbR^{n\t
m}$,  $u(\cd)$ is the
control (valued in $\dbR^m$) and $x(\cd)$ is the
state (valued in $\dbR^n$).

\begin{definition}\label{def sto ode1}
The system \eqref{sto ode sy2} is called exactly
controllable (at time $T$) if for any $y_0\in
\dbR^n$ and $y_T\in L^2_{\cF_T}(\Omega;\dbR^n)$,
there exists a control $ u(\cd) \in
L^2_\dbF(0,T;\dbR^m)$ such that  the
corresponding solution $y(\cd)\in
L^2_{\dbF}(\Omega;C([0,T];\dbR^n))$ to \eqref{sto
ode sy2} satisfies that $y(T)=y_T$.
\end{definition}

Define a (deterministic) function $\eta(\cd)$ on $[0,T]$ by
\begin{equation}\label{6.18-eq13}
\eta(t) = \left\{
\begin{array}{ll}
\ds 1, &\mbox{ for } t\in \ds\Big[\Big(1-\frac{1}{2^{2i}}\Big)T, \Big(1-\frac{1}{2^{2i+1}}\Big)T \Big), \q i=0,1,2,\cds,\\
\ns\ds  -1,  &\mbox{ otherwise}.
\end{array}
\right.
\end{equation}
One can show that {\rm (\cite{Peng})}  there exists a constant
$\b>0$ such that
\begin{equation}\label{sto ode f1}
\int_t^T |\eta(s)-c|^2 ds \geq 4\b (T-t),
\q\mbox{ for any } (c,t)\in \dbR\t [0,T].
\end{equation}

One has the following result,
which provides a necessary condition for the
exact controllability of \eqref{sto ode sy2}.
\begin{proposition}\label{sto ode nece}
{\rm (\cite{Peng})} If the system \eqref{sto ode sy2} is exactly
controllable, then $\rank D=n$.
\end{proposition}

{\it Proof}\,: We use the contradiction
argument. Assume that the system \eqref{sto ode
sy2} was exactly controllable for some matrix
$D$ with $\rank D< n$. Then, we
would find a vector $v\in \dbR^n$ with
$|v|_{\dbR^n}=1$ such that $v\cd D = 0$.

Let  $y_T = \int_0^T \eta(t)dW(t) v$ (recall
(\ref{6.18-eq13}) for $\eta(\cd)$). Since
\eqref{sto ode sy2} was exactly controllable,
there would exist a control $u\in
L^2_{\dbF}(0,T;\dbR^m)$ such that
$$
y_T = y_0 + \int_0^T \big[Ay(t)+Bu(t)\big]dt +
\int_0^T\big[  Cy(t) + Du(t) \big]dW(t),
$$
which implies that
\begin{equation}\label{sto ode eq1}
\int_0^T \eta(t)dW(t) = v\cd y_0 + \int_0^T
v\cd\big[Ay(t)+Bu(t)\big] dt + \int_0^T v\cd
Cy(t) dW(t).
\end{equation}
Hence,
$$
\int_0^T \big[\eta(t) - v \cd Cy(t)\big]dW(t) =
v\cd y_0 + \int_0^T v\cd\big[Ay(t)+Bu(t)\big]dt.
$$
Therefore,
$$
\begin{array}{ll}\ds
\int_0^t \big[\eta(s) - v \cd Cy(s)\big]dW(s)\\
\ns\ds= v\cd y_0 + \int_0^t
v\cd\big[Ay(s)+Bu(s)\big]ds+ \dbE\Big(\int_t^T
v\cd\big[Ay(s)+Bu(s)\big]ds \Big| \cF_t \Big).
\end{array}
$$
This gives that
$$
\begin{array}{ll}\ds
\int_t^T \big[\eta(s) - v \cd Cy(s)\big]dW(s)\\
\ns\ds = \int_t^T v\cd\big[Ay(s)+Bu(s)\big]ds-
\dbE\Big(\int_t^T v\cd\big[Ay(s)+Bu(s)\big]ds
\Big| \cF_t \Big),
\end{array}
$$
which implies that
\begin{equation}\label{sto ode eq2}
\begin{array}{lll}\ds
\dbE\int_t^T \big|\eta(s) - v \cd Cy(s)\big|^2ds\\
\ns\ds = \dbE\Big[\int_t^T v\cd\big[Ay(s)+Bu(s)\big]ds- \dbE\Big(\int_t^T v\cd\big[Ay(s)+Bu(s)\big]ds \Big| \cF_t \Big)\Big]^2\\
\ns\ds \leq \dbE\Big[\int_t^T v\cd\big[Ay(s)+Bu(s)\big]ds\Big]^2\\
\ns\ds \leq (T-t)\int_t^T
\big|v\cd\big[Ay(s)+Bu(s)\big]\big|^2ds.
\end{array}
\end{equation}
On the other hand, by the inequality \eqref{sto
ode f1}, we have that
\begin{equation}\label{sto ode eq3}
\begin{array}{ll}\ds
\dbE\int_t^T \big|\eta(s) - v \cd Cy(s)\big|^2ds\\
\ns\ds \geq \frac{1}{2} \dbE\int_t^T \big|\eta(s) - v \cd Cy(T)\big|^2ds - \dbE\int_t^T \big|v \cd Cy(T) - v \cd Cy(s)\big|^2ds\\
\ns\ds \geq 2\b (T-t) - \dbE\int_t^T \big|v \cd
Cy(T) - v \cd Cy(s)\big|^2ds.
\end{array}
\end{equation}
By virtue of that $y(\cd)\in
L^2_{\dbF}(\Omega;C([0,T];\dbR^n))$,  there is a
$\tilde t\in [0,T)$ such that
$$
\dbE\big|v\cd Cy(T) - v\cd Cy(s)\big| \leq \b,
\q \mbox{ for all } s\in [\tilde t,T).
$$
This, together with \eqref{sto ode eq3}  implies
that
\begin{equation}\label{sto ode eq4}
\dbE\int_t^T \big|\eta(s) - v \cd Cy(s)\big|^2ds
\geq \b (T-t), \mbox{ for all } t\in [\tilde
t,T).
\end{equation}
From \eqref{sto ode eq2} and  \eqref{sto ode
eq4}, we have that
$$
\b \leq \int_t^T
\big|\,v\cd\big[Ay(s)+Bu(s)\big]\big|^2ds,
\mbox{ for all } t\in [\tilde t,T),
$$
which leads to a contradiction.
\endpf

%
%

%
\begin{proposition}\label{sto ode nece1}
If the system \eqref{sto ode sy2} is exactly
controllable at time $T$, then $(A,B)$ fulfills the Kalman
rank condition.
\end{proposition}

{\it Proof}\,: Let $\tilde y = \mE y$, where $y$
is a solution to \eqref{sto ode sy2} with some
$y_0\in \dbR^n$ and $u(\cd)\in
L^2_{\dbF}(0,T;\dbR^m)$. Then $\tilde y$ solves
\begin{equation}\label{4.11-eq2}
\left\{
\begin{array}{lll}\ds
\frac{d\tilde y}{dt} = A\tilde y+B\mE u &\mbox{ in } [0,T],\\
\ns\ds \tilde y(0)=y_0.
\end{array}
\right.
\end{equation}
Since \eqref{sto ode sy2} is exactly
controllable, we see that \eqref{4.11-eq2} is
exactly controllable. Hence, $(A,B)$
fulfills the Kalman rank condition.
\endpf

\ms

By means of Propositions \ref{sto ode nece}--\ref{sto ode nece1},  it follows that we should
assume that $\rank D=n$ and $(A,B)$ fulfills the
Kalman rank condition if we expect the exact
controllability of the system \eqref{sto ode
sy2} in the sense of Definition \ref{def sto ode1}.

Since  $\rank D=n$, it is easy to see that
$n\leq m$, and we can find two matrices $K_1\in
\dbR^{m\t m}$ and $K_2\in \dbR^{m\t n}$ such
that $DK_1 = (I_{n},0)$ and that $DK_2 = -C$.
Introducing a simple linear transformation
$$
u = K_1\left(\begin{array}{ll}\ds v_2\\ \ns\ds v_1
\end{array}\right)+K_2y,
$$
where $v_1\in
L^2_{\dbF}(0,T;\dbR^{m-n})$ and $v_2\in L^2_{\dbF}(0,T;\dbR^n)$ , we see that the
system \eqref{sto ode sy2} is reduced to the
following system
\begin{equation}\label{sto ode sy3}
\left\{
\begin{array}{ll}\ds
dy = (A_1 y + A_2 v_2 + B_1 v_1) dt + v_2dW(t) &\mbox{ in } [0,T],\\
\ns\ds y(0)=y_0,
\end{array}
\right.
\end{equation}
where
$$
A_1 = A + BK_2,  \; A_2 \in \dbR^{n\t n},\; B_1
\in \dbR^{n\t (m-n)} \mbox{ and } A_2 v_2 +
B_1 v_1  = BK_1\left(\begin{array}{ll}\ds v_2\\
\ns\ds v_1 \end{array}\right).
$$

In order to deal with the exact controllability problem for
\eqref{sto ode sy3}, we consider the following controlled backward stochastic differential system:
\begin{equation}\label{sto ode sy4}
\left\{
\begin{array}{ll}\ds
dy = (A_1 y + A_2 Y + B_1 v) dt + YdW(t) &\mbox{ in } [0,T],\\
\ns\ds y(T)=y_T,
\end{array}
\right.
\end{equation}
where $y_T\in L^2_{\cF_T}(\Omega;\dbR^n)$, $v\in
L^2_{\dbF}(0,T;\dbR^{m-n})$ is the control variable.

\begin{definition}
The system \eqref{sto ode sy4} is called
exactly controllable (at time $0$) if for any $y_T\in
L^2_{\cF_T}(\Omega;$ $\dbR^n)$ and $y_0\in \dbR^n$,
there is a control $v\in
L^2_{\dbF}(0,T;\dbR^{n\t (m-n)})$ such that the corresponding
solution $(y(\cd), Y(\cd))\in
L^2_{\dbF}(\Omega;C([0,T];\dbR^n))\times L^2_{\dbF}(0,T;\dbR^n)$ to  \eqref{sto ode sy4} satisfies $y(0)=y_0$.
\end{definition}

It is easy to show the following result:

\begin{proposition}
The system \eqref{sto ode sy3} is
exactly controllable at time $T$ if and only if the system \eqref{sto ode sy4} is
exactly controllable at time $0$.
\end{proposition}

The dual equation of the system \eqref{sto ode sy4} is the
following  stochastic ordinary
differential equation:
\begin{equation}\label{sto ode sy5}
\left\{
\begin{array}{ll}\ds
dz = -A_1^\top z dt - A_2^\top z dW(t) &\mbox{ in } [0,T],\\
\ns\ds z(0) = z_0\in \dbR^n.
\end{array}
\right.
\end{equation}

Similar to Theorem \ref{20191002th1}, one can show the following result:

\begin{theorem}\label{20161003e7}
The following statements are equivalent:

{\rm 1)} The system \eqref{sto ode sy4} is
exactly controllable at time $0$;

{\rm 2)} All solutions to \eqref{sto ode sy5} satisfy the following observability estimate:
 \bel{20161003e23}
 |z_0|^2\le \cC\mE\int_0^T|B_1^\top z(t)|^2dt,\qq \forall\;z_0\in \dbR^n;
 \ee

 {\rm 3)} Solutions to \eqref{sto ode sy5} enjoy the following observability:
  \bel{20161003e3}
 B_1^\top z(\cd)\equiv 0 \hbox{ in }(0,T),\ \as \Rightarrow z_0=0;
 \ee

 {\rm 4)} The following rank condition holds:
  \begin{equation}\label{rank condition}
\rank\,[B_1,\, A_1 B_1,\,  A_2 B_1,\, A_1^2B_1,\,
A_1 A_2 B_1,\, A_2^2B,\, A_2 A_1B_1,\cdots]=n.
\end{equation}
\end{theorem}

{\it Proof}\,: By means of the classical duality argument, it is easy to show that ``1)$\Longleftrightarrow$2)". The proof of ``2)$\Longleftrightarrow$3)" is easy.

\ms

``4)$\Longrightarrow$3)". We use an idea from the proof of \cite[Theorem 3.2]{Peng}. Let us assume that  $B_1^\top z(\cd)\equiv 0 \hbox{ in }(0,T),\ \as$ for some $z_0\in \dbR^n$.
Then,
$$
B_1^\top z(t) = B_1^\top z_0+ \int_0^t B_1^\top
A_1^\top z(s)ds + \int_0^t B_1^\top A_2^\top
z(s)dW(s)=0,\qq \forall\;t\in (0,T).
$$
Therefore, we have that
\begin{equation}\label{rank1}
B_1^\top z_0 = 0,\q B_1^\top A_1^\top z\equiv 0,\q
B_1^\top A_2^\top z\equiv 0.
\end{equation}
Hence  $B_1^\top A_1^\top z_0 =
B_1^\top A_2^\top z_0=0$.

Noticing that $z(\cd)$  solves \eqref{sto ode sy5}, we have that
$$
z(t) = z_0 + \int_0^t A_1^\top z(s)ds + \int_0^t
A_2^\top z(s)dW(s).
$$
This together with \eqref{rank1} implies that
$$
B_1^\top A_1^\top z = B_1^\top A_1^\top z_0 +
\int_0^t B_1^\top A_1^\top A_1^\top z(s)ds +
\int_0^t B_1^\top A_1^\top A_2^\top z(s)dW(s)=0,
$$
and
$$
B_1^\top A_2^\top z = B_1^\top A_2^\top z_0 +
\int_0^t B_1^\top A_2^\top A_1^\top z(s)ds +
\int_0^t B_1^\top A_2^\top A_2^\top z(s)dW(s)=0,
$$
which are equivalent to
$$
B_1^\top A_1^\top A_1^\top z\equiv B_1^\top
A_1^\top A_2^\top z\equiv B_1^\top A_2^\top
A_1^\top z\equiv B_1^\top A_2^\top A_2^\top
z\equiv 0,
$$
and implies that $B_1^\top A_1^\top A_2^\top z_0=B_1^\top A_1^\top A_1^\top z_0=
B_1^\top A_2^\top A_1^\top z_0=B_1^\top A_2^\top A_2^\top z_0= 0$.

Utilizing the above argument, by induction, we
can conclude that
\bel{20161003e5}
z_0^\top [B_1,\, A_1 B_1,\,  A_2 B_1,\, A_1^2B_1,\,
A_1 A_2 B_1,\, A_2^2B,\, A_2 A_1B_1,\cdots]=0.
\ee
By \eqref{rank condition} and \eqref{20161003e5}, it follows that $z_0 = 0$.

\ms

``3)$\Longrightarrow$4)". We use the contradiction argument. Assume that \eqref{rank condition} was false. Then, we could find a nonzero $z_0\in\dbR^n$ satisfying \eqref{20161003e5}. For this $z_0$, denote by $z(\cd)$  the corresponding solution to \eqref{sto ode sy5}. Clearly, $z(\cd)$ can be approximated (in $L^2_{\dbF}(\Omega;C([0,T];\dbR^n))$) by the Picard sequence $\{z_k(\cd)\}_{k=0}^\infty$ defined  as follows
 \bel{20161003e70}
 \left\{
 \ba{ll}
 z_0(\cd)=z_0,\\
 \ns
 \ds
 z_k(\cd)=z_0+\int_0^\cd A_1^\top z_{k-1}(s)ds + \int_0^\cd
A_2^\top z_{k-1}(s)dW(s),\q k\in\dbN.
 \ea
 \right.
 \ee
By \eqref{20161003e5} and \eqref{20161003e70}, via a direct computation, one can show that
 \bel{20161003e6}
 B_1^\top z_k(\cd)=0,\qq k=0,1,2,\cdots.
 \ee
By \eqref{20161003e6}, we deduce that $B_1^\top z(\cd)\equiv 0 \hbox{ in }(0,T)$. Hence, by \eqref{20161003e3}, it follows that $z_0=0$, which is a contradiction.
\endpf

\ms

As a consequence of Theorem \ref{20161003e7}, we have the following characterization for the exact controllability of \eqref{sto ode sy3} (and hence also for that of \eqref{sto ode sy2}).

\begin{corollary}\label{c2-exth} {\rm (\cite{Peng})}  The system \eqref{sto ode sy3} is
exactly controllable at time $T$ if and only if  the rank condition \eqref{rank condition} holds.
\end{corollary}

In the above, we introduce two (different)
controls $v_1$ and $v_2$ in the system
\eqref{sto ode sy3}, and both $v_1$ and $v_2$
are $L^2$-(in time). Is it possible to introduce
only one control or to use other class of
controls?

We consider the simplest one-dimensional
controlled ``stochastic" differential equation
as follows
\bel{1.23}
\left\{\ba{ll}
dy(t)=u(t)dt,\\
\ns
y(0)=y_0.
\ea\right.\ee
We say that the system
(\ref{1.23}) is exactly controllable if for any
$y_0\in\dbR$ and $y_T\in L^2_{{\cal F}_T}(\Om)$, there exists a
control $u(\cd)\in L^1_\dbF (0,T;L^2(\Om))$ such that the
corresponding solution $y(\cd)$ satisfies $y(T)=y_T$.

It is showed in \cite{LYZ} that the
system (\ref{1.23}) is exactly controllable at any time $T>0$ (by
means of $L^1_\dbF (0,T;L^2(\Om))$-controls).

On the other hand, surprisingly, in virtue of
Proposition \ref{sto ode nece}, the system
(\ref{1.23}) is NOT exactly controllable if one
is confined to use admissible controls $u(\cd)$
in $L^2_\dbF (0,T;L^2(\Om))$! Further, the
authors in \cite{LYZ} showed that the system
(\ref{1.23}) is NOT exactly controllable, either
provided that one uses admissible controls
$u(\cd)$ in $L^p_\dbF(0,T;L^2(\Om))$ for any
$p\in(1,\infty]$.

To the best of our knowledge, unlike the
deterministic case, there exists no universally
accepted notion for stochastic controllability
so far. Motivated by the above example, we
introduced a corrected formulation for the exact
controllability of stochastic differential
equations.

\begin{definition}\label{def sto-ode0}
The system \eqref{sto ode sy2} is called
exactly controllable if for any $y_0\in \dbR^n$ and  $y_T\in
L^2_{{\cal F}_T}(\Om;\dbR^n)$, one can find a control $u(\cd)\in
L^1_\dbF (0,T;L^2(\Om;\dbR^m))$ such that $Du(\cd,\omega)\in L^2(0,T;\dbR^n),$ $\ae\omega\in\Omega$, and the corresponding solution $y(\cd)$ to
\eqref{sto ode sy2} satisfies $y(T)=y_T$.
\end{definition}

The above definition seems to be a reasonable
notion for exact controllability of stochastic
differential equations. Nevertheless, a complete
study on this problem is still under
consideration and it does not seem to be easy.

One may think that the requirement of exact
controllability for \eqref{sto ode sy2} is too
strong. How about the null/approximate
controllability? Consider the following two
weaker notions of controllability.

\begin{definition}\label{def sto-ode2}
The system \eqref{sto ode sy2} is called null controllable
(at time $T$) if for any $y_0\in \dbR^n$, there
exists a control $ u(\cd)\in
L^2_\dbF(0,T;\dbR^m)$ such that the
corresponding solution $y(\cd)$ to \eqref{sto
ode sy2} satisfies $y(T)=0$.
\end{definition}
\begin{definition}\label{def sto ode3}
The system \eqref{sto ode sy2} is called approximately
controllable (at time $T$) if for any $y_0\in
\dbR^n$, $y_T\in L^2_{\cF_T}(\Om;\dbR^n)$ and
$\e>0$, there exists a control $u(\cd)\in
L^2_\dbF(0,T;\dbR^m)$ such that the
corresponding solution $y(\cd)$ to \eqref{sto
ode sy2} satisfies
$|y(T)-y_T|_{L^2_{\cF_T}(\Om;\dbR^n)}<\e$.
\end{definition}

We shall show below that there exists no any
rank condition for the null/approximate
controllability of  \eqref{sto ode sy2}. In
fact, if there is a such kind of rank condition,
then it should has the following properties:
\begin{itemize}
  \item It is robust with
respect to perturbations small enough;
  \item The system \eqref{sto ode sy2} is null/approximately
  controllable at time $T$ for any $T>0$.
\end{itemize}
However, as pointed in \cite{LZ4}, such
properties cannot be held. In fact, consider the
following $2$-dimensional stochastic
differential system:
\begin{equation}\label{5.12-eq1}
\begin{cases}\ds
dy_1 = y_2dt + \e y_2dW(t) &\mbox{ in } [0,T],\\
\ns\ds dy_2 = udt  &\mbox{ in } [0,T],\\
\ns\ds y_1(0)=y_{10},\;y_2(0)=y_{20},
\end{cases}
\end{equation}
where $( y_{10},y_{20})\in \dbR^2$, $u(\cd)\in
L^1_\dbF (0,T;L^2(\Om))$ is the control variable, $\e$ is a parameter. Clearly, if $\e=0$, then \eqref{5.12-eq1} is
null controllable. If the above two properties held,
then there would exist an $\e_0>0$ such that for all
$\e\in [-\e_0,\e_0]$ and $T>0$, \eqref{5.12-eq1}
is null controllable. Let us take $y_{10}=0$,
$y_{20}=1$, $\e=\e_0$ and $T=\frac{\e_0^2}{2}$.
Since \eqref{5.12-eq1} is null controllable at
$T=\frac{\e_0^2}{2}$, then
$$
y_1\(\frac{\e_0^2}{2}\)=\int_0^{\frac{\e_0^2}{2}}y_2
dt + \e_0\int_0^{\frac{\e_0^2}{2}}y_2 dW(t)=0.
$$
Thus,
\begin{equation}\label{5.12-eq2}
\mE\Big|\int_0^{\frac{\e_0^2}{2}}y_2 dt\Big|^2 =
\mE\Big|\e_0\int_0^{\frac{\e_0^2}{2}}y_2
dW(t)\Big|^2 =
\e_0^2\int_0^{\frac{\e_0^2}{2}}\mE|y_2|^2dt.
\end{equation}
On the other hand,
\begin{equation}\label{5.12-eq3}
\mE\Big|\int_0^{\frac{\e_0^2}{2}}y_2
dt\Big|^2\leq \mE\Big|
\(\int_0^{\frac{\e_0^2}{2}}1dt\)\(\int_0^{\frac{\e_0^2}{2}}|y_2|^2dt\)\Big|\leq
\frac{\e_0^2}{2}\int_0^{\frac{\e_0^2}{2}}\mE|y_2|^2dt.
\end{equation}
It follows from \eqref{5.12-eq2} and
\eqref{5.12-eq3} that
$\int_0^{\frac{\e_0^2}{2}}\mE|y_2|^2dt=0$, which
contradicts the choice of $y_2(0)$.

Next, we consider the approximate
controllability. For this purpose, we introduce the following
backward stochastic differential equation:
\begin{equation}\label{8.19-eq1}
\begin{cases}\ds
dz_1 = Z_1dW(t) &\mbox{ in } [0,T],\\
\ns\ds dz_2 = -(z_1+\e Z_1)dt + Z_2dW(t)  &\mbox{ in } [0,T],\\
\ns\ds z_1(T)=z_{1T},\;z_2(T)=z_{2T},
\end{cases}
\end{equation}
where $( z_{1T},z_{2T})\in L^2_{\cF_T}(\Omega;\dbR^2)$. By the classical duality argument, it is easy to show that the
approximate controllability of \eqref{5.12-eq1}
is equivalent to the following observability of
\eqref{8.19-eq1}: {\it If $z_2(\cd)=0$, then
$(z_1(\cd),Z_1(\cd),z_2(\cd),Z_2(\cd))=(0,0,0,0)$}.

If $\e=0$ and $z_2(\cd)=0$, then we
$$
-\int_0^t z_1(s) ds + \int_0^tZ_2(s)dW(s)=0,
\qq\mbox{ for all }t\in [0,T].
$$
This, together with the uniqueness of the
decomposition of semimartingale (See \cite[page
358]{Rogers}), implies that
$z_1(\cd)=Z_2(\cd)=0$. Then, by the first
equation in \eqref{8.19-eq1}, we see that
$Z_1(\cd)=0$. Therefore, we conclude that
\eqref{5.12-eq1} is approximately controllable
if $\e=0$.

However, if $\e\neq 0$, then it is
easy to check that
$$
(z_1(t),Z_1(t),z_2(t),Z_2(t))=\(\exp\left\{-\frac{W(t)}{\e}-\frac{t}{2\e^2}\right\},
-\frac{1}{\e}\exp\left\{-\frac{W(t)}{\e}-\frac{t}{2\e^2}\right\},0,0\)
$$
is a solution to \eqref{8.19-eq1} with $(
z_{1T},z_{2T})=\(\exp\left\{-\frac{W(T)}{\e}-\frac{T}{2\e^2}\right\},
-\frac{1}{\e}\exp\left\{-\frac{W(T)}{\e}-\frac{T}{2\e^2}\right\}\)$.
Hence, the above observability of
\eqref{8.19-eq1} does not hold. Therefore,
\eqref{5.12-eq1} is not approximately
controllable whenever $\e \neq0$.

Generally speaking, when $n>1$, the controllability for the linear system \eqref{sto ode sy2} is far from well-understood. Actually, in our opinion, compared to the deterministic case, the controllability/observability for stochastic differential equations
is at its ¡°enfant¡± stage.


\section{Pontryagin-type maximum principle for
controlled stochastic (ordinary) differential
equations}\label{s3}


The first order necessary optimality condition,
i.e., Pontryagin-type maximum principle, for
optimal control problems for stochastic
(ordinary) differential equations is by now
well-understood (at least when there exist no
endpoint constraints). When $\mathbf{F}$ is the
natural filtration,  the general stochastic
maximum principle was established in
\cite{Peng}. In this section, we do not assume
that $\mathbf{F}$ is the natural filtration.
Thus, we cannot use the classical well-posedness
theory of backward stochastic differential
equations. A key point is that we need to use
the stochastic transposition method, developed
in \cite{LZ}.

Let $U$ be a separable metric space with its
metric $\mathbf{d}(\cd,\cd)$. Put
$$\cU[0,T] \triangleq \Big\{u(\cdot):\,[0,T]\to U\;\Big|\; u(\cd) \mbox{ is $\mathbf{F}$-adapted} \Big\}.$$
We assume the following condition.

\ms

\no{\bf (A1)} {\it Suppose that
$a(\cd,\cd,\cd):[0,T]\times \dbR^n\times U\to
\dbR^n$ and $b(\cd,\cd,\cd):[0,T]\times \dbR^n
\times U\to \dbR^n$ are two functions
satisfying: i) For any $(x,u)\in \dbR^n\times
U$, the functions $a(\cd,x,u):[0,T]\to \dbR^n$
and $b(\cd,x,u):[0,T]\to \dbR^n$ are Lebesgue
measurable; ii) For any $(t,x)\in [0,T]\times
\dbR^n$, the functions $a(t,x,\cd):U\to \dbR^n$
and $b(t,x,\cd):U\to \dbR^n$ are continuous; and
iii) There is a constant $\cC_L>0$ such that
 for all $(t,x_1,x_2,u)\in [0,T]\times
\dbR^n\times \dbR^n\times U$,
\begin{equation}\label{9.25-eq2}
\left\{
\begin{array}{ll}\ds
|a(t,x_1,u) - a(t,x_2,u)|_{\dbR^n}+|b(t,x_1,u) -
b(t,x_2,u)|_{\dbR^n} \leq
\cC_L|x_1-x_2|_{\dbR^n},\\
\ns\ds |a(t,0,u)|_{\dbR^n} +|b(t,0,u)|_{\dbR^n}
\leq \cC_L.
\end{array}
\right.
\end{equation}}

Let us consider the following controlled stochastic differential equation:
\begin{equation}\label{9.25-eq1}
\left\{
\begin{array}{lll}\ds
dx =  a(t,x,u) dt + b(t,x,u)dW(t) &\mbox{ in }[0,T],\\
\ns\ds x(0)=x_0,
\end{array}
\right.
\end{equation}
where $u\in \cU[0,T]$ and $x_0\in
L^{p}_{\cF_0}(\Omega;\dbR^n)$ for a given $p\geq
2$.  Under the assumption (A1), it is easy to show that the equation
\eqref{9.25-eq1} is well-posed in the sense of adapted solutions in the space $L^{p}_\dbF(\Omega;C([0,T];\dbR^n))$.

Also, we need the following condition:

\ms

\no{\bf (A2)} {\it Suppose that
$g(\cd,\cd,\cd):[0,T]\times \dbR^n\times U\to
\dbR$ and $h(\cd):\dbR^n\to \dbR$ are two
functions satisfying: i) For any $(x,u)\in
\dbR^n\times U$, the function
$g(\cd,x,u):[0,T]\to \dbR$ is Lebesgue
measurable; ii) For any $(t,x)\in [0,T]\times
\dbR^n$, the function $g(t,x,\cd):U\to \dbR$ is
continuous; and iii) For all $(t,x_1,x_2,u)\in [0,T]\times
\dbR^n\times \dbR^n\times U$,
\begin{equation}\label{9.25-eq3}
\left\{
\begin{array}{ll}\ds
|g(t,x_1,u) - g(t,x_2,u)|_{\dbR^n} +|h(x_1) -
h(x_2)|_{\dbR^n}
\leq \cC_L|x_1-x_2|_{\dbR^n},\\
\ns\ds |g(t,0,u)|_{\dbR^n} +|h(0)|_{\dbR^n} \leq
\cC_L.
\end{array}
\right.
\end{equation}}

\ms

Define a cost functional $\cJ(\cdot)$ (for the
controlled system \eqref{9.25-eq1}) as
follows:
\begin{equation}\label{jk1}
\cJ(u(\cdot))\triangleq \dbE\Big[\int_0^T
g(t,x(t),u(t))dt + h(x(T))\Big],\q\forall\,
u(\cdot)\in \cU[0,T],
\end{equation}
where $x(\cd)$ is the corresponding solution to
\eqref{9.25-eq1}.

Let us consider the following optimal control
problem for the system \eqref{9.25-eq1}:

\ms

\no {\bf Problem (OPF)} {\it Find a $\bar
u(\cdot)\in \cU[0,T]$ such that
\begin{equation}\label{9.25-eq4}
\ds\cJ (\bar u(\cdot)) = \inf_{u(\cdot)\in
\cU[0,T]} \cJ (u(\cdot)).
\end{equation}
Any $\bar u(\cdot)$ satisfying (\ref{9.25-eq4}) is
called an {\it optimal control}. The
corresponding state process $\bar x(\cdot)$ is
called an {\it optimal state (process)}, and
$(\bar x(\cdot),\bar u(\cdot))$ is called an
{\it optimal pair}.}

\ms

Furthermore, we impose the following assumption.

\ms

\no {\bf(A3)} {\it The functions
$a(t,x,u),b(t,x,u),g(t,x,u)$ and $h(x)$ are
$C^2$ in $x$, and for $\f(t,x,u)=b(t,x,u),$ $\si(t,x,u),
f(t,x,u),h(x)$ and any $t\in[0,T]$, $x,\h
x\in\dbR^n$ and $u,\h u\in U$, it holds that
$$
\left\{
\begin{array}{ll}
|\f(t,x,u)-\f(t,\h x,\h u)|\le  \cC_L\big(|x-\h
x|+ \mathbf{d}(u,\h u)\big),\\
\ns\ds |\f(t,0,u)|\leq  \cC_L,\\\ns\ds
|\f_{x}(t,x,u)-\f_{x}(t,\h x,\h u)|\le
\cC_L\big(|x-\h
x|+ \mathbf{d}(u,\h u)\big),\\
\ns\ds|\f_{xx}(t,x,u)-\f_{xx}(t,\h x,\h u)|\le
\cC_L\big(|x-\h x|+ \mathbf{d}(u,\h u)\big).
\end{array}
\right.
$$}

Suppose that $(\bar x(\cd), \bar u(\cd))$ is a
given optimal pair. Similar to the corresponding
deterministic setting, one introduces the
following first order adjoint equation (which is
however a backward stochastic differential equation in the stochastic case):
\begin{equation}\label{3.8}
\left\{\3n \begin{array}{ll}\ds
dy(t)=\!-\big[a_x(t,\bar x(t),\bar u(t))^\top
y(t)\!+\!b_x(t,\bar x(t),\bar u(t))^\top  Y(t)
\!-\!g_x(t,\bar x(t),\bar
u(t))\big]dt\\
\ns\ds \hspace{8.6cm}+ Y(t)dW(t)\qq\mbox{in }
[0,T],\\
\ns\ds y(T)=-h_x(\bar x(T)).
\end{array}
\right.
\end{equation}

Next, to establish the desired maximum principle
for stochastic controlled systems with
control-dependent diffusion and possibly
nonconvex control domains,  one has to introduce
an additional second order adjoint equation as
follows:
\begin{equation}\label{3.9}
\!\!\!\left\{\ba{ll} \!\!
\!dP(t)\!=\!-\Big[a_x(t,\bar x(t),\bar
u(t))^\top P(t)\!+\!P(t) a_x(t,\bar x(t),\bar
u(t))\! +\!b_x(t,\bar x(t),\bar u(t))^\top
P(t)b_x(t,\bar x(t),
\bar u(t))\\
\ns
\qq\qq  +b_x(t,\bar x(t),\bar u(t))^\top
Q(t)\!+\!Q(t) b_x(t,\bar x(t),\bar u(t))
\!+\!\dbH_{xx}(t,\bar x(t),\bar
u(t),y(t),Y(t))\Big]dt\\
\ns\ds \qq\qq+ Q(t) dW(t) \hspace{9.9cm} \mbox{ in } [0,T),\\
\ns
\!\!\!P(T)=-h_{xx}(\bar x(T)).\ea\right.
\end{equation}
In (\ref{3.9}), the {\it Hamiltonian}
$\mathbb{H}(\cd,\cd,\cd,\cd,\cd)$ is defined by
$$
\begin{array}{ll}
\dbH(t,x,u,y_1,y_2)=\lan
y_1,a(t,x,u)\ran_{\dbR^n}+\lan y_2, b(t,x,u)
\ran_{\dbR^n}-g(t,x,u),\\\ns
\hspace{4.5cm}(t,x,u,y_1,y_2)\in[0,T] \times
\dbR^n \times U \times  \dbR^n \times \dbR^n.
\end{array}
$$

Since we do not assume that $\mathbf{F}$ is the
natural filtration, the
equations \eqref{3.8}/\eqref{3.9} may not
have classical adapted solutions. We need to
introduce below the notion of transposition
solutions to the following backward stochastic differential equation:
\begin{equation}\label{9.25-eq5}
\left\{
\begin{array}{lll}
\ds dy (t) = f(t,y(t),Y(t))dt + Y(t) dW(t) & \mbox{ in } [0,T], \\
\ns\ds y(T) = y_T,
\end{array}
\right.
\end{equation}
where $y_T\in L_{\cF_T}^p(\Om;\dbR^n)$,
$f(\cd,\cdot,\cdot)$ satisfies $f(\cd,0,0)\in
L^p_{\dbF}(\Om;L^1(0,T;\dbR^n))$, and
\begin{equation}\label{9.25-eq6}
\begin{array}{ll}\ds
|f(t,p_1,q_1)-f(t,p_2,q_2)|_{\dbR^n}\le
\cC_L(|p_1-p_2|_{\dbR^n}+|q_1-q_2|_{\dbR^n}),\\
\ns\ds \hspace{4.3cm} t\in [0,T],\,
\dbP\hb{-a.s.}, \forall\;p_1,p_2,q_1,q_2\in
\dbR^n.
\end{array}
\end{equation}

In order to define the transposition solution to
\eqref{9.25-eq5}, for any $t\in [0,T]$, we consider
the following linear  stochastic
differential equation
\begin{equation}\label{9.25-eq7}
\left\{
\begin{array}{lll}
\ds dz(\tau)  = u(\tau)d\tau + v(\tau)dW(\tau), & \tau\in (t,T], \\
 \ns\ds z(t) = \eta.
\end{array}
\right.
\end{equation}
For any given $u(\cdot)\in
L^1_{\dbF}(t,T;L^q(\Om;\dbR^n))$, $v(\cdot)\in
L^q_{\dbF}(\Om;L^2(t,T;\dbR^n))$ and $\eta\in
 L^q_{\cF_t}(\Om;\dbR^n)$, the equation \eqref{9.25-eq7} admits a unique adapted solution $z(\cdot)\in
L^{q}_{\dbF}(\Om;C([t,T];\dbR^n))$. Now, if the
equation \eqref{9.25-eq5} admits an adapted
solution $(y(\cdot), Y(\cdot)) \in
L^{p}_{\dbF}(\Om;C([0,T];\dbR^n)) \t
L^{p}_{\dbF}(0,T;L^2(\Om;$ $\dbR^n))$, then,
applying It\^o's formula to $\lan
z(t),y(t){\ran}_{\dbR^n}$, it is easy to check
that
\begin{equation}\label{sol1}
\begin{array}{ll}
\ds\mathbb{E}\lan z(T),y_T{\ran}_{\dbR^n} - \mathbb{E}\lan
\eta,y(t){\ran}_{\dbR^n}\\
\ns \ds = \mathbb{E}\int_t^T \lan
z(\tau),f(\tau,y(\tau),Y(\tau)){\ran}_{\dbR^n} d\tau +
\mathbb{E}\int_t^T \lan u(\tau),y(\tau){\ran}_{\dbR^n}
d\tau + \mathbb{E} \int_t^T\lan v(\tau),
Y(\tau){\ran}_{\dbR^n}d\tau.
\end{array}
\end{equation}
This inspires us to introduce the following new
notion of solution to the equation
\eqref{9.25-eq5}.

\begin{definition}\label{def of solution}
We call $(y(\cdot), Y(\cdot)) \!\in\!
D_{\dbF}([0,T];L^{p}(\Om;\dbR^n)) \t
L^p_{\dbF}(\Om;L^2(0,T;\dbR^n))$ a transposition
solution to \eqref{9.25-eq5} if the identity
(\ref{sol1}) holds for any $t\in [0,T]$,
$u(\cdot)\in L^1_{\dbF}(t,T;L^q(\Om;\dbR^n))$,
$v(\cdot)\in
 L^q_{\dbF}(\Om;L^2(t,T;\dbR^n))$ and $\eta\in
 L^q_{\cF_t}(\Om;\dbR^n)$.
\end{definition}

We have the following well-posedness result for
\eqref{9.25-eq5} in the sense of transposition
solution.

\begin{theorem}\label{9.25-lm2}
{\rm (\cite{LZ})} For any given $y_T \in
L^p_{\cF_T}(\Om;\dbR^m)$, the equation
\eqref{9.25-eq5} admits a unique transposition
solution $(y(\cdot), Y(\cdot))\in
D_{\dbF}([0,T];L^{p}(\Om;\dbR^m)) \t
L^{p}_{\dbF}(\Om;L^2(0,T;\dbR^m))$. Furthermore,
\begin{equation}\label{9.25-lm2-eq1}
\begin{array}{ll}\ds
|(y(\cdot),   Y(\cdot))|_{
D_{\dbF}([0,T];L^{p}(\Om;\dbR^m)) \t
L^{p}_{\dbF}(\Om;L^2(0,T;\dbR^m))}\\\ns \ds\le
\cC\left[|f(\cd,0,0)|_{
L^p_{\dbF}(\Om;L^1(0,T;\dbR^m))}+|y_T|_{
L^p_{\cF_T}(\Om;\dbR^m)}\right].
\end{array}
\end{equation}
\end{theorem}

By means of the transposition solutions $(y(\cdot),
Y(\cdot))$ and  $(P(\cdot),
Q(\cdot))$ respectively to \eqref{3.8} and
\eqref{3.9} (guaranteed by Theorem \ref{9.25-lm2}), we can
establish the following Pontryagin-type maximum
principle for Problem
(OPF).

\begin{theorem}\label{maximum principle}
{\rm (\cite{LZ})} Let (A1)--(A3) hold and $x_0\in  \dbR^n$. Let
$(\bar x(\cd),\bar u(\cd))$ be an optimal pair
of Problem (OPF). Then
\begin{equation}\label{MP-eq0}
\begin{array}{ll}\ds
\mathbb{H}(t,\bar x(t),\bar
u(t),y(t),Y(t))-\mathbb{H}(t,\bar
x(t),u,y(t),Y(t))\\
\ns\ds \q-{\frac{1}{2}}\big\langle
P(t)\big[b(t,\bar x(t),\bar u(t))-b(t,\bar
x(t),u) \big], b(t,\bar x(t),\bar
u(t))-b(t,\bar x(t),u)\big\rangle_{\dbR^n}\\
\ns\ds  \geq 0,\qq\qq\q\forall\, u\in U,\q \ae
t\in[0,T],\q \dbP\mbox{-}\as
\end{array}
\end{equation}
\end{theorem}

{\bf Sketch of the proof of Theorem \ref{maximum
principle}}\,: Since the detailed proof of
Theorem \ref{maximum principle} is too lengthy,
we shall give below only a sketch to show some
key points for establishing the stochastic
maximum principle.

Fix any $u(\cd)\in\cU[0,T]$ and $\e>0$, let
$$u^\e(t)=\left\{
\begin{array}{ll}
\bar u(t),\qq t\in[0,T]\setminus E_\e,\\
\ns\ds u(t),\qq t\in E_\e,
\end{array}
\right.
$$
where $E_\e\subseteq[0,T]$ is a measurable set
with Lebesgue measure $|E_\e|=\e$. For $\f=a,b$
and $f$, we set
\begin{equation}\left\{
\begin{array}{ll}\ds
\f_x(t) =\f_x(t,\bar x(t),\bar
u(t)),\q\f_{xx}(t) =\f_{xx}(t, \bar x(t),\bar
u(t)),\\
\ns\ds \d\f(t) =\f(t,\bar x(t),u(t))-\f(t,\bar
x(t),\cl u(t)).
\end{array}
\right.
\end{equation}
Let $x_1^\e(\cd)$ and $x_2^\e(\cd)$ solve
respectively the following stochastic
differential equations
\begin{equation}\label{MP-equ2}
\left\{
\begin{array}{ll}\ds
dx_1^\e(t)=a_x(t)x_1^\e(t)dt+
\big[b_x(t)x_1^\e(t)
+\chi_{E_\e}(t)\d b(t)\big]dW(t) &\mbox{ in } [0,T],\\
\ns\ds x_1^\e(0)=0,
\end{array}
\right.
\end{equation}
and
\begin{equation}\label{MP-equ3}
\left\{
\begin{array}{ll}\ds
dx_2^\e(t)=\[a_x(t)x_2^\e(t)+\chi_{E_\e}(t)\d
a(t)
+\frac{1}{2} a_{xx}(t)\big(x_1^\e(t), x_1^\e(t) \big)\]dt\\
\ns\ds \qq\qq+ \[b_x(t)x_2^\e(t)+
\chi_{E_\e}(t)\d b_x(t)x_1^\e(t) +\frac{1}{2}
b_{xx}(t)\big(x_1^\e(t), x_1^\e(t) \big)\]dW(t)
&
\mbox{ in } [0,T], \\
\ns\ds x_2^\e(0)=0.
\end{array}
\right.
\end{equation}
Then, by some lengthy but direct computations, one can
obtain that
\begin{equation}\label{MP-eq1}
\begin{array}{ll}\ds \cJ(u^\e(\cd)) -
\cJ(\cl u(\cd))\\
\ns\ds =\mE\big\langle h_x(\bar
x(T)),x_1^\e(T)+x_2^\e(T)\big\rangle_{\dbR^n}\!
+\frac{1}{2}\mE\big\langle h_{xx}(\cl
x(T))x_1^\e(T),x_1^\e(T)\big\rangle_{\dbR^n}\\
\ns \ds \q + \mE\!\int_0^T\!\Big\{\big\langle
g_x(t),x_1^\e(t)\!+\!x_2^\e(t)\big\rangle_{\dbR^n}\!\!+\!\frac{1}{2}\big\langle
g_{xx}(t)
x_1^\e(t),x_1^\e(t)\big\rangle_{\dbR^n}\! +\!
\chi_{E_\e}(t)\d g(t)\Big\}dt  +o(\e).
\end{array}
\end{equation}
By means of the fact that $(y(\cd),Y(\cd))$ is
the transposition solution to the equation
\eqref{3.8} with $p=2$, we find that
\begin{equation}\label{MP-eq2}
-\mE\big\langle h_x(\bar
x(T)),x_1^\e(T)\big\rangle_{\dbR^n} = \mE
\int_0^T \big[\big\langle g_x(t),
x_1^\e(t)\big\rangle_{\dbR^n} +
\chi_{E_\e}(t)\big\langle\d b(t), Y(t)
\big\rangle_{\dbR^n}\big]dt,
\end{equation}
and
\begin{equation}\label{MP-eq3}
\begin{array}{ll}\ds
-\mE\langle h_x(\bar x(T)),x_2^\e(T)\rangle_{\dbR^n} \\
\ns\ds = \mE\int_0^T\Big\{\big\langle
g_x(t),x_2^\e(t)\big\rangle_{\dbR^n}
+\frac{1}{2}\big[\big\langle
y(t),a_{xx}(t)\big(x_1^\e(t),x_1^\e(t)\big)\big\rangle_{\dbR^n}
+ \big\langle Y(t),
b_{xx}(t)\big(x_1^\e(t),x_1^\e(t)\big)\big\rangle_{\dbR^n}\big]\\
\ns \ds  \qq\qq + \chi_{E_\e}(t)\big[\big\langle
y(t),\d a(t)\big\rangle_{\dbR^n} + \big\langle
Y(t),\d b_x(t)x_1^\e(t)
\big\rangle_{\dbR^n}\big]\Big\}dt.
\end{array}
\end{equation}
Further, put $x_3^\e(t) = x_1^\e(t)
x_1^\e(t)^\top(\in\dbR^{n\times n})$. A direct
computation shows that $x_3^\e(\cd)$ solves
\begin{equation}\label{MP-equ4}
\left\{
\begin{array}{ll}\ds
dx_3^\e (t)
\!=\!\ds\Big\{a_x(t)x_3^\e(t)+x_3^\e(t)a_x(t)^\top
+ b_x(t)x_3^\e(t)b_x(t)^\top + \chi_{E_\e}(t)\d
b(t)\d b(t)^\top\\
\ns \ds \qq\qq \!\!\!\!\!\!\!+
\chi_{E_\e}(t)\[b_x(t)x_1^\e(t)\d b(t)^\top
+\d b(t)x_1^\e(t)^\top b_x(t)^\top\]\Big\}dt\\
\ns \ds \qq\qq \!\!\!\!\!\!\!+
\[b_x(t)x_3^\e(t)\!+\!x_3^\e(t)b_x(t)^\top
\!\!\!+\!\! \chi_{E_\e}\!(t) \(\d
b(t)x_1^\e(t)^\top\!\!\!+\! x_1^\e(t)\d
b(t)^\top\)\]
dW(t) & \mbox{in } (0,T],\\
 \ns\ds
x_3^\e(0)=0.
\end{array}
\right.
\end{equation}
Utilizing the fact that $(P(\cd),Q(\cd))$ is the transposition
solution to the equation \eqref{3.9} with $p=4$,
and noting that  the inner product defined in
$\dbR^{n\times n}$ is $\tr(P_1P_2^\top)$ for
$P_1,P_2\in \dbR^{n\times n}$, we find that
$$
\begin{array}{ll}\ds
-\mE \tr\big[h_{xx}(\bar x(T))x_3^\e(T)]\\
\ns\ds = \mE \int_0^T \tr \big[\chi_{E_\e}(t)\d
b(t)^\top P(t)\d b(t) - \mathbb{H}_{xx}(t,\bar
x(t),\bar u(t),y(t),Y(t))x_3^\e(t) \big]dt +
o(\e),
\end{array}
$$
which gives that
\begin{equation}\label{MP-eq4}
 \begin{array}{ll}\ds -\mE \big\langle
h_{xx}(\bar x(T))x_1^\e(T),x_1^\e(T)
\big\rangle_{\dbR^n}\\
\ns\ds = \mE \int_0^T
\big[\chi_{E_\e}(t)\big\langle P(t)\d b(t), \d
b(t)\big\rangle_{\dbR^n}   - \big\langle
\mathbb{H}_{xx}(t,\bar x(t),\bar u(t),y(t),Y(t))
x_1^\e(t),x_1^\e(t)\big\rangle_{\dbR^n} \big]dt
+ o(\e).
\end{array}
\end{equation}
From \eqref{MP-eq1}--\eqref{MP-eq4}, we obtain
that
\begin{equation}\label{MP-eq5}
\begin{array}{ll}\ds
\cJ(u^\e(\cd)) - \cJ(\bar u(\cd)) \\
\ns\ds = \mE\int_0^T
\chi_{E_\e}(t)\Big\{\big[\mathbb{H}(t,\bar
x(t),u(t),y(t),Y(t))-\mathbb{H}(t,\bar
x(t),\bar u(t),y(t),Y(t))\big]\\
\ns\ds \q  -  \frac{1}{2} \big\langle
P(t)\big[b(t,\bar x(t),u(t)) - b(t,\bar
x(t),\bar u(t)) \big], b(t,\bar x(t),u(t)) -
b(t,\bar x(t),\bar u(t))\big\rangle_{\dbR^n}
\Big\}dt   + o(\e).
\end{array}
\end{equation}
Since $\bar u(\cd)$ is the optimal control, we
have $\cJ(u^\e(\cd)) - \cJ(\bar u(\cd))\geq 0$.
This, together with \eqref{MP-eq5}, yields that
\begin{equation}\label{MP-eq5zz}
\begin{array}{ll}\ds
\mE\int_0^T
\chi_{E_\e}(t)\Big\{\big[\mathbb{H}(t,\bar
x(t),u(t),y(t),Y(t))-\mathbb{H}(t,\bar
x(t),\bar u(t),y(t),Y(t))\big]\\
\ns\ds \q - \frac{1}{2} \big\langle
P(t)\big[b(t,\bar x(t),u(t))\!-\!b(t,\bar
x(t),\bar u(t)) \big], b(t,\bar
x(t),u(t))\!-\!b(t,\bar
x(t),\bar u(t))\big\rangle_{\dbR^n} \Big\}dt \\
\ns\ds  \geq o(\e),
\end{array}
\end{equation}
which leads to \eqref{MP-eq0}.
\endpf

\ms

For some optimal controls, it may happen that
the first-order necessary conditions turn out to
be trivial. When an optimal control is singular,
the first-order necessary condition cannot
provide enough information for the theoretical
analysis and numerical computation, and
therefore one needs to study the second order
necessary conditions. Quite different from the
deterministic setting, there exist some
essential difficulties in deriving the pointwise
second-order necessary condition from an
integral-type one when the diffusion term of the
control system contains the control variable,
even for the case of convex control constraint
(see the first four paragraphs of subsection 3.2
in \cite{ZZ1} for a detailed explanation). In
\cite{ZZ1, ZZ2}, these difficulties were
overcome by means of some technique from the
Malliavin calculus, and some pointwise
second-order necessary conditions for stochastic
optimal controls were established, even for the
general case when the control region is
nonconvex but the  full picture is still quite
unclear (see \cite{FZZ} for some recent
progresses).


\section{Controllability of stochastic differential equations in infinite dimensions: An analysis of a typical equation}\label{s5}


This section is devoted to studying the
controllability of stochastic differential
equations in infinite dimensions. Since the
stochastic controllability problem is even less
understood in finite dimensions, we shall
concentrate only on a typical equation, i.e., a
stochastic parabolic system. Our main results
can be described as follows:

\begin{itemize}
\item When the coefficients of the underlying
system are space-independent,  using the
spectral method, we show the null/approximate
controllability using only one control applied
to the drift term;

\item The null/approximate controllability of general stochastic parabolic systems with two controls are shown by means of duality argument.

\end{itemize}
\no In each of the above cases, we shall explain
the main differences between the deterministic
problem and its stochastic counterpart.

\subsection{Formulation of the problem}\label{sto heat se1}

Throughout this section, we assume that $\mathbf{F}$ is the natural filtration generated
by $W(\cd)$, $G \subset \dbR^{n}$
($n\in \dbN$) is a given bounded domain with a
$C^\infty$ boundary $\Gamma $, and $G_0$ is a
given nonempty open subset of $G$. Denote by
$\chi_{G_0}$ the characteristic function of
$G_0$ in $G$. Put $ Q\= (0,T)\times G$,
$\Si\=(0,T)\t\Gamma$ and $Q_0\=(0,T)\times G_0$.
Also, we assume that $a^{jk}: \; \cl{G}\to
\dbR^{n\times n}\;$ ($j,k=1,2,\cdots,n$)
satisfies $a^{jk}\in C^1(\ol G))$,
$a^{jk}=a^{kj}$, and for some $s_0>0$,
\begin{equation}\label{h1}
\sum_{j,k=1}^n a^{jk}(x)\xi_{j}\xi_{k} \geq
s_0|\xi|^{2},\q\forall\; (x,\xi)\equiv
(x_1,\cdots,x_n,\xi_{1},\cdots,\xi_{n}) \in
G\times \dbR^{n}.
\end{equation}

Let us fix an $m\in\dbN$ and consider the following
controlled stochastic parabolic system:
\begin{equation}\label{heat 1.1}
\left\{
\begin{array}{ll}
\ds dy-\sum_{j,k=1}^n(a^{jk}y_{x_j})_{x_k}dt =
\Big(\sum_{j=1}^na_{1j} y_{x_j} +a_2 y+\chi_{G_0}u\Big)dt+(a_3y+v)\, dW(t)&\hb{ in }Q,\\
\ns
\ds y=0&\hb{ on }\Si,\\
\ns
 y(0)=y_0&\hb{ in } G,
 \end{array}
 \right.
 \end{equation}
where
\bel{20160618e1} \left\{\begin{array}{ll} a_{1j}
\in
L^{\infty}_{\dbF}(0,T;W^{1,\infty}(G;\dbR^{m\times m})),\q j=1,2,\cdots,n,\\
a_2\in
L^{\infty}_{\dbF}(0,T;L^{\infty}(G;\dbR^{m\times
m})), \q a_3\in L^{\infty}_{\dbF}(0,T;
L^{\infty}(G;\dbR^{m\times m})).
\end{array}
\right. \ee In the system \eqref{heat 1.1}, the
initial state $y_0\in
L_{\cF_0}^2(\Omega;L^2(G;\dbR^m))$, $y$ is the
state variable, and the control variable
consists of a pair $(u, v)\in
L^2_{\dbF}(0,T;L^2(G_0;$ $\dbR^m))\times
L^2_{\dbF}(0,T;L^2(G;\dbR^m))$.

We first recall the following well-posedness
result for the equation \eqref{heat 1.1}. The
proof can be found in \cite[Chapter 6]{Prato}
and \cite[Chapter 3]{LZ3}.
\begin{lemma}\label{wellposed of heat 1.1}
Let $a_{1j} \in
L^{\infty}_{\dbF}(0,T;L^\infty(G;\dbR^{m\times
m}))$ for $j=1,2,\cdots,n$, and $a_2$ and $a_3$
be given as in \eqref{20160618e1}. Then, for any
$y_0 \in L^2(G;\dbR^m)$ and $(u, v)\in
L^2_{\dbF}(0,T;L^2(G_0;\dbR^m))\times
L^2_{\dbF}(0,T;$ $L^2(G;\dbR^m))$, the system
\eqref{heat 1.1} admits a unique weak solution
$y\in
L_{\dbF}^2(\Omega;C([0,T];L^2(G;\dbR^m)))\cap
L_{\dbF}^2(0,T;H_0^1(G;\dbR^m))$. Moreover,
\begin{equation}\label{20160616e3} \ba{ll}\ds |y|_{
L_{\dbF}^2(\Omega;C([0,T];L^2(G;\dbR^m)))\cap
L_{\dbF}^2(0,T;H_0^1(G;\dbR^m))}\\ \ns \ds\le
\cC \big(|y_0|_{L^2(G;\dbR^m)}+|(u,
v)|_{L^2_{\dbF}(0,T;L^2(G_0;\dbR^m))\times
L^2_{\dbF}(0,T;L^2(G;\dbR^m))}\big). \ea
\end{equation}
\end{lemma}
\begin{definition}
The system \eqref{heat 1.1} is said to be null
controllable if for any $y_0\in L^2(G;\dbR^m)$,
there exists a pair of $(u,v)\in  L^2_{\dbF}(0,T;L^2(G_0;\dbR^m))\times L^2_{\dbF}(0,T;L^2(G;\dbR^m))$
such that the corresponding solution to \eqref{heat 1.1} fulfills
that $y(T)=0$, $\dbP$-a.s.
\end{definition}

Note that we introduce two controls $u$ and $v$
in \eqref{heat 1.1}. In view of the
controllability result for the deterministic
parabolic equation, it is more natural to use
only one control and consider the following
controlled stochastic parabolic system (which is
a special case of  \eqref{heat 1.1} with
$v\equiv0$):
\begin{equation}\label{heat---1.1}
\left\{
\begin{array}{ll}
\ds dy-\sum_{j,k=1}^n(a^{jk}y_{x_j})_{x_k}dt =
\Big(\sum_{j=1}^na_{1j} y_{x_j} +a_2 y+\chi_{G_0}u\Big)dt+a_3y dW(t)\q&\hb{ in }Q,\\
\ns
\ds y=0&\hb{ on }\Si,\\
\ns
 y(0)=y_0&\hb{ in } G.
 \end{array}
 \right.
 \end{equation}
It is easy to see that, the dual system of both
\eqref{heat 1.1} and \eqref{heat---1.1} is  the
following backward stochastic parabolic system:
\begin{equation}\label{dual heat 1.1}
\left\{
\begin{array}{ll}
\ds dz+\sum_{j,k=1}^n(a^{jk}z_{x_j})_{x_k}dt
=\Big[\sum_{j=1}^n\big(a_{1j}^\top
z\big)_{x_j}-a_2^\top z- a_3^\top Z\Big]dt+Z
dW(t) \q&\hb{ in }Q,\\
\ns
z=0&\hb{ on }\Si,\\
\ns z(T)=z_T&\hb{ in }G.
\end{array}
\right.
\end{equation}

We have the following well-posedness result for
the equation \eqref{dual heat 1.1} (See
\cite[Chapter 4]{LZ3} for example).
\begin{proposition}\label{20160619c1}
Under the condition \eqref{20160618e1}, for any
$z_T \in L_{\cF_T}^2(\Omega;L^2(G;$ $\dbR^m))$,
the system \eqref{dual heat 1.1} admits one and
only one weak solution $(z,Z)\in
\big(L^2_{\dbF}(\Omega;C([0,T];L^2(G;\dbR^m)))$
$\bigcap L_{\dbF}^2(0,T;H_0^1(G;\dbR^m))\big)$
$\times L_{\dbF}^2(0,T;L^2(G;\dbR^m))$.
Moreover, for any $t\in [0,T]$,
\begin{equation}\label{20160619e1}
\begin{array}{ll}
\ds|(z(\cd),Z(\cd))|_{
\left(L^2_{\dbF}(\Omega;C([0,t];L^2(G;\dbR^m)))\cap
L_{\dbF}^2(0,t;H_0^1(G;\dbR^m))\right)\times L_{\dbF}^2(0,t;L^2(G;\dbR^m))}\le \cC |z(t)|_{L_{\cF_t}^2(\Omega;L^2(G;\dbR^m))}.
\end{array}
\end{equation}
\end{proposition}

In order to obtain the null controllability of
\eqref{heat---1.1},  we need to prove that
solutions to the system \eqref{dual heat 1.1}
satisfy the following observability estimate:
\begin{equation}\label{btrulyheat obser1}
|z(0)|_{L_{\cF_0}^2(\Omega;L^2(G;\dbR^m))}\le
\cC |z|_{L^2_{\dbF}(0,T;L^2(G_0;\dbR^m))},\q
\forall\;z_T\in
L^2_{\cF_T}(\Omega;L^2(G;\dbR^m)).
\end{equation}
Unfortunately, at this moment, we are not able
to prove the observability estimate
\eqref{btrulyheat obser1} for the general case. Instead, we obtain a
weak version of \eqref{btrulyheat obser1}, i.e.,
a weak observability estimate (for the system
\eqref{dual heat 1.1}) in Theorem \ref{b heat
ob1} (See Subsection \ref{sto heat se2}). By
duality, Theorem \ref{b heat ob1} implies the
null controllability of  \eqref{heat 1.1}.

There exists a main difficulty to establish
\eqref{btrulyheat obser1}, that is, though the
correction   term ``$Z$" plays a ``coercive"
role for the well-posedness of \eqref{dual heat
1.1}, it seems to be a ``bad" (non-homogeneous)
term when one tries to prove \eqref{btrulyheat
obser1} using the global Carleman estimate.

Nevertheless, based on the spectral method, for some special case, we are able to show the controllability of \eqref{heat---1.1}.

\subsection{Controllability of a class of
stochastic parabolic systems}\label{null
particular1}

In this subsection,  we show that when the
coefficients of the stochastic parabolic system
are space-independent, it is
null/approximately controllable using only one
control applied to the drift term. These results
were first proved in \cite{Lu1}.

We consider the following stochastic parabolic
system:
\begin{equation}\label{sheatsystem1}
\left\{
\begin{array}{lll}
\ds dy - \sum_{j,k=1}^n (a^{jk}y_{x_j})_{x_k}dt = [a(t)y + \chi_{E}(t)\chi_{G_0}(x) u]dt+b(t)y dW(t) \q& \mbox{ in } Q, \\
\ns\ds y = 0  & \mbox{ on } \Si, \\
 \ns\ds y(0) = y_0 &\mbox{ in } G,
\end{array}
\right.
\end{equation}
where  $a(\cd)\in
L_{\dbF}^{\infty}(0,T;\dbR^{m\times m})$ and
$b(\cd)\in L_{\dbF}^{\infty}(0,T;\dbR^{m\times
m})$ are given, $E$ is a fixed Lebesgue
measurable subset in $(0,T)$ with a positive
Lebesgue measure (i.e., $\dbm(E)>0$). In
\eqref{sheatsystem1}, $y$ is the state variable
(valued in $L^2(G;\dbR^m)$),  $y_0 \in
L^2(G;\dbR^m)$ is the initial state, $u$ is the
control variable, and the control space is $
L^{\infty}_\dbF(0,T;L^2(\Omega;L^{2}(G;\dbR^m)))$.

\begin{definition}
The system \eqref{sheatsystem1} is said to be
null controllable at time $T$ if for any $y_0\in
L^2(G;\dbR^m)$, there exists a $u\in
L^{\infty}_\dbF(0,T;L^2(\Omega;L^{2}(G;\dbR^m)))$
such that the corresponding solution to \eqref{sheatsystem1} fulfills
that $y(T)=0$, $\dbP$-a.s.
\end{definition}
We have the following null controllability
result for the system \eqref{sheatsystem1}.

\begin{theorem} \label{thnull}
The system (\ref{sheatsystem1}) is null
controllable at time $T$.
\end{theorem}

\br\label{20160520r1} When $E=(0,T)$, one can
use the global Carleman estimate to prove the
corresponding null controllability result for
the deterministic counterpart of
(\ref{sheatsystem1}). However, at least at this
moment we do not know how to use a similar
method to prove Theorem \ref{thnull} even for
the same case that $E=(0,T)$. \er

Next, we consider the  approximate
controllability for the system
\eqref{sheatsystem1} under a stronger assumption
on the controller $E\times G_0$ than that for
the null controllability.

\begin{definition}
The system \eqref{sheatsystem1} is said to be
approximately controllable at time $T$ if for any
initial datum $y_0\in L^2(G;\dbR^m)$, any final
state $y_T\in L^2_{\cF_T}(\Om;L^2(G;\dbR^m))$ and
any $\varepsilon>0$, there exists a control
$u\in
L^{\infty}_\dbF(0,T;L^2(\Omega;L^{2}(G;\dbR^m)))$
such that the corresponding solution to
\eqref{sheatsystem1} satisfies that
$|y(T)-y_T|_{L^2_{\cF_T}(\Om;L^2(G;\dbR^m))}\leq
\varepsilon$.
\end{definition}

\begin{theorem}\label{thapp}
The system (\ref{sheatsystem1}) is approximately
controllable at time $T$ if and only if
$\dbm((s,T)\cap E)>0$ for any $s\in [0,T)$.
\end{theorem}

At the first glance, it seems that Theorem
\ref{thapp} is unreasonable.  If $b(\cd)\equiv
0$, then the system (\ref{sheatsystem1}) is like
a deterministic parabolic equation with a random
parameter. The readers may guess that one can
obtain the approximate controllability by only
assuming that $\dbm(E)>0$. However, this is not
the case. The reason for this comes from our
definition of the approximate controllability
for the system (\ref{sheatsystem1}). We expect
that any element belonging to
$L^2_{\cF_T}(\Omega;L^2(G))$ rather than
$L^2_{\cF_s}(\Omega;L^2(G))$ ($s<T$) can be
attached as close as one wants. Hence we need to
put the control $u$ to be active until the time
$T$.

In some sense, it is surprising that one needs a
little more assumption in Theorem \ref{thapp}
for the approximate controllability of
\eqref{sheatsystem1} than that in Theorem
\ref{thnull} for the null controllability.
Indeed, it is well-known that in the
deterministic setting, the null controllability
is usually stronger than the approximate
controllability. But this does not remain to be
true in the stochastic case. Actually, from
Theorem \ref{thapp}, we see that the additional
condition (compared to the null controllability)
that $\dbm((s,T)\bigcap E)>0$ for any $s\in
[0,T)$ is not only sufficient but also necessary
for the approximate controllability of
(\ref{sheatsystem1}). Therefore, in the setting
of stochastic distributed parameter systems, the
null controllability does NOT imply the
approximate controllability. This indicates that
there exists some essential difference between
the controllability theory of the
deterministic parabolic equations and its
stochastic counterpart.

\subsubsection{Some preliminaries}

Before proving Theorems \ref{thnull} and
\ref{thapp}, we give some preliminary results.
To begin with, we recall the following known
property about Lebesgue measurable sets.

\begin{lemma} {\rm (\cite[pp. 256--257]{J.Lions})} \label{mea lemma}
For a.e. $\tilde t \in E$, there exists a
sequence of numbers $\{ t_{i} \}_{i =
1}^{\infty}\subset (0,T)$ such that
 \begin{eqnarray}
&&t_{1} < t_{2} < \cdots < t_{i} < t_{i+1} <
\cdots <\tilde  t, \qq t_{i} \to\tilde t \mbox{
as } i \to \infty, \label{t2.5}
\\
\ns\ds && \dbm(E\cap [t_{i},t_{i+1}]) \geq
\rho_1 (t_{i+1}-t_{i}),\quad i = 1,2,\cdots,
\label{t2.6}
\\
\ns\ds &&\frac{t_{i+1}-t_{i}}{t_{i+2}-t_{i+1}}
\leq \rho_2, \quad i = 1,2,\cdots, \label{t2.7}
\end{eqnarray}
where $\rho_1$ and $\rho_2$ are two positive
constants which are independent of $i$.
 \end{lemma}

Next, we give the following result (which is a
Riesz-type Representation Theorem for the dual
of space $L^p_\dbF(0,T;L^q(\Omega;H))$). Its
proof can be found in \cite{LYZ} or
\cite[Chapter 1]{LZ3}.

\begin{lemma}\label{rep}
Suppose $1\le p,q<\infty$, and that $H$ is a
Hilbert space. Then
\begin{equation}
L^p_\dbF(0,T;L^q(\Omega;H))^*=L^{p'}_\dbF(0,T;L^{q'}(\Omega;H)).
\end{equation}
Here, $p'$ and $q'$ are respectively the (usual H\"older) conjugate numbers of $p$ and $q$.
\end{lemma}

Next, let us define an unbounded operator $A$ on
$L^2(G)$ as follows:
\begin{equation}\label{defA}
\left\{
\begin{array}{ll}
\ds D(A)=  H^2(G)\cap H_0^1(G), \\
\ds Ah = -\sum^n_{j,k=1}(a^{jk}
h_{x_j})_{x_k},\q\forall\;h\in D(A).
\end{array}
\right.
\end{equation}
Let $\{\lambda_i\}^\infty_{i=1}$ be the
eigenvalues of $A$, and $\{e_i\}^{\infty}_{i=1}$
be the corresponding eigenfunctions satisfying
$|e_i|_{L^2(G)} = 1,\,i = 1,2,3 \cdots $. It
holds that $0< \lambda_1\le \lambda_2\le
\lambda_3\le \cdots\le \lambda_k\le \cdots\to
\infty$. For any $r\ge \l_1$, write
$\L_r=\big\{i\in\dbN\;\big|\;\l_i\le r\big\}$.
We recall the following observability estimate
(for partial sums of the eigenfunctions of $A$),
established in \cite[Theorem 1.2]{Luqi5} (See
also \cite[Theorem 3]{Lebeau-Zuazua1} for a
special case of this result).

\begin{lemma} \label{ob ellip}
There exist two positive constants $\cC_{1}\ge
1$ and $\cC_{2}\ge 1$ such that
\begin{equation}\label{lr1}
\sum_{i\in \L_r}|a_{i}|^{2} \leq
\cC_{1}e^{\cC_{2}\sqrt{r}}\int_{G_0}\Big|
\sum_{i\in \L_r}a_{i}e_{i}(x) \Big|^2dx
\end{equation}
holds for any $r\ge \l_1$ and $a_i\in \dbC$ with
$i\in \L_r$.
\end{lemma}

Further, for any $s_1$ and $s_2$ satisfying
$0\leq s_1<s_2\leq T$, we introduce the
following backward stochastic parabolic system:
\begin{equation}\label{sbheatsystem2}
\left\{
\begin{array}{lll}
\ds dz + \sum_{j,k=1}^n (a^{jk}z_{x_j})_{x_k}dt
= -[a(t)^\top z+b(t)^\top Z]dt + ZdW(t) & \mbox{ in }\, (s_1,s_2)\times G, \\
\ns\ds z = 0  & \mbox{ on }\, (s_1,s_2)\times \Gamma, \\
 \ns\ds z(s_2)=\eta  & \mbox{ in }\, G,
\end{array}
\right.
\end{equation}
where $\eta\in
L^2_{\cF_{s_2}}(\Omega;L^2(G;\dbR^m))$.

Put
 $$r_0 =2
|a|_{L_{\dbF}^{\infty}(0,T;\dbR^{m\times
m})}+|b|^2_{L_{\dbF}^{\infty}(0,T;\dbR^{m\times
m})}.
 $$
For each $r\ge\l_1$, we set $H_{r} =
\mbox{\span}\{e_i\;|\;\l_i\le r\}$ and denote by
$\Pi_{r}$ the orthogonal projection from
$L^{2}(G)$ to $H_{r}$. Write
 \bel{20160605e2}
 H_{r}^m=\overbrace{H_r\times
H_r\times \cdots\times H_r}^{m\hb{ \tiny
times}}.
 \ee
To simplify the notation, we also denote by
$\Pi_{r}$ the orthogonal projection from
$L^{2}(G;\dbR^m)$ to $H_{r}^m$. We need the
following observability result for
\eqref{sbheatsystem2} with the final data
belonging to $L^2_{\cF_{s_2}}(\Omega; H_{r}^m)$,
a proper subspace of
$L^2_{\cF_{s_2}}(\Omega;L^2(G;\dbR^m))$.

\begin{proposition}\label{conprop1}
For each $r\ge\l_1$, the solution to  the system
(\ref{sbheatsystem2}) with $\eta \in
L^2_{\cF_{s_2}}(\Omega; H_{r}^m)$ satisfies that
\bel{coneq11}
\mathbb{E}\big|z(s_1)\big|^2_{L^{2}(G;\dbR^m)}
\leq \frac{\cC_{1}e^{\cC_2\sqrt{r}+r_0
(s_2-s_1)}}{(\dbm(E\cap[s_1,s_2]))^2} \big|
\chi_{E} \chi_{G_0} z
\big|^{2}_{L^{1}_\dbF(s_1,s_2;L^2(\Omega;L^{2}(G;\dbR^m)))},
\ee
whenever $\dbm(E\cap[s_1,s_2])\not=0$.
\end{proposition}

{\it Proof}\,: Each $\eta \in
L^2_{\cF_{s_2}}(\Omega; H_{r}^m)$ can be written
as $ \ds \eta = \sum_{i\in \L_r}\eta_{i}e_{i}(x)
$ for some $\eta_i \in L^2_{\cF_{s_2}}(\Omega;
\dbR^m)$ with $i\in \L_r$. The solution $(z,Z)$
to (\ref{sbheatsystem2}) can be expressed as
$$
z= \sum_{i\in \L_r}z_{i}(t)e_{i},\qq
Z=\sum_{i\in \L_r}Z_{i}(t)e_{i},
$$
where $z_i(\cdot)\in
C_{\dbF}([s_1,s_2];L^2(\Omega;\dbR^m))$ and
$Z_i(\cdot)\in L^2_{\dbF}(s_1,s_2;\dbR^m)$, and
satisfy the following equation
\begin{eqnarray*}
\left\{
\begin{array}{lll}\ds
dz_i -\l_i z_i dt = -[a(t)^\top z_i+b(t)^\top Z_i]dt + Z_idW(t) \quad & \mbox{ in } \,[s_1,s_2],\\
\ns\ds z_i(T) = \eta_i.
\end{array}
\right.
\end{eqnarray*}
By Lemma \ref{ob ellip}, for any $t \in
[s_1,s_2]$, we have \bel{conprop1eq1} \ba{ll}
\ds\mathbb{E}\int_G|z(t)|_{\dbR^m}^{2}dx
\3n&\ds=\mathbb{E}\sum_{i\in
\L_r}|z_{i}(t)|_{\dbR^m}^{2} \ds\leq
\cC_{1}e^{\cC_{2}\sqrt{r}}\mathbb{E}\int_{G_0}\|\sum_{i\in
\L_r}
   z_{i}(t)e_i\|_{\dbR^m}^{2}dx
\\
\ns&\ds =
  \cC_{1}e^{\cC_{2}\sqrt{r}}\mathbb{E}\int_{G_0}|z(t)|_{\dbR^m}^{2}dx.
\ea \ee

By It\^{o}'s formula, we find that
$$
d(e^{r_0  t}|z|_{\dbR^m}^2) = r_0  e^{r_0
t}|z|_{\dbR^m}^2+e^{r_0  t}\big({\lan
dz,z\ran}_{\dbR^m}+{\lan z,dz\ran}_{\dbR^m}\big)
+ e^{r_0 t}|dz|_{\dbR^m}^2.
$$
Hence,
\begin{equation}\label{conprop1eq21}
\begin{array}{ll}\ds
\mathbb{E}\(e^{r_0 t}\int_G|z(t)|_{\dbR^m}^2dx\)
- \mathbb{E}\(e^{r_0  s_1}\int_G|z(s_1)|_{\dbR^m}^2dx\) \\
\ns\ds =r_0\mathbb{E}\int_{s_1}^t\int_G e^{r_0
s}|z(s)|_{\dbR^m}^2dxds +2 \sum_{i\in
\L_r}\mathbb{E}\int_{s_1}^te^{r_0
s}\l_i|z_{i}(s)|_{\dbR^m}^{2} ds
\\
\ns\ds\q+  \mathbb{E}\int_{s_1}^t\int_G e^{r_0  s}\big(-{\langle a(s)^\top z(s)+b(s)^\top Z(s),z(s)\rangle}_{\dbR^m}\\
\ns\ds\q-{\langle z(s),a(s)^\top z(s)+b(s)^\top Z(s)\rangle}_{\dbR^m}+|Z(s)|_{\dbR^m}^2\big)dxds  \\
\ns\ds  \geq 2\sum_{i\in
\L_r}\mathbb{E}\int_{s_1}^te^{r_0
s}\l_i|z_{i}(s)|_{\dbR^m}^{2} ds \geq 0.
\end{array}
\end{equation}
From  (\ref{conprop1eq1}) and
(\ref{conprop1eq21}), we obtain that, for any $t
\in [s_1,s_2]$,
\bel{160605e3}
\mathbb{E}\int_{G}|z(s_1,x)|_{\dbR^m}^2dx \leq
\cC_{1}e^{\cC_{2}\sqrt{r}+r_0
(s_2-s_1)}\mathbb{E}\int_{G_0}|z(t,x)|_{\dbR^m}^{2}dx.
\ee
By \eqref{160605e3}, it follows that
$$
\begin{array}{ll}\ds
\int_{E\cap[s_1,s_2]}\[\mathbb{E}\int_{G}|z(s_1,x)|_{\dbR^m}^2dx\]^{\frac{1}{2}}dt
\\
\ns\ds\leq  \big(\cC_{1}e^{\cC_{2}\sqrt{r}+r_0
(s_2-s_1)}\big)^{\frac{1}{2}}\int_{E\cap[s_1,s_2]}\[\mathbb{E}\int_{G_0}|z(t,x)|_{\dbR^m}^{2}dx\]^{\frac{1}{2}}dt.
\end{array}
$$
Hence, when $\dbm(E\cap[s_1,s_2])\not=0$, we
obtain that for each $\eta \in
L^2_{\cF_{s_2}}(\Omega; H_{r}^m)$,
 $$
\begin{array}{ll}\ds
\mathbb{E}\int_{G} |z(s_1,x)|_{\dbR^m}^2dx \\
\ns\ds\leq {\frac{\cC_{1}e^{\cC_{2}\sqrt{r}+r_0
(s_2-s_1) }}{(\dbm(E\cap[s_1,s_2]))^{2}}}
\Big\{\int^{s_2}_{s_1} \[\mathbb{E}\int_{G}
|{\chi_{E}(t)}{\chi_{G_0}(x){z(t,x)}}|_{\dbR^m}^{2}dx\]^{\frac{1}{2}}dt\Big\}^{2}
 \\
\ns \ds= {\frac{\cC_{1}e^{\cC_{2}\sqrt{r}+r_0
(s_2-s_1)
}}{(\dbm(E\cap[s_1,s_2]))^{2}}}|{\chi_{E}}{\chi_{G_0}}{z}|^2_{L^{1}_\dbF(s_1,s_2;L^2(\Omega;L^{2}(G;\dbR^m)))},
\end{array}
$$
which gives (\ref{coneq11}).
\endpf

\medskip

By means of the usual duality argument,
Proposition \ref{conprop1} yields a partial
controllability result for the following
controlled system:
\begin{equation}\label{sheatsystem101}
\left\{
\begin{array}{lll}
\ds dy - \sum_{j,k=1}^n (a^{jk}y_{x_j})_{x_k}dt
 =  [a(t)y + \chi_{E}\chi_{G_0} u]dt+b(t)y dW(t) \q & \mbox{ in } (s_1,s_2)\times G, \\
\ns\ds y = 0  & \mbox{ on } (s_1,s_2)\times \Gamma, \\
 \ns\ds y(s_1) = y_{s_1} &\mbox{ in } G,
\end{array}
\right.
\end{equation}
where $y_{s_1}\in
L^2_{\cF_{s_1}}(\Omega;L^2(G;\dbR^m))$. That is,
we have the following result.

\begin{proposition}\label{conprop2}
If $\dbm(E\cap[s_1,s_2])\not=0$, then for every
$r\ge \l_1$ and $y_{s_1}\in
L^2_{\cF_{s_1}}(\Omega;L^2(G;\dbR^m))$, there
exists a control $u_r \in
 L^{\infty}_\dbF(s_1,s_2;L^2(\Omega;L^{2}(G;\dbR^m)))$
such that the solution $y$ to the system
\eqref{sheatsystem101} with  $u = u_{r}$
satisfies that $\Pi_{r}(y(s_2)) = 0 $,
$\dbP$-a.s. Moreover, $u_{r}$ verifies that
\begin{equation}\label{conprop2est11}
|u_{r}|^2_{L^{\infty}_\dbF(s_1,s_2;L^2(\Omega;L^2(G;\dbR^m)))}
\leq {\frac{\cC_{1}e^{\cC_{2}\sqrt{r}+r_0
(s_2-s_1)}}{(\dbm(E\cap[s_1,s_2]))^{2}}}
|y_{s_1}|_{L^2_{\cF_{s_1}}(\Omega;L^2(G;\dbR^m))}^{2}.
\end{equation}
\end{proposition}

{\it Proof}\,: Define a subspace $H$ of
$L^{1}_\dbF(s_1,s_2;L^2(\Omega;L^{2}(G;\dbR^m)))$:
$$
H=\left\{f=\chi_E\chi_{G_0}z\;\Big|\;\ (z,Z)
\mbox{ solves }\eqref{sbheatsystem2}\mbox{ for
some }\eta\in L^2_{\cF_{s_2}}(\Omega; H_{r}^m)
\right\}
$$
and a linear functional $\mathcal{L}$ on $H$:
$$
\mathcal{L}(f)=-\mathbb{E}\int_{G} {\lan
y_{s_1},z(s_1)\ran}_{\dbR^m}dx.
$$
By Proposition \ref{conprop1}, it is easy to
check that $\mathcal{L}$ is a bounded linear
functional on $H$ and
$$
|\mathcal{L}|^2\leq
{\frac{\cC_{1}e^{\cC_{2}\sqrt{r}+r_0
(s_2-s_1)}}{(\dbm(E\cap[s_1,s_2]))^{2}}}
|y_{s_1}|_{L^2_{\cF_{s_1}}(\Omega;L^2(G;\dbR^m))}^{2}.
$$
By the Hahn-Banach Theorem, $\mathcal{L}$ can be
extended to a  bounded linear functional
$\wt{\mathcal{L}}$ (satisfying
$\big|\wt{\mathcal{L}}\big|$=$|\mathcal{L}|$) on
$L^{1}_\dbF(s_1,s_2;$
$L^2(\Omega;L^{2}(G;\dbR^m)))$. By Lemma
\ref{rep}, there exists a control $u_r\in
L^{\infty}_\dbF(s_1,s_2;$
$L^2(\Omega;L^{2}(G;\dbR^m)))$ such that
$$
\mathbb{E}\int_{s_1}^{s_2}\int_{G} {\lan
u_r,f\ran}_{\dbR^m}dxdt=\wt{\mathcal{L}}(f),\quad\forall\;
f\in
L^{1}_\dbF(s_1,s_2;L^2(\Omega;L^{2}(G;\dbR^m))).
$$
In particular, for any $\eta\in
L^2_{\cF_{s_2}}(\Omega; H_{r}^m)$, the
corresponding solution $(z,Z)$ to
\eqref{sbheatsystem2} satisfies
\bel{20160503eq1}
\mathbb{E}\int_{s_1}^{s_2}\int_{G} {\lan
u_r,\chi_E\chi_{G_0}z\ran}_{\dbR^m}dxdt=-\mathbb{E}\int_{G}
{\lan y_{s_1},z(s_1)\ran}_{\dbR^m}dx. \ee
Applying It\^o's formula to ${\lan
y,z\ran}_{\dbR^m}$, where $y$ solves the system
\eqref{sheatsystem101} with  $u = u_{r}$, we
obtain that
 \bel{20160503eq2}
 \ba{ll}\ds
 \mathbb{E}\int_{G} {\lan  y(s_2),\eta\ran}_{\dbR^m}dx-\mathbb{E}\int_{G} {\lan  y_{s_1},z(s_1)\ran}_{\dbR^m}dx=\mathbb{E}\int_{s_1}^{s_2}\int_{G} {\lan
\chi_E\chi_{G_0} u_r,z\ran}_{\dbR^m}dxdt.
 \ea
 \ee
Combining \eqref{20160503eq1} and
\eqref{20160503eq2}, we arrive at
 $$
\mathbb{E}\int_{G} {\lan
y(s_2),\eta\ran}_{\dbR^m}dx=0,\qq\forall\;\eta\in
L^2_{\cF_{s_2}}(\Omega; H_{r}^m),
 $$
which implies that $\Pi_{r}(y(s_2)) = 0 $,
$\dbP$-a.s. Moreover,
$|u_r|_{L^{\infty}_\dbF(s_1,s_2;L^2(\Omega;L^{2}(G;\dbR^m)))}$
$=|\mathcal{L}|$, which yields
\eqref{conprop2est11}.\endpf

\ms

Finally, for any $s\in[0,T)$, we consider the
following equation:
\begin{eqnarray}\label{sheatsystem11}
\left\{
\begin{array}{lll}
\ds dy - \sum_{j,k=1}^n (a^{jk}y_{x_j})_{x_k}dt
=  a(t)ydt +b(t)y dW(t) \q & \mbox{ in }  (s,T)\times G, \\
\ns\ds y = 0  & \mbox{ on } (s,T)\times \Gamma, \\
 \ns\ds y(s) = y_s &\mbox{ in } G,
\end{array}
\right.
\end{eqnarray}
where $y_s\in
L^2_{\cF_s}(\Omega;L^2(G;\dbR^m))$. Let us show
the following decay result for the system
\eqref{sheatsystem11}.

\begin{proposition}\label{energydecay}
Let $r\ge \l_1$. Then, for any $y_s\in
L^2_{\cF_s}(\Omega;L^2(G;\dbR^m))$ with
$\Pi_{r}(y_s)=0$, $\dbP$-a.s., the corresponding
solution $y$ to (\ref{sheatsystem11}) satisfies
that
\begin{equation}\label{aaasss}
\mathbb{E}|y(t)|_{L^2(G;\dbR^m)}^2 \leq
e^{-(2r-r_0 )(t-s)}
|y_s|_{L^2_{\cF_s}(\Omega;L^2(G;\dbR^m))}^2,
\q\forall\, t\in[ s,T].
\end{equation}
\end{proposition}

{\it Proof}\,: Since $y_s\in
L^2_{\cF_s}(\Omega;L^2(G;\dbR^m))$ satisfying
$\Pi_{r}(y_s)=0$, we see that $y_s =
\sum_{i\in\dbN\setminus \L_r}y_s^ie_i$ for some
$ y_s^i\in L^2_{\cF_s}(\Omega;\dbR^m)$ with
$i\in\dbN\setminus \L_r$. Clearly, the solution
$y$ to (\ref{sheatsystem11}) can be expressed as
$ y= \sum_{i\in\dbN\setminus \L_r}
y^{i}(t)e_{i}$, where $y^{i}(\cdot)\in
C_{\dbF}([s,T];L^2(\Omega;\dbR^m)$ solves the
following stochastic differential equation:
\begin{eqnarray*}
\left\{
\begin{array}{lll}\ds
dy^{i} +\l_i y^{i} dt = a(t)y^idt +b(t)y^{i}dW(t) & \mbox{ in } \,[s,T],\\
\ns\ds y^{i}(s) = y_s^{i}.
\end{array}
\right.
\end{eqnarray*}

By It\^{o}'s formula, we have that
$$
\begin{array}{ll}\ds
 d(e^{(2r-r_0 )(t-s)}|y|_{\dbR^m}^2)\3n&\ds = e^{(2r-r_0 )(t-s)} \big({\lan dy,y\ran}_{\dbR^m}+{\lan y,dy\ran}_{\dbR^m}\big)\\
\ns&\ds \q  + e^{(2r-r_0 )(t-s)}|dy|_{\dbR^m}^2+
(2r-r_0 )e^{(2r-r_0 )(t-s)}|y|_{\dbR^m}^2.
\end{array}
$$
Hence, by $\l_i>r$ for each $i\in\dbN\setminus
\L_r$ and recalling that $r_0 =2
|a|_{L_{\dbF}^{\infty}(0,T;\dbR^{m\times
m})}+|b|^2_{L_{\dbF}^{\infty}(0,T;\dbR^{m\times
m})}$, we arrive at
$$
\begin{array}{ll}\ds
\q \mathbb{E}\int_{G}e^{(2r-r_0
)(t-s)}|y(t)|_{\dbR^m}^2dx -
\mathbb{E}\int_{G} |y(s)|_{\dbR^m}^2dx \\
\ns\ds  = -2\!\sum_{i\in\dbN\setminus
\L_r}\!\l_i\mathbb{E}\int_{s}^t\! e^{(2r-r_0
)(\si-s)} |y^i(\si)|_{\dbR^m}^2d\si \!  +\!
\mathbb{E}\int_{s}^t\int_{G}\!
e^{(2r-r_0 )(\si-s)}\big({\lan ay,y\ran}_{\dbR^m}\!\!+\!{\lan y,ay\ran}_{\dbR^m}\big) dxd\si\\
\ns\ds   \q + \mathbb{E}\int_{s}^t\int_{G}
e^{(2r-r_0 )(\si-s)}|b(\si)y(\si)|_{\dbR^m}^2 dxd\si+(2r-r_0 ) \mathbb{E}\int_{s}^t e^{(2r-r_0 )(\si-s)}  |y(\si)|_{\dbR^m}^2 dxd\si\\
  \leq 0,
\end{array}
$$
which gives the desired estimate (\ref{aaasss})
immediately.
\endpf

\subsubsection{Proof of the null
controllability result}

Now we are in a position to prove Theorem
\ref{thnull}.

\ms

{\it Proof of Theorem
\ref{thnull}}\,: We borrow some idea in
\cite{Lebeau-Robbiano}. For simplicity, we
assume that $m=1$. By Lemma {\ref{mea lemma}},
we may take a number $\tilde{t} \in E$ with
$\tilde{t} < T$ and a sequence $\{ t_{N} \}_{N =
1}^{\infty}\subset(0,T)$ such that
(\ref{t2.5})--(\ref{t2.7}) hold for some
positive numbers $\rho_1$ and $\rho_2$.

Write $\tilde{y}_0=\psi(t_1)$, where $\psi(\cd)$
solves the following stochastic  parabolic
system:
 \begin{eqnarray*}
\left\{
\begin{array}{lll}\ds
\ds d\psi - \sum_{j,k=1}^n (a^{jk}\psi_{x_j})_{x_k}dt = a(t)\psi dt+b(t)\psi dW(t)  \q & \mbox{ in } (0,t_1)\t G, \\
\ns\ds \psi= 0 & \mbox{ on } (0,t_1)\t \Gamma ,
\\
\ns\ds  \psi(0) =  y_0 &\mbox{ in } G.
\end{array}
\right.
 \end{eqnarray*}
Let us  consider the following controlled
stochastic parabolic system:
\begin{equation}\label{eq2.8}
\left\{
\begin{array}{lll}\ds
\ds d\tilde{y} - \sum_{j,k=1}^n
(a^{jk}\tilde{y}_{x_j})_{x_k}dt
= [a(t)\tilde{y} + \chi_{E}\chi_{G_0} \tilde{u}]dt+b(t)\tilde{y} dW(t)  \q& \mbox{ in } (t_1,\tilde{t}\,)\t G, \\
\ns\ds \tilde{y} = 0  & \mbox{ on }
(t_1,\tilde{t}\,)\t \Gamma ,
\\
\ns\ds  \tilde{y} (t_1) = \tilde{y}_{0} &\mbox{
in } G.
\end{array}
\right.
\end{equation}
It suffices to find a control $\tilde{u}\in
L^{\infty}_\dbF(t_{1},\tilde{t};L^2(\Omega;L^2(G)))$
with
\bel{20160504e5}
|\tilde{u}|^2_{L^{\infty}_\dbF(t_{1},\tilde{t};L^2(\Omega;L^2(G)))}
\leq \cC\mathbb{E}| \tilde{y}_{0}
|^2_{L^{2}(\Omega)}, \ee
such that the solution $\tilde{y}$ to
\eqref{eq2.8} satisfies $\tilde y(\tilde{t}\,) =
0$ in $G$, $\dbP$-a.s.

Set $I_{N} = [t_{2N-1 },t_{2N}]$ and $J_{N} =
[t_{2N},t_{2N+1}]$ for $N \in\dbN$. Then $
[t_{1},\tilde{t}\,) = \bigcup_{N =
1}^{\infty}(I_{N}\cup J_{N})$.
Clearly, $\dbm(E \cap I_{N}) > 0$ and $\dbm(E
\cap J_{N}) > 0$. We will introduce a suitable
control on each $I_N$ and allow the system to
evolve freely on every $J_N$. Also, we fix a
suitable, strictly increasing sequence
$\{r_N\}_{N=1}^\infty$ of positive integers (to
be given later) satisfying that $\l_1\leq
r_1<r_2<\cdots<r_N\to\infty$ as $N\to\infty$.

We consider first the controlled stochastic
parabolic system on the interval $I_{1} =
[t_{1},t_{2}]$ as follows:
\begin{equation}\label{system t1}
\left\{
\begin{array}{lll}\ds
\ds dy_1 - \sum_{j,k=1}^n
(a^{jk}y_{1,x_j})_{x_k}dt   = [a(t)y_1 +
\chi_{E}\chi_{G_0} u_1]dt+b(t) y_1 dW(t) \q &
\mbox{ in } (t_1,t_2)\t G,
\\
\ns\ds y_1 = 0 & \mbox{ on } (t_1,t_2)\t \Gamma
,
\\
\ns\ds  y_1 (t_1) =  \tilde{y}_{0} &\mbox{ in }
G.
\end{array}
\right.
\end{equation}
By Proposition \ref{conprop2}, there exists a
control $u_{1} \in
L^{\infty}_\dbF(t_{1},t_{2};L^2(\Omega;L^{2}(G)))$
with the estimate:
$$
|u_{1}|^2_{L^{\infty}_\dbF(t_{1},t_{2};L^2(\Omega;L^{2}(G)))
} \leq
{\frac{\cC_{1}e^{\cC_{2}\sqrt{r_1}+r_0T}}{(\dbm(E
\cap
[t_1,t_2]))^{2}}}\mathbb{E}|\tilde{y}_{0}|_{L^2(G)}^2,
$$
such that $\Pi_{r_1}(y(t_2)) = 0$ in $G$,
$\dbP$-a.s. By \eqref{t2.6}, we see that
\bel{20160501e3} |u_1
|^2_{L^{\infty}_\dbF(t_{1},t_{2};L^2(\Omega;L^{2}(G)))
} \leq
{\frac{\cC_{1}e^{\cC_{2}\sqrt{r_1}+r_0T}}{{\rho_1^2}(t_2-t_1)^2}}\mathbb{E}|\tilde{y}_{0}|_{L^2(G)}^2.
\end{equation}
Applying It\^o's formula to
$e^{-(r_0+1)t}|y_1(t)|^2_{L^2(G)} $, similar to
the proof of \eqref{conprop1eq21}, we obtain
that
\begin{equation*}
\begin{array}{ll}\ds
e^{-(r_0+1)t_2}\mathbb{E}|y_1(t_2)|^2_{L^2(G)} \\
\ns\ds =e^{-(r_0+1)t_1}  \mathbb{E}| y_1(t_1)
|^2_{L^2(G)} - (r_0+1)\mathbb{E}\int_{t_1}^{t_2}e^{-(r_0+1)s}\int_G |y_1|^2dxds\\
\ns\ds\q-
2\sum_{j,k=1}^n\mathbb{E}\int_{t_1}^{t_2}e^{-(r_0+1)s}
\int_G
  a^{jk}y_{1,x_j} y_{1,x_k}dxds  \\
\ns\ds\q + \mathbb{E}\!\int_{t_1}^{t_2}\!e^{-(r_0+1)s}\int_G[2a(s)|y_1|^2\!+\! |b(s)y_1|^2]dxds+  2\mathbb{E}\!\int_{t_1}^{t_2}\!e^{-(r_0+1)s}\int_G \!\chi_{E}\chi_{G_0}u_1 y_1 dxds \\
\ns\ds\leq e^{-(r_0+1)t_1}  \mathbb{E}| y_1(t_1)
|^2_{L^2(G)} +  \mathbb{E}\int_{t_1}^{t_2}e^{-(r_0+1)s}\int_G |u_1|^2  dxds\\
 \ns\ds \leq  e^{-(r_0+1)t_1}  \mathbb{E}| {\tilde{y}_0} |^2_{L^2(G)} +
{\frac{e^{-(r_0+1)t_1}-e^{-(r_0+1)t_2}}{r_0+1 }
|u_1
|^2_{L^{\infty}_\dbF(t_1,t_2;L^2(\Omega;L^{2}(G)))}}.
\end{array}
\end{equation*}
Hence, in view of \eqref{20160501e3},
 \bel{20160501e1}
 \mathbb{E}|y_1(t_2)|^2_{L^2(G)} \leq
\frac{\cC_3 e^{\cC_3{\sqrt{r_1}}}}{(t_2-t_1)^2}
\mathbb{E}| {\tilde{y}}_0 |^2_{L^2({G})}.
 \ee
where
$\cC_3=\max(2\rho_1^{-2}\cC_1e^{(2r_0+1)T},\cC_2)$.

Then, on the interval $J_1 \equiv [t_2,t_3],$ we
consider the following stochastic parabolic
system without control:
\begin{equation*}
\left\{
\begin{array}{lll}\ds
\ds dz_1 - \sum_{j,k=1}^n (a^{jk}z_{1,x_j})_{x_k}dt = a(t)z_1 dt+ b(t)z_1 dW(t)\q  & \mbox{ in } (t_2,t_3)\t G, \\
\ns\ds z_1 = 0 & \mbox{ on } (t_2,t_3)\t \Gamma
,
\\
\ns\ds  z_1 (t_2) =  y_1(t_2) &\mbox{ in } G.
\end{array}
\right.
\end{equation*}
Since $\Pi_{r_1}(y_1(t_2)) = 0$, $\dbP$-a.s., by
Proposition \ref{energydecay}, we have
\begin{equation}\label{zt1}
\begin{array}{ll}\ds
\mathbb{E}|z_1(t_3)|_{L^2(G)}^2  \leq
e^{(-2r_1+r_0 ) (t_3 - t_2)}\mathbb{E}| y_1(t_2)
|_{L^2(G)}^2
 \leq   {\frac{\cC_3
e^{\cC_3{\sqrt{r_1}}}}{(t_2-t_1)^2}}
e^{(-2r_1+r_0 )(t_3 -
t_2)}\mathbb{E}|\tilde{y}_0 |_{L^2(G)}^2.
\end{array}
\end{equation}

Generally, on the interval $I_N$ with
$N\in\dbN\setminus\{1\}$, we consider a
controlled stochastic parabolic system as
follows:
\begin{eqnarray*}
\left\{\!
\begin{array}{lll}\ds
\ds dy_N -\!\sum_{j,k=1}^n
(a^{jk}y_{N,x_j})_{x_k}dt \! = \![a(t)y_N \!+\!
\chi_{E}\chi_{G_0} u_N]dt\!+\!b(t) y_N dW(t) \q
& \mbox{ in } (t_{2N-1},t_{2N})\t G,
\\
\ns\ds  y_N = 0 & \mbox{ on }
(t_{2N-1},t_{2N})\t \Gamma ,
\\
\ns\ds  y_N (t_{2N-1}) =  z_{N-1}(t_{2N-1})
&\mbox{ in } G.
\end{array}
\right.
 \end{eqnarray*}
Similar to the above argument (See the proof of
\eqref{20160501e3} and \eqref{20160501e1}), one
can find a control $u_{N} \in
L^{\infty}_\dbF(t_{2N-1},t_{2N};L^2(\Omega;L^{2}(G)))$
with the estimate:
\begin{equation}\label{f22}\ds
|u_N
|^2_{L^{\infty}_\dbF(t_{2N-1},t_{2N};L^2(\Omega;L^{2}(G)))
} \leq
{\frac{\cC_{1}e^{\cC_{2}\sqrt{r_N}+r_0T}}{\rho_1^2(t_{2N}-t_{2N-1})^2}}\mathbb{E}|z_{N-1}(t_{2N-1})|_{L^2(G)}^2.
\end{equation}
such that $\Pi_{r_N}(y_N(t_{2N})) = 0$ in $G$,
$\dbP$-a.s. Moreover,
 \bel{20160501e2}
 \mathbb{E}|y_N(t_{2N})|^2_{L^2(G)} \leq
\frac{\cC_3
e^{\cC_3{\sqrt{r_N}}}}{(t_{2N}-t_{2N-1})^2}
\mathbb{E}| z_{N-1}(t_{2N-1}) |^2_{L^2({G})}.
 \ee
 On the interval $J_N,$ we consider the following stochastic parabolic system without
 control:
 \begin{eqnarray*}
\left\{
\begin{array}{lll}\ds
\ds dz_N - \sum_{j,k=1}^n
(a^{jk}z_{N,x_j})_{x_k}dt
= a(t)z_N dt+b(t)z_N dW(t)  \q& \mbox{ in } (t_{2N},t_{2N+1})\t G, \\
\ns\ds z_N = 0 & \mbox{ on } (t_{2N},t_{2N+1})\t
\Gamma ,
\\
\ns\ds  z_N (t_{2N}) =  y_N(t_{2N}) &\mbox{ in }
G.
\end{array}
\right.
 \end{eqnarray*}
Since $\Pi_{r_N}(y_N(t_{2N})) = 0$, $\dbP$-a.s.,
by Proposition \ref{energydecay} and similar to
\eqref{zt1}, and recalling that $y_N (t_{2N-1})
=  z_{N-1}(t_{2N-1})$ in $G$, we have
\begin{equation}\label{ztNNN}
\ba{ll}\ds \mathbb{E}|z_N(t_{2N+1})|_{L^2(G)}^2
\3n&\ds\leq   {\frac{\cC_3
e^{\cC_3{\sqrt{r_N}}}}{(t_{2N}-t_{2N-1})^2}}
e^{(-2r_N+r_0 )(t_{2N+1}-t_{2N})}\mathbb{E}|y_N
(t_{2N-1}) |_{L^2(G)}^2\\\ns&\ds\leq
{\frac{\cC_4
e^{\cC_4{\sqrt{r_N}}}}{(t_{2N}-t_{2N-1})^2}}
e^{-2r_N(t_{2N+1}-t_{2N})}\mathbb{E}|z_{N-1}(t_{2N-1})
|_{L^2(G)}^2, \ea
\end{equation}
where $\cC_4=\cC_3e^{r_0T}$.

Inductively, by \eqref{t2.7} and \eqref{ztNNN},
we conclude that, for all $N \geq 1$,
\begin{equation}\label{2016050e1}
\ba{ll}\ds
\mathbb{E}|z_N(t_{2N+1})|_{L^2(G)}^2 \\
\ns\ds\leq   \frac{\cC_4^N
e^{\cC_4(\sqrt{r_N}+\sqrt{r_{N-1}}+\cdots+\sqrt{r_1})}}
{(t_{2N}-t_{2N-1})^2(t_{2N-2}-t_{2N-3})^2\cdots(t_{2}-t_{1})^2}\\
\ns\ds\q\times
\mbox{exp}\Big\{-2r_N(t_{2N+1}-t_{2N})-2r_{N-1}(t_{2N-1}-t_{2N-2})
-\cdots-2r_1(t_{3}-t_{2})\Big\}\mathbb{E}|\tilde
y_0
|_{L^2(G)}^2\\
\ns\ds\leq   \frac{\cC_4^N
\exp\left\{\cC_4N\sqrt{r_N}-2r_N(t_{2N+1}-t_{2N})\right\}}
{(t_{2N}-t_{2N-1})^2(t_{2N-2}-t_{2N-3})^2\cdots(t_{2}-t_{1})^2}
\mathbb{E}|\tilde y_0
|_{L^2(G)}^2\\
\ns\ds\leq   \frac{\cC_4^N\rho_2^{2N(N-1)}
\exp\left\{\cC_4N\sqrt{r_N}-2(t_2-t_1)\rho_2^{1-2N}r_N\right\}}
{(t_{2}-t_{1})^{2N}} \mathbb{E}|\tilde y_0
|_{L^2(G)}^2. \ea
\end{equation}
By  \eqref{t2.7},
\eqref{f22}--\eqref{20160501e2} and
\eqref{2016050e1}, we see that
\begin{equation}\label{f22e1}
\ba{ll}\ds |u_N
|^2_{L^{\infty}_\dbF(t_{2N-1},t_{2N};L^2(\Omega;L^{2}(G)))
} \\
\ns\ds\leq
\frac{\cC_1\cC_4^{N-1}\rho_2^{2N(N-1)}}
{\rho_1^2(t_{2}-t_{1})^{2N}}
\exp\Big\{\cC_2\sqrt{r_N}+\!r_0T\!+\cC_4(N\!-\!1)\sqrt{r_{N-1}}\!-2(t_2\!-t_1)\rho_2^{3-2N}r_{N-1}\Big\}
\mathbb{E}|\tilde y_0 |_{L^2(G)}^2. \ea
\end{equation}
and
 \bel{20160501e2e1}
\ba{ll}\ds \mathbb{E}|y_N(t_{2N})|^2_{L^2(G)} \\
\ns\ds\leq
\frac{\cC_3\cC_4^{N-1}\rho_2^{2N(N-1)}}
{(t_{2}-t_{1})^{2N}}
\exp\Big\{\cC_3\sqrt{r_N}+\cC_4(N-1)\sqrt{r_{N-1}}-2(t_2-t_1)\rho_2^{3-2N}r_{N-1}\Big\}
\mathbb{E}|\tilde y_0 |_{L^2(G)}^2. \ea
\end{equation}
We now choose $r_N=\max(2^{N^2},[\l_1]+1)$. From
\eqref{f22e1}--\eqref{20160501e2e1}, it is easy
to see that, whenever $N$ is large enough,
 \bel{20160504e1}
|u_N
|^2_{L^{\infty}_\dbF(t_{2N-1},t_{2N};L^2(\Omega;L^{2}(G)))
}\le \frac{1}{2^N}\mathbb{E}|\tilde y_0
|_{L^2(G)}^2 \ee and
 \bel{20160504e13}
\mathbb{E}|y_N(t_{2N})|^2_{L^2(G)}\le
\frac{1}{2^N}\mathbb{E}|\tilde y_0 |_{L^2(G)}^2.
\ee

We now construct a control $\tilde{u}$ by
setting
\begin{eqnarray}
\tilde{u}(t,x) = \left\{
\begin{array}{lll}\ds
u_N(t,x),\ \ & (t,x) \in I_N\times G, \q N \geq
1,
\\
\ns\ds 0, & (t,x) \in J_N\times G, \q N \geq 1.
\end{array}
\right. \label{2.21}
\end{eqnarray}
By \eqref{20160504e1}, we see that $\tilde{u}\in
L^{\infty}_\dbF(t_1,\tilde{t}\,;L^2(\Omega;L^2(G)))$
satisfies \eqref{20160504e5}. Let $\tilde{y}$ be
the solution to the system \eqref{eq2.8}
corresponding to the control constructed in
\eqref{2.21}. Then $\tilde{y}(\cdot) =
y_N(\cdot)$ on $I_N\times G$. By
\eqref{20160504e13} and recalling that
$t_{2N}\to\tilde t$ as $N\to\infty$, we deduce
that $\tilde{y}(\tilde{t}\,) = 0$, $\dbP$-a.s.
This completes the proof of Theorem
\ref{thnull}.
\endpf

\subsubsection{Proof of the approximate
controllability result}

To begin with, we  show  the following two
preliminary results, which have some independent
interests.

\begin{proposition}\label{ob est}
If $\dbm((s,T)\bigcap E)>0$ for any $s\in
[0,T)$, then for any given $\eta\in
L^2_{\cF_T}(\Omega;$ $L^2(G;\dbR^m))$, the
corresponding solution to \eqref{sbheatsystem2}
with $s_1=0$ and $s_2=T$ satisfies
\begin{equation}\label{ok1}
|z(s)|_{L^2_{\cF_s}(\Omega;L^2(G;\dbR^m))}\leq
\cC(s)|\chi_E\chi_{G_0}
z|_{L^1_\dbF(s,T;L^2(\Omega;L^2(G;\dbR^m)))}.
\end{equation}
Here and henceforth, $\cC(s)>0$ is a generic
constant depending on $s$.
\end{proposition}

{\it Proof}\,:  We consider the following
controlled stochastic parabolic system:
\bel{asystem3} \left\{ \ba{ll} \ds dy -
\sum_{j,k=1}^n (a^{jk}y_{x_j})_{x_k}dt
= [a(t)y + \chi_{(s,T)\cap E}\chi_{G_0} u]dt+b(t) y dW(t) \q& \mbox{ in } (s,T)\times G , \\
\ns\ds y = 0  & \mbox{ on } (s,T)\times  \Gamma, \\
 \ns\ds y(s) = y_s &\mbox{ in } G,
\ea \right. \ee where $y$ is the state variable,
$u$ is the control variable, the initial state
$y_s\1n\in\1n L^2_{\cF_s}(\Omega;
L^2(G;\dbR^m))$ and the control $u(\cd)\in
L^\infty_\dbF(s,T;L^2(\Omega;L^2(G;\dbR^m)))$.
By the proof of Theorem \ref{thnull}, it is easy
to show that the system \eqref{asystem3} is null
controllable, i.e., for any $y_s\in
L^2_{\cF_s}(\Omega; L^2(G;\dbR^m))$, there
exists a control $u\in
L^\infty_\dbF(s,T;L^2(\Omega;L^2(G;\dbR^m)))$
such that $y(T)=0$ in $G$, $\dbP$-a.s., and
\begin{equation}\label{cost1}
|u|_{L^\infty_\dbF(s,T;L^2(\Omega;L^2(G;\dbR^m)))}^2
\leq \cC(s)|y_s|^2_{L^2_{\cF_s}(\Omega;
L^2(G;\dbR^m))}.
\end{equation}
Applying It\^{o}'s formula to ${\lan
y,z\ran}_{\dbR^m}$, where $y$ and $(z,Z)$ solve
respectively  \eqref{asystem3} and
\eqref{sbheatsystem2} with $s_1=0$ and $s_2=T$,
and noting that $y(T)=0$ in $G$, $\dbP$-a.s., we
obtain that
$$
- \mathbb{E}\int_G {\lan y_s,
z(s)\ran}_{\dbR^m}dx = \mathbb{E}\int_{(s,T)\cap
E}\int_{G_0}{\lan  u,z\ran}_{\dbR^m} dxdt.
$$
Choosing $y_s=-z(s)$ in \eqref{asystem3}, we
then have
$$
\begin{array}{ll}\ds
\mathbb{E}\int_G |z(s)|_{\dbR^m}^2 dx=
\mathbb{E}\int_{(s,T)\cap E}
\int_{G_0}{\lan u, z\ran}_{\dbR^m}dxdt\\
\ns\ds \leq
|u|_{L^\infty_\dbF(s,T;L^2(\Omega;L^2(G;\dbR^m)))}|\chi_{(s,T)\cap
E}\chi_{G_0}
z|_{L^1_\dbF(s,T;L^2(\Omega;L^2(G;\dbR^m)))} \\
\ns\ds \leq \cC(s)\Big( \mathbb{E}\int_G
|z(s)|_{\dbR^m}^2 dx
\Big)^{\frac{1}{2}}|\chi_{(s,T)\cap E}\chi_{G_0}
z|_{L^1_\dbF(s,T;L^2(\Omega;L^2(G;\dbR^m)))},
\end{array}
$$
which gives immediately the desired estimate
\eqref{ok1}.
\endpf

\ms

As an easy consequence of Proposition \ref{ob
est}, we have the following unique continuation
property for solutions to \eqref{sbheatsystem2}
with $s_1=0$ and $s_2=T$.

\begin{corollary}\label{heat1.1ucp}
If $\dbm((s,T)\cap E)>0$ for any $s\in [0,T)$,
then any solution $(z,Z)$ to
\eqref{sbheatsystem2} with $s_1=0$ and $s_2=T$
vanishes identically in $Q$, $\dbP$-a.s.
provided that $z=0$ in $G_0\t E$, $\dbP$-a.s.
\end{corollary}

{\it Proof}\,: Since $z=0$ in $G_0\t E$,
$\dbP$-a.s., by Proposition \ref{ob est}, we see
that $z(s)=0$ in $G$, $\dbP$-a.s., for any $
s\in [0,T)$. Therefore, $z\equiv 0$ in $Q$,
$\dbP$-a.s.
\endpf

\begin{remark}\label{rmkapp1}
If the condition $\dbm((s,T)\cap E)>0$ for any
$s\in [0,T)$ was not assumed, the conclusion in
Corollary \ref{heat1.1ucp} might fail to be
true. This can be shown by the following
counterexample. Let $E$ satisfy that $\dbm(E)>0$
and $\dbm((s_0,T)\cap E)=0$ for some $s_0\in
[0,T)$. Let $(z_1,Z_1) = 0$ in $(0,s_0)\times
G$, $\dbP$-a.s. and $\xi_2$ be a nonzero process
in $L^2_{\dbF}(s_0,T;\dbR^m)$ (Then $Z_2\equiv
\xi_2e_1$ is a nonzero process in
$L^2_{\dbF}(s_0,T;L^2(G;\dbR^m))$). Solving the
following forward stochastic differential
equation:
\begin{eqnarray*}
\left\{
\begin{array}{lll}\ds
d\zeta_1 -\l_1 \zeta_1 dt =- [a(t)^\top \zeta_1+b(t)^\top \xi_2]dt + \xi_2dW(t) & \mbox{ in } \,[s_0,T],\\
\ns\ds \zeta_1(s_0) = 0,
\end{array}
\right.
\end{eqnarray*}
we find a nonzero $\zeta_1\in
L_{\dbF}^2(\Omega;C([s_0,T];\dbR^m))$. In this
way, we find a nonzero solution $(z_2,Z_2)\equiv
(\zeta_1e_1,\xi_2e_1)\in
L_{\dbF}^2(\Omega;C([s_0,T];L^2(G;\dbR^m)))\times
L^2_{\dbF}(s_0,T;L^2(G;$ $\dbR^m))$  to the
following forward stochastic partial
differential equation:
\begin{equation}\label{asystem4}
\left\{
\begin{array}{lll}
\ds dz_2 + \sum_{j,k=1}^n
(a^{jk}z_{2,x_j})_{x_k}dt
= -[a(t)^\top z_2+b(t)^\top Z_2]dt  + Z_2dW(t) \quad & \mbox{ in }\, (s_0,T)\times G, \\
\ns\ds z_2 = 0  & \mbox{ on } (s_0,T)\times \Gamma , \\
 \ns\ds z_2(s_0)= 0  & \mbox{ in }\, G.
\end{array}
\right.
\end{equation}
(Note however that one cannot solve the system
\eqref{asystem4} directly because this system is
not well-posed). Put
\begin{eqnarray*}
(z,Z)=\left\{
\begin{array}{lll}\ds
(z_1,Z_1),\q & \mbox{ in } \; (0,s_0)\t G,\\
\ns\ds (z_2,Z_2), & \mbox{ in } \; (s_0,T)\t G.
\end{array}
\right.
\end{eqnarray*}
Then, $(z,Z)$ is a nonzero solution to
\eqref{sbheatsystem2} with $s_1=0$ and $s_2=T$,
and $z=0$ in $G_0\t E$, $\dbP$-a.s. Note also
that, the nonzero solution constructed for the
system (\ref{asystem4}) indicates that, in
general, the forward uniqueness does NOT hold
for backward stochastic differential equations.
\end{remark}

We are now in a position to prove Theorem
\ref{thapp}.

\ms

{\it Proof of Theorem \ref{thapp}}\,: The ``if"
part follows from Corollary \ref{heat1.1ucp}. To
prove the ``only if" part, we use the
contradiction argument. Assume that
$\dbm((s_0,T)\cap E)=0$ for some $s_0\in [0,T)$.
Since the system (\ref{sheatsystem1}) is
approximately controllable at time $T$,  we
deduce that any solution $(z,Z)$ to
\eqref{sbheatsystem2} with $s_1=0$ and $s_2=T$
vanishes identically in $Q$ provided that $z=0$
in $G_0\t E$, $\dbP$-a.s. This contradicts the
counterexample in Remark \ref{rmkapp1}.
\endpf


\subsection{Null controllability of stochastic parabolic systems}
\label{sto heat se2}


In this subsection, we deal with the null
controllability for \eqref{heat 1.1}. The results in this subsection are taken from
\cite{Tang-Zhang1}.

We have the following result.

\begin{theorem}\label{f heat control1}
Let the condition \eqref{20160618e1} be
satisfied.   Then the system \eqref{heat 1.1} is
null controllable at time $T$.
\end{theorem}

In order to prove Theorem \ref{f heat control1}, by means of the standard duality argument, it suffices to establish
the following observability result for
\eqref{dual heat 1.1}:
\begin{theorem}\label{b heat ob1}
Let the condition \eqref{20160618e1} be
satisfied. Then,  for all $z_T\in
L^2_{\cF_T}(\Omega;L^2(G;\dbR^m))$, solutions
$(z,Z)\in L^2_{\dbF}(\Omega;C([0,T];L^2(G;$
$\dbR^m)))\times L^2_{\dbF}(0,T;L^2(G;\dbR^m))$
to the system \eqref{dual heat 1.1} satisfy
\begin{equation}\label{b heat obser1}
\begin{array}{ll}\ds
|z(0)|_{L_{\cF_0}^2(\Omega;L^2(G;\dbR^m))}\le
\ds \cC \big(|z|_{L^2_{\dbF}(0,T;L^2(G_0;\dbR^m))}+|Z|_{L^2_{\dbF}(0,T;L^2(G;\dbR^m))}\big).
\end{array}
\end{equation}
\end{theorem}

\br
In Theorem \ref{b heat ob1}, we assume that
$a_{1j} \in
L^{\infty}_{\dbF}(0,T;W^{1,\infty}(G;$
$\dbR^{m\times m}))$ for $j=1,2,\cdots,n$ (See
the condition \eqref{20160618e1}). It seems that
this assumption can be weakened as $a_{1j} \in
L^{\infty}_{\dbF}(0,T;L^{\infty}(G;$
$\dbR^{m\times m}))$.
\er

The rest of this subsection is devoted to giving a
proof of Theorem \ref{b heat ob1}. For
simplicity, we consider only the case $m=1$.

\subsubsection{A weighted identity and Carleman
estimate  for a stochastic parabolic-like
operator}

In order to prove Theorem \ref{b heat ob1}, we need to derive a
weighted identity and Carleman estimate  for  a
stochastic parabolic-like operator.

We assume that
\begin{equation}\label{heat 1.1 1}
b^{jk}=b^{kj}\in
L_{\dbF}^2(\Omega;C^1([0,T];W^{2,\infty}(G))),\qq
j,k=1,2,\cdots,n,
\end{equation}
$\ell\in C^{1,3}(Q)$ and $\Psi\in C^{1,2}(Q)$.
Write
\begin{equation}\label{c1e15}
\left\{
 \begin{array}{ll}
 \ds \cA =-\sum_{j,k=1}^n
 \big(b^{jk}\ell_{x_j}\ell_{x_k}-b^{jk}_{x_k}\ell_{x_j}
 -b^{jk}\ell_{x_jx_k}\big)-\Psi-\ell_t,\\
  \ns
 \ds
 \cB=2\[\cA \Psi-
 \sum_{j,k=1}^n\big(\cA b^{jk}\ell_{x_j}\big)_{x_k}\] -\cA _t-\sum_{j,k=1}^n (b^{jk}\Psi_{x_k})_{x_j},\\
  \ns
 \ds
 c^{jk}=\sum_{j',k'=1}^n\[2b^{jk'}\big(b^{j'k}\ell_{x_{j'}}\big)_{x_{k'}}
-
\big(b^{jk}b^{j'k'}\ell_{x_{j'}}\big)_{x_{k'}}\]-\frac{b_t^{jk}}{2}+\Psi
b^{jk}.
  \end{array}
\right.
\end{equation}

First, we establish a fundamental weighted
identity for the stochastic parabolic-like
operator
``$dh-\sum_{j,k=1}^n(b^{jk}h_{x_j})_{x_k}dt$"\footnote{Since
only the symmetry condition \eqref{heat 1.1 1}
is assumed for the coefficient matrix
$\big(b^{jk}\big)$, we call ``$
dh-\sum_{j,k=1}^n\big(b^{jk}h_{x_j}\big)_{x_k}dt$"
a stochastic parabolic-like operator.}.

\begin{theorem}\label{c1t1}
Let $h$ be an $H^2(G)$-valued continuous
semi-martingale. Set $\th=e^{\ell }$ and $w=\th
h$. Then, for any $t\in [0,T]$ and a.e.
$(x,\om)\in G\times\Omega$,
\begin{equation}\label{c1e2a}
\begin{array}{ll} \displaystyle
2 \th\[-\sum_{j,k=1}^n (b^{jk}w_{x_j})_{x_k}+\cA
w\]\[dh-\sum_{j,k=1}^n(b^{jk}h_{x_j})_{x_k}dt\]+2
\sum_{j,k=1}^n
(b^{jk}w_{x_j}dw)_{x_k}\\
\ns\ds\q\!\!\!+ 2\!
\sum_{j,k=1}^n\!\!\[\!\sum_{j',k'=1}^n\!\!\!\big(2b^{jk}
b^{j'k'}\ell_{x_{j'}}w_{x_j}w_{x_{k'}} \!\!
-\!b^{jk}b^{j'k'}\ell_{x_j}w_{x_{j'}}w_{x_{k'}}\!\big)
\!+\!\Psi
 b^{jk}w_{x_j}w\!-\! b^{jk}\(\cA \ell_{x_j}\!\!+\!\frac{\Psi_{x_j}}{2}\)w^2\]_{x_k}\!dt\\
\ns\ds =2 \!\sum_{j,k=1}^n
c^{jk}w_{x_j}w_{x_k}dt \!+\! \cB
w^2dt+d\(\sum_{j,k=1}^n
\!b^{jk}w_{x_j}w_{x_k}+\cA w^2\) \! +2 \[\!-\sum_{j,k=1}^n \big(b^{jk}w_{x_j}\big)_{x_k}+\cA w\]^2dt\\
\ns\ds\q - \th^2\sum_{j,k=1}^n
 b^{jk}( dh_{x_j} + \ell_{x_j}dh)(dh_{x_k}+\ell_{x_k}dh)- \th^2\cA (dh)^2.
 \end{array}
 \end{equation}
\end{theorem}

{\it Proof}\,: The proof is divided into four
steps.

\ms

{\bf Step 1.} Recalling that $\th=e^\ell$ and $w=\th
h$, one has $dh=\th^{-1}(dw-\ell_twdt)$ and
$h_{x_j}=\th^{-1}(w_{x_j}-\ell_{x_j}w)$ for
$i=1,2,\cds,m$. By \eqref{heat 1.1 1}, it is
easy to see that
$\sum_{j,k=1}^nb^{jk}(\ell_{x_j}w_{x_k}+\ell_{x_k}w_{x_j})=2\sum_{j,k=1}^nb^{jk}\ell_{x_j}w_{x_k}$.
Hence,
\begin{equation}\label{t2}
\begin{array}{ll}
\ds \th
\sum_{j,k=1}^n(b^{jk}h_{x_j})_{x_k}=\th\sum_{j,k=1}^n[\th^{-1}b^{jk}(w_{x_j}-\ell_{x_j}
w)]_{x_k}
\\ \ns \ds  =\sum_{j,k=1}^n[b^{jk}(w_{x_j}-\ell_{x_j} w)]_{x_k}-\sum_{j,k=1}^nb^{jk}(w_{x_j}-\ell_{x_j} w)\ell_{x_k}\\
\ns \ds
=\sum_{j,k=1}^n\[(b^{jk}w_{x_j})_{x_k}-b^{jk}(\ell_{x_j}w_{x_k}+\ell_{x_k}w_{x_j})
 +(b^{jk}\ell_{x_j}\ell_{x_k}-b^{jk}_{x_k}\ell_{x_j}
-b^{jk}\ell_{x_jx_k})w\]\\
\ns \ds
=\sum_{j,k=1}^n\[(b^{jk}w_{x_j})_{x_k}-2b^{jk}\ell_{x_j}w_{x_k}
+(b^{jk}\ell_{x_j}\ell_{x_k}-b^{jk}_{x_k}\ell_{x_j}
-b^{jk}\ell_{x_jx_k})w\].
\end{array}
\end{equation}

Put
\begin{equation}\label{t4}
\left\{
\begin{array}{ll}
\ds I\=-\sum_{j,k=1}^n (b^{jk}w_{x_j})_{x_k}+\cA w,\q\; I_1\=\[-\sum_{j,k=1}^n (b^{jk}w_{x_j})_{x_k}+\cA w\]dt,\\
 \nm
\ds
I_2\=dw+2\sum_{j,k=1}^nb^{jk}\ell_{x_j}w_{x_k}dt,
 \qq I_3\=\Psi wdt.
\end{array}
\right.
\end{equation}
By (\ref{t2}) and (\ref{t4}), it follows that
$$
\th\[dh-\sum_{j,k=1}^n(b^{jk}h_{x_j})_{x_k}\,
dt\]=I_1+I_2+I_3.
$$
Hence,
\begin{equation}\label{heat 1.1 c1e5}
\begin{array}{ll}
\displaystyle 2 \th\[-\sum_{j,k=1}^n
(b^{jk}w_{x_j})_{x_k}+\cA
w\]\[dh-\sum_{j,k=1}^n(b^{jk}h_{x_j})_{x_k}dt\]
=2 I(I_1+I_2+I_3).
\end{array}
\end{equation}
\ms

{\bf Step 2.} Let us compute $2 II_2$. Utilizing
(\ref{heat 1.1 1}) again, and noting that
$$
 \sum_{j,k,j',k'=1}^n
 (b^{jk}b^{j'k'}\ell_{x_{j'}}w_{x_j}w_{x_k})_{x_{k'}}= \sum_{j,k,j',k'=1}^n
 (b^{jk}b^{j'k'}\ell_{x_j}w_{x_{j'}}w_{x_{k'}})_{x_k},$$
we get
\begin{equation}\label{2c2t8}
\begin{array}{ll}
\ds 2 \sum_{j,k,j',k'=1}^n  b^{jk}b^{j'k'}\ell_{x_{j'}}w_{x_j}w_{x_kx_{k'}}\\
\ns \ds = \sum_{j,k,j',k'=1}^n
b^{jk}b^{j'k'}\ell_{x_{j'}}(w_{x_j}w_{x_kx_{k'}}+w_{x_k}w_{x_jx_{k'}})
 = \sum_{j,k,j',k'=1}^n
b^{jk}b^{j'k'}\ell_{x_{j'}}(w_{x_j}w_{x_k})_{x_{k'}}\\
\noalign{\ss} \displaystyle =
\sum_{j,k,j',k'=1}^n
(b^{jk}b^{j'k'}\ell_{x_{j'}}w_{x_j}w_{x_k})_{x_{k'}}-
\sum_{j,k,j',k'=1}^n
(b^{jk}b^{j'k'}\ell_{x_{j'}})_{x_{k'}}w_{x_j}w_{x_k}\\
\noalign{\ss} \displaystyle =
\sum_{j,k,j',k'=1}^n
(b^{jk}b^{j'k'}\ell_{x_j}w_{x_{j'}}w_{x_{k'}})_{x_k}-
\sum_{j,k,j',k'=1}^n
(b^{jk}b^{j'k'}\ell_{x_{j'}})_{x_{k'}}w_{x_j}w_{x_k}.
\end{array}
\end{equation}
Hence, by \eqref{2c2t8}, and noting that
$$
 \sum_{j,k,j',k'=1}^n
b^{jk}(b^{j'k'}\ell_{x_{j'}})_{x_k}
w_{x_j}w_{x_{k'}}=  \sum_{j,k,j',k'=1}^n
b^{jk'}(b^{j'k}\ell_{x_{j'}})_{x_{k'}}
w_{x_j}w_{x_k},
$$
we obtain that
\begin{equation}\label{heat1.1 1c1.3}
\begin{array}{ll} \displaystyle
 4\[-\sum_{j,k=1}^n (b^{jk}w_{x_j})_{x_k}+\cA w\]\sum_{j,k=1}^nb^{jk}\ell_{x_j}w_{x_k}\\
\ns \ds =- 4 \sum_{j,k,j',k'=1}^n (b^{jk}
b^{j'k'}\ell_{x_{j'}}w_{x_j}w_{x_{k'}})_{x_k}
+4 \sum_{j,k,j',k'=1}^n  b^{jk}(b^{j'k'}\ell_{x_{j'}})_{x_k} w_{x_j}w_{x_{k'}}\\
\ns \ds \q+4 \sum_{j,k,j',k'=1}^n b^{jk}
b^{j'k'}\ell_{x_{j'}}w_{x_j}w_{x_kx_{k'}}+
2\cA\sum_{j,k=1}^nb^{jk}\ell_{x_j}(w^2)_{x_k}\\
\ns \ds = -
2\sum_{j,k=1}^n\[\sum_{j',k'=1}^n\big(2b^{jk}
b^{j'k'}\ell_{x_{j'}}w_{x_j}w_{x_{k'}}
-b^{jk}b^{j'k'}\ell_{x_j}w_{x_{j'}}w_{x_{k'}}\big)-\cA b^{jk}\ell_{x_j} w^2\]_{x_k}\\
\ns \ds \q+2 \sum_{j,k,j',k'=1}^n
\[2b^{jk'}(b^{j'k}\ell_{x_{j'}})_{x_{k'}} -
(b^{jk}b^{j'k'}\ell_{x_{j'}})_{x_{k'}}\]w_{x_j}w_{x_k}-
2\sum_{j,k=1}^n(\cA b^{jk}\ell_{x_j})_{x_k} w^2.
\end{array}
\end{equation}
Using It\^o's formula, we have
\begin{equation}\label{2c1}
\begin{array}{ll}
\ds
2 \[-\sum_{j,k=1}^n (b^{jk}w_{x_j})_{x_k}+\cA w\]dw\\
\ns \ds =-2
\sum_{j,k=1}^n(b^{jk}w_{x_j}dw)_{x_k}+2
\sum_{j,k=1}^n b^{jk}w_{x_j}dw_{x_k}+2 \cA wdw
\\
\ns \ds =-2 \sum_{j,k=1}^n
(b^{jk}w_{x_j}dw)_{x_k}+d\(\sum_{j,k=1}^n
b^{jk}w_{x_j}w_{x_k}+\cA w^2\) \\
\ns \ds \q- \sum_{j,k=1}^n
b_t^{jk}w_{x_j}w_{x_k}dt - \cA _tw^2dt -
\sum_{j,k=1}^n b^{jk}dw_{x_j}dw_{x_k}- \cA
(dw)^2.
\end{array}
\end{equation}
Now, from (\ref{t4}), (\ref{heat1.1 1c1.3}) and
(\ref{2c1}), we arrive at
\begin{equation}\label{2c11}
\begin{array}{ll}
\ds 2 II_2\3n&\ds= - 2
\sum_{j,k=1}^n\[\sum_{j',k'=1}^n\big(2b^{jk}
b^{j'k'}\ell_{x_{j'}}w_{x_j}w_{x_{k'}}
-b^{jk}b^{j'k'}\ell_{x_j}w_{x_{j'}}w_{x_{k'}}\big)-\cA b^{jk}\ell_{x_j} w^2\]_{x_k}dt\\
\ns&\ds \q  -2 \sum_{j,k=1}^n
(b^{jk}w_{x_j}dw)_{x_k}+d\(\sum_{j,k=1}^n
b^{jk}w_{x_j}w_{x_k}+\cA w^2\)\\
\ns &\ds \q+2 \sum_{j,k=1}^n
\[\sum_{j',k'=1}^n\left(2b^{jk'}(b^{j'k}\ell_{x_{j'}})_{x_{k'}}
-(b^{jk}b^{j'k'}\ell_{x_{j'}})_{x_{k'}}\right)-\frac{ b_t^{jk}}{2}\]w_{x_j}w_{x_k}dt\\
\ns &\ds\q-
\[\cA_t+
2\sum_{j,k=1}^n(\cA b^{jk}\ell_{x_j})_{x_k}\]
w^2dt- \sum_{j,k=1}^n
 b^{jk}dw_{x_j}dw_{x_k}- \cA
(dw)^2.
\end{array}
\end{equation}

{\bf Step 3.} Let us compute $2 II_3$. By
(\ref{t4}), we get
\begin{equation}\label{heat 1.1 2c2t11}
\begin{array}{ll}
\displaystyle 2 II_3 =2 \[-\sum_{j,k=1}^n
(b^{jk}w_{x_j})_{x_k}+\cA w\]\Psi wdt\\
\ns\ds= \[ -2\sum_{j,k=1}^n \big(\Psi
b^{jk}w_{x_j}w\big)_{x_k}+ 2\Psi \sum_{j,k=1}^n
b^{jk}w_{x_j}w_{x_k} +\sum_{j,k=1}^n b^{jk}\Psi_{x_k}(w^2)_{x_j}+ 2\cA \Psi w^2\]dt\\
\ns\ds=  \Big\{-\sum_{j,k=1}^n \(2\Psi
b^{jk}w_{x_j}w - b^{jk}\Psi_{x_j} w^2\)_{x_k} +
2\Psi
 \sum_{j,k=1}^n
b^{jk}w_{x_j}w_{x_k}\\
\ns\ds\q+\[-\sum_{j,k=1}^n
(b^{jk}\Psi_{x_k})_{x_j}+2\cA \Psi\]
w^2\Big\}dt.
\end{array}
\end{equation}
\ms

{\bf Step 4.} Finally, combining the equalities
(\ref{heat 1.1 c1e5}), (\ref{2c11}) and
(\ref{heat 1.1 2c2t11}), and noting that
$$
\begin{array}{ll}
\ds\sum_{j,k=1}^n b^{jk} dw_{x_j}dw_{x_k}+ \cA
(dw)^2 = \th^2\sum_{j,k=1}^n b^{jk}
(dh_{x_j}+\ell_{x_j}dh)(dh_{x_k}+\ell_{x_k}dh)+
\th^2\cA (dh)^2,
\end{array}
$$
we conclude the desired equality (\ref{c1e2a})
immediately. \endpf

\medskip

Next, we shall derive a Carleman estimate for
the stochastic parabolic-like operator
``$dh-\sum_{j,k=1}^n(b^{jk}h_{x_j})_{x_k}dt$",

For any fixed nonnegative and nonzero function
$\psi\in C^4(\cl{G})$, and (large) parameters
$\l>1$ and $\mu>1$, we choose
 \be
 \label{alphad}
 \th=e^{\ell }, \q\ell=\l\a,\q\a(t,x)=\frac{e^{\mu\psi(x)}-e^{2\mu|\psi|_{C(\cl{G })}}}{ t(T-t)},\q
 \varphi(t,x)=\frac{e^{\mu\psi(x)}}{t(T-t)},
 \end{equation}
and
\begin{equation}\label{h5}
\Psi=2\sum_{j,k=1}^nb^{jk}\ell_{x_jx_k}.
\end{equation}
In what follows, for a positive integer $r$, we
denote by $O(\mu^r)$ a function of order $\mu^r$
for large $\mu$ (which is independent of $\l$);
by $O_{\mu}(\l^r)$ a function of order $\l^r$
for fixed $\mu$ and for large $\l$. In a similar
way, we use the notation
$O(e^{\mu|\psi|_{C(\cl{G })}})$ and so on. For
$j,k=1,2,\cdots,n$, it is easy to check that
 \begin{equation}\label{heat 1.1 h2}
 \ell_t=\l\a_t,\q \ell_{x_j}=\l\mu\f\psi_{x_j},\q
 \ell_{x_jx_k}=\l\mu^2\f\psi_{x_j}\psi_{x_k}+\l\mu\f\psi_{x_jx_k}
 \end{equation}
and that
 \begin{equation}\label{h4}
\a_t=\f^2O(e^{2\mu |\psi|_{C(\cl{G})}}),\qq \f_t
=\f^2O(e^{\mu |\psi|_{C(\cl{G})}}).
 \end{equation}

We have the following result.

\begin{theorem}\label{c1t2}
Assume that either $(b^{jk})_{n\t n}$ or
$-(b^{jk})_{n\t n}$ is a uniformly positive
definite matrix, and  its smallest eigenvalue is
bigger than a constant $s_0>0$. Let $h$ and $
w=\th h$ be that in Theorem \ref{c1t1} with
$\th$ being given in \eqref{alphad}. Then, the
equality \eqref{c1e2a} holds for any $t\in [0,
T]$ and a.e. $(x,\om)\in G\times \Omega$.
Moreover, for $\cA$, $\cB$ and $c^{jk}$ appeared
in \eqref{c1e2a} (and given by \eqref{c1e15}),
when $|\n \psi(x)|>0$,  $\l$ and $\mu$ are large
enough, it holds that
 \begin{equation}\label{h6}
 \left\{
 \begin{array}{ll}\ds
\cA
=-\l^2\mu^2\f^2\sum_{j,k=1}^nb^{jk}\psi_{x_j}\psi_{x_k}
  +\l\f^2O(e^{2\mu
|\psi|_{C(\cl{G})}}), \\
 \ns
 \ds \cB\ge \ds 2s_0^2\l^3\mu^4\f ^3|\n \psi|^4
  +\l^3\f^3O(\mu^3)
   +\l^2\f^3 O(\mu^2e^{2\mu
|\psi|_{C(\cl{G})}})+\l\f^3O(e^{2\mu
|\psi|_{C(\cl{G})}}), \\
 \ns
 \ds
 \sum_{j,k=1}^n c^{jk}w_{x_j}w_{x_k}\ge [s_0^2\l\mu^2\f |\n\psi|^2+\l\f  O(\mu)]|\n
w|^2
  \end{array}\right.
 \end{equation}
for any $t\in [0,T]$, $\dbP$-a.s.
\end{theorem}

{\it Proof}\,: By Theorem \ref{c1t1}, it remains
to prove the estimates in \eqref{h6}.

Noting (\ref{h5})--(\ref{heat 1.1 h2}), from
(\ref{c1e15}), we have $\ell_{x_jx_k}=\l\mu^2\f
\psi_{x_j}\psi_{x_k}+\l\f O(\mu)$ and that
$$
\begin{array}{ll}
\ds \sum_{j,k=1}^n c^{jk}w_{x_j}w_{x_k}\\
\ns\ds=\sum_{j,k=1}^n\Big\{\sum_{j',k'=1}^n\!
\[2b^{jk'}b^{j'k}\ell_{x_{j'}x_{k'}}\! +\!
b^{jk}b^{j'k'}\ell_{x_{j'}x_{k'}}\!+\!
2b^{jk'}b_{x_{k'}}^{j'k}\ell_{x_{j'}} \!-\!
(b^{jk}b^{j'k'})_{x_{k'}}\ell_{x_{j'}}\] \!- \!\frac{b_t^{jk}}{2} \Big\}w_{x_j}w_{x_k}\\
\ns \ds
=\sum_{j,k=1}^n\!\Big\{\sum_{j',k'=1}^n\!\[2\l\mu^2\f
b^{jk'}b^{j'k}\psi_{x_{j'}}\psi_{x_{k'}}
\!+\!\l\mu^2\f
b^{jk}b^{j'k'}\psi_{x_{j'}}\psi_{x_{k'}} \!+\!\l\f  O(\mu)\]\Big\}w_{x_j}w_{x_k}\\
\ns \ds=2\l\mu^2\f
\(\sum_{j,k=1}^nb^{jk}\psi_{x_j}w_{x_k}\)^2+\l\mu^2\f
\(\sum_{j,k=1}^n
b^{jk}\psi_{x_j}\psi_{x_k}\)\(\sum_{j,k=1}^nb^{jk}w_{x_j}w_{x_k}\)+\l\f |\n w|^2 O(\mu)\\
\ns \ds\ge [s_0^2\l\mu^2\f |\n\psi|^2+\l\f
O(\mu)]|\n w|^2,
\end{array}
$$
which gives the last inequality in (\ref{h6}).

Similarly, by the definition of $\cA$ in
\eqref{c1e15}, and noting \eqref{h4},  we see
that
$$
\begin{array}{ll}
\ds \cA =-\sum_{j,k=1}^n \big(b^{jk}\ell_{x_j}\ell_{x_k}-b_{x_k}^{jk}\ell_{x_j}+b^{jk}\ell_{x_jx_k}\big)-\ell_t\\
\ns \ds \q
=-\l\mu\sum_{j,k=1}^n\[b^{jk}\l\mu\f^2\psi_{x_j}\psi_{x_k}-b_{x_k}^{jk}\f\psi_{x_j}+b^{jk}\big(\mu\f\psi_{x_j}\psi_{x_k}
+\f\psi_{x_jx_k}\big)\] +\l\f^2O(e^{2\mu
|\psi|_{C(\cl{G})}})\\
\ns \ds \q
=-\l^2\mu^2\f^2\sum_{j,k=1}^nb^{jk}\psi_{x_j}\psi_{x_k}+\l\f^2O(e^{2\mu
|\psi|_{C(\cl{G})}}).
\end{array}
$$
Hence, we get the first estimate in \eqref{h6}.

Now, let us estimate $\cB$ (recall \eqref{c1e15}
for the definition of $\cB$). For this, by
\eqref{heat 1.1 h2}, and recalling the
definitions of $\Psi$ (in \eqref{h5}), we see
that
$$
\ds\Psi =2\l\mu\sum_{j,k=1}^nb^{jk}(\mu\f
\psi_{x_j}\psi_{x_k}+\f
\psi_{x_jx_k})=2\l\mu^2\f
\sum_{j,k=1}^nb^{jk}\psi_{x_j}\psi_{x_k}+\l\f
O(\mu);
$$
$$
\ell_{x_{j'}x_{k'}x_k}=\l\mu^3\f\psi_{x_{j'}}\psi_{x_{k'}}\psi_{x_k}+\l\f
O(\mu^2),
$$
$$
\begin{array}{ll}\ds
\ell_{x_{j'}x_{k'}x_jx_k}=\l\mu^4\f\psi_{x_{j'}}\psi_{x_{k'}}\psi_{x_j}\psi_{x_k}+\l\f O(\mu^3),\\
\ns \ds
\Psi_{x_k}=2\sum_{j',k'=1}^n\big(b^{j'k'}\ell_{x_{j'}x_{k'}}\big)_{x_k}=2\sum_{j',k'=1}^n(b^{j'k'}_{x_k}\ell_{x_{j'}x_{k'}}+b^{j'k'}\ell_{x_{j'}x_{k'}x_k})
\\\ns\ds \q\;\;\;=2\l\mu^3\f\sum_{j',k'=1}^nb^{j'k'}\psi_{x_{j'}}\psi_{x_{k'}}\psi_{x_k}+\l\f O(\mu^2),\\
\end{array}
$$
$$
\begin{array}{ll}\ds
\Psi_{x_jx_k}=2\sum_{j',k'=1}^n\big(b^{j'k'}_{x_jx_k}\ell_{x_{j'}x_{k'}}+b^{j'k'}\ell_{x_{j'}x_{k'}x_jx_k}+2b^{j'k'}_{x_k}\ell_{x_{j'}x_{k'}x_j}\big)
\\ \ns\ds \q\;\;\;\;\;\;=2\l\mu^4\f\sum_{j',k'=1}^nb^{j'k'}\psi_{x_{j'}}\psi_{x_{k'}}\psi_{x_j}\psi_{x_k}+\l\f O(\mu^3),
\end{array}
$$
$$
-\sum_{j,k=1}^n\!\big(b^{jk}\Psi_{x_k}\big)_{x_j}
\!=\!-\sum_{j,k=1}^n\!\big(b^{jk}_{x_j}\Psi_{x_k}+b^{jk}\Psi_{x_jx_k}\big)
\!=\!-2\l\mu^4\f\(\sum_{j,k=1}^n\!\!b^{jk}\psi_{x_j}\psi_{x_k}\)^2
\!+\!\l\f O(\mu^3).
$$
Hence, recalling  the definition of $\cA$ (in
\eqref{c1e15}), and using \eqref{heat 1.1 h2}
and \eqref{h4}, we have that
$$
\cA\Psi=-2\l^3\mu^4\f^3\(\sum_{j,k=1}^n
b^{jk}\psi_{x_j}\psi_{x_k}\)^2+\l^3\f^3O(\mu^3)+\l^2\f^3O(\mu^2e^{2\mu
|\psi|_{C(\cl{G})}}),
$$
$$
\begin{array}{ll}
\ds
\cA_{x_k}=-\sum_{j',k'=1}^n\big(b_{x_k}^{j'k'}\ell_{x_{j'}}\ell_{x_{k'}}
+2b^{j'k'}\ell_{x_{j'}}\ell_{x_{k'}x_k}-b_{x_{k'}x_k}^{j'k'}\ell_{x_{j'}}
 - b_{x_{k'}}^{j'k'}\ell_{x_{j'}x_k}
\\
\ns\ds \qq\q  + b_{x_k}^{j'k'}\ell_{x_{j'}x_{k'}} + b^{j'k'}\ell_{x_{j'}x_{k'}x_k}\big)-\ell_{tx_k}\\
\ns \ds \;\q\,\;\;
=-\sum_{j',k'=1}^n\big(2b^{j'k'}\ell_{x_{j'}}\ell_{x_{k'}x_k}+b^{j'k'}\ell_{x_{j'}x_{k'}x_k}\big)-\ell_{tx_k}+\big(\l\f+\l^2\f^2
\big)O(\mu^2)\\
\ns \ds \q\;\,\;\;
=-2\l^2\mu^3\f^2\sum_{j',k'=1}^nb^{j'k'}\psi_{x_{j'}}\psi_{x_{k'}}\psi_{x_k}+\l^2\f^2
O(\mu^2)+\l\f^2O(\mu e^{2\mu
|\psi|_{C(\cl{G})}}),
\end{array}
$$
$$
\sum_{j,k=1}^n\cA_{x_k}b^{jk}\ell_{x_j}=-2\l^3\mu^4\f^3
\(\sum_{j,k=1}^nb^{jk}\psi_{x_j}\psi_{x_k}\)^2
+\l^3\f^3 O(\mu^3)+\l^2\f^3O(\mu^2e^{2\mu
|\psi|_{C(\cl{G})}}),
$$
$$
\begin{array}{ll}
\ds
\sum_{j,k=1}^n\big(\cA b^{jk}\ell_{x_j}\big)_{x_k}=\sum_{j,k=1}^n\cA_{x_k}b^{jk}\ell_{x_j}+\cA\sum_{j,k=1}^n\big(b_{x_k}^{jk}\ell_{x_j}+b^{jk}\ell_{x_jx_k}\big)\\
\ns \ds \qq=-3\l^3\mu^4\f^3
\(\sum_{j,k=1}^nb^{jk}\psi_{x_j}\psi_{x_k}\)^2
+\l^3\f^3 O(\mu^3)+\l^2\f^3O(\mu^2e^{2\mu
|\psi|_{C(\cl{G})}}),
\end{array}
$$
and that
 $$
  \begin{array}{ll}
  \ds
  \ds \cA_t=-\sum_{j,k=1}^n
 \(b^{jk}\ell_{x_j}\ell_{x_k}-b^{jk}_{x_k}\ell_{x_j}
 +b^{jk}\ell_{x_jx_k}-\ell_t\)_t\\
  \ns
 \ds\;\q\, =-\sum_{j,k=1}^n
 \[b^{jk}\big(\ell_{x_j}\ell_{x_k}\big)_t-b^{jk}_{x_k}\ell_{x_jt}
 +b^{jk}\ell_{x_jx_kt}\] +\l^2\f^2 O(\mu^2)+\l\f^3O(e^{2\mu
|\psi|_{C(\cl{G})}})\\
  \ns
 \ds \;\q\,=\l^2\f^3 O(\mu^2e^{2\mu
|\psi|_{C(\cl{G})}})+\l\f^3O(e^{2\mu
|\psi|_{C(\cl{G})}}).
  \end{array}
  $$
From the definition of $\cB$ (See
\eqref{c1e15}), we have that
 $$
 \begin{array}{ll}
 \cB \3n&=\ds -4\l^3\mu^4\f^3\(\sum_{j,k=1}^nb^{jk}\psi_{x_j}\psi_{x_k}\)^2+\l^3\f^3O(\mu^3)+\l^2\f^3O(\mu^2e^{2\mu
|\psi|_{C(\cl{G})}})\\
 \ns
 &\ds\q +6\l^3\mu^4\f^3\(\sum_{j,k=1}^nb^{jk}\psi_{x_j}\psi_{x_k}\)^2+\l^3\f^3
O(\mu^3)+\l^2\f^3O(\mu^2e^{2\mu
|\psi|_{C(\cl{G})}})\\
 \ns&\ds\q+\l^2\f^3 O(\mu^2e^{2\mu
|\psi|_{C(\cl{G})}})+\l\f^3O(e^{2\mu
|\psi|_{C(\cl{G})}})
-2\l\mu^4\f\(\sum_{j,k=1}^nb^{jk}\psi_{x_j}\psi_{x_k}\)^2
+\l\f O(\mu^3)\\
 \ns
  &=\ds 2\l^3\mu^4\f ^3\(\sum_{ij} b^{jk}\psi_{x_j}\psi_{x_k}\)^2+\l^3\f^3O(\mu^3)
  +\l^2\f^3 O(\mu^2e^{2\mu
|\psi|_{C(\cl{G})}})+\l\f^3O(e^{2\mu
|\psi|_{C(\cl{G})}}),\end{array}
 $$
which leads to the second estimate in
(\ref{h6}).\endpf


\subsubsection{Global Carleman estimate for
backward stochastic parabolic equations}


As a key preliminary to prove Theorem \ref{b heat ob1}, we need to establish a global
Carleman estimate for the following backward
stochastic parabolic equation:
\begin{equation}\label{zzh9}
\left\{
\begin{array}{ll}
\ds
dz+\sum_{j,k=1}^n(a^{jk}z_{x_j})_{x_k}dt=fdt+ZdW(t)\q&\hb{ in } Q,\\
\ns \ds
z=0&\hb{ on }\Si,\\
\ns \ds z(T)=z_T&\hb{ in }G.
\end{array}
\right.
\end{equation}

We begin with the following known technical
result (See \cite[p. 4, Lemma 1.1]{FI} and
\cite{Luqi5} for its proof), which shows the
existence of a nonnegative function with an
arbitrarily given critical point location in
$G$.

\begin{lemma}\label{hl1}
For any nonempty open subset $G_1$ of $G$, there
is a  $\psi\in C^\infty({\overline{G}})$ such
that $\psi>0$ in $G$, $\psi=0$ on $\Gamma$, and
$| \nabla\psi(x)|>0$ for all $x\in \cl{G
\setminus G _1}$.
\end{lemma}

Let us choose $\th$ and $\ell$ as that in
\eqref{alphad}, and $\psi$ given by Lemma
\ref{hl1} with $G_1$ being any fixed nonempty
open subset of $G$ such that $\cl{G_1}\subset
G_0$. The desired global Carleman estimate for
\eqref{zzh9} is stated as follows:

\begin{theorem}\label{c1t4}
There is a constant $\mu_0=\mu_0(G,
G_0,(a^{jk})_{n\t n},T)>0$ such that for all
$\mu\ge \mu_0$, one can find two constants
$\cC=\cC(\mu)>0$ and $\l_0=\l_0(\mu)>0$ such
that for all $\l\ge \l_0$, $f\in
L_{\dbF}^2(0,T;L^2(G))$ and $z_T\in
L^2_{\cF_T}(\Omega;L^2(G))$, the solution
$(z,Z)\in
C_{\dbF}([0,T];L^2(\Omega;L^2(G)))\times
L^2_{\dbF}(0,T;L^2(G))$  to \eqref{zzh9}
satisfies that
\begin{equation}\label{h5.2}
\begin{array}{ll}\ds
\ds
\l^3\mu^4\mathbb{E}\int_Q\th^2\f^3z^2dxdt+\l\mu^2\mathbb{E}\int_Q\th^2\f|\n
z|^2dxdt\\
\ns \le\ds  \cC \( \l^3\mu^4\mathbb{E}
\int_{Q_0} \th^2\f^3 z^2dxdt + \mathbb{E}\int_Q
\th^2f^2dxdt+ \l^2\mu^2\mathbb{E} \int_Q
\th^2\f^2Z^2dxdt \).
\end{array}
\end{equation}
\end{theorem}

{\it Proof}\,:  We use Theorem \ref{c1t2} with
$b^{jk}$ and $h$ replaced respectively by
$-a^{jk}$ and $z$ (and hence $w=\th z$).

Integrating the equality \eqref{c1e2a} (with
$b^{jk}$ replaced by $-a^{jk}$) on $G$, taking
mean value in both sides, and noting \eqref{h6},
we conclude that
\begin{equation}\label{6h10}
\begin{array}{ll} \ds
2\mathbb{E}\int_Q\th\[\sum_{j,k=1}^n
(a^{jk}w_{x_j})_{x_k}+ \cA
w\]\[dz+\sum_{j,k=1}^n (a^{jk}z_{x_j})_{x_k}dt
\]dx-2\mathbb{E}\int_Q
\sum_{j,k=1}^n (a^{jk}w_{x_j}dw)_{x_k}dx\\
\ns \ds\q+
2\mathbb{E}\int_Q\sum_{j,k=1}^n\[\sum_{j',k'=1}^n\(2a^{jk}
a^{j'k'}\ell_{x_{j'}}w_{x_j}w_{x_{k'}}
-a^{jk}a^{j'k'}\ell_{x_j}w_{x_{j'}}w_{x_{k'}}\)-\Psi
a^{jk}w_{x_j}w\\
\ns \ds\q+ a^{jk}\(\cA\ell_{x_j}+\frac{\Psi_{x_j}}{2}\)w^2\]_{x_k}dxdt\\
\ns \ds \geq 2s_0^2\mathbb{E}\int_Q\Big[\f
\big(\l\mu^2|\n\psi|^2+\l O(\mu)\big)|\n
w|^2+\f^3\Big(\l^3\mu^4|\n
\psi|^4+\l^3O(\mu^3)\\
\noalign{\ss} \displaystyle\q
+\l^2O(\mu^2e^{2\mu |\psi|_{C(\cl{G})}})+\l
O(e^{2\mu |\psi|_{C(\cl{G})}})\Big)w^2\Big] dxdt
+2\mathbb{E} \int_Q \|\sum_{j,k=1}^n
(a^{jk}w_{x_j})_{x_k}+ \cA
w\|^2dxdt\\
\ns \ds\q +\mathbb{E}\int_Q\th^2 \sum_{j,k=1}^n
a^{jk}
(dz_{x_j}+\ell_{x_j}dz)(dz_{x_k}+\ell_{x_k}dz)dx-
\mathbb{E} \int_Q \th^2\cA(dz)^2dx,
\end{array}
\end{equation}
where
 $$
\cA=\sum_{j,k=1}^n
(a^{jk}\ell_{x_j}\ell_{x_k}-a^{jk}_{x_k}\ell_{x_j}
+a^{jk}\ell_{x_jx_k})-\ell_t,\qq
\Psi=-2\sum_{j,k=1}^na^{jk}\ell_{x_jx_k}.
$$

It follows from (\ref{zzh9}) that
\begin{equation}\label{6h11}
\begin{array}{ll} \ds
2\mathbb{E}\int_Q\th\[\sum_{j,k=1}^n (a^{jk}w_{x_j})_{x_k}+\cA w\]\[dz+\sum_{j,k=1}^n(a^{jk}z_{x_j})_{x_k}dt\]dx\\
\ns \ds=2\mathbb{E}\int_Q\th\[\sum_{j,k=1}^n
(a^{jk}w_{x_j})_{x_k}+\cA
w\]\big(fdt+ZdW(t)\big)dx
\\\ns\ds =2\mathbb{E}\int_Q\th\[-\sum_{j,k=1}^n
(a^{jk}w_{x_j})_{x_k}+\cA w\]fdtdx\\
\ns \ds\le \mathbb{E}\int_Q\|\sum_{j,k=1}^n
(a^{jk}w_{x_j})_{x_k}+\cA
w\|^2dtdx+\mathbb{E}\int_Q\th^2 f^2dtdx.
\end{array}
\end{equation}

It is clear that the term
``$ \mathbb{E}\int_Q\th^2 \sum_{j,k=1}^n a^{jk}
(dz_{x_j}+\ell_{x_j}dz)(dz_{x_k}+\ell_{x_k}dz)dx
$"
in \eqref{6h10} is nonnegative. Hence, by
\eqref{6h10}--\eqref{6h11},  one can show that
\begin{equation}\label{6h16}
\begin{array}{ll}
\displaystyle 2s_0^2\mathbb{E}\int_Q\Big[\f
\big(\l\mu^2|\n\psi|^2+\l O(\mu)\big)|\n
w|^2+\f^3\Big(\l^3\mu^4|\n
\psi|^4+\l^3O(\mu^3)\\
\noalign{\ss} \displaystyle\q
+\l^2O(\mu^2e^{2\mu |\psi|_{C(\cl{G})}})+\l
O(e^{2\mu
|\psi|_{C(\cl{G})}})\Big)w^2\Big]dxdt\\
\ns \ds\leq \mathbb{E}\int_Q\th^2(f^2+\cA Z^2)
dxdt.
\end{array}
\end{equation}

From \eqref{6h16}, we conclude that there is a
$\mu_0>0$ such that for all $\mu\ge \mu_0$, one
can find a constant $\l_0=\l_0(\mu)$ so that for
any $\l\ge \l_0$, it holds that
\begin{equation}\label{nh20}
\begin{array}{ll}\ds
\l\mu^2\mathbb{E}\int_Q\th^2\f\big(|\n
z|^2+\l^2\mu^2\f^2z^2\big)dxdt\\
\ns \ds\leq
\cC\[\mathbb{E}\int_Q\th^2(f^2+\l^2\mu^2\f^2
Z^2)
dxdt+\l\mu^2\mathbb{E}\int_0^T\int_{G_1}\th^2\f\big(|\n
z|^2+\l^2\mu^2\f^2 z^2\big)dxdt\].
\end{array}
\end{equation}
Choose a cut-off function $\zeta\in
C_0^{\infty}(G_0;[0,1])$ so that $\zeta\equiv 1$
in $G_1$. By
$ d(\th^2\f h^2)=h^2(\th^2\f)_tdt+2\th^2\f h
dh+\th^2\f(dh)^2$,
recalling $\lim_{t\to0^+}\f(t,\cd)$ $=\lim_{t\to
T^-}\f(t,\cd)\equiv 0$ and using \eqref{zzh9},
we find that
$$
\begin{array}{ll} 0\3n
&\ds=\mathbb{E}\int_{Q_0}\th^2\Big[\zeta^2z^2(\f_t+2\l\f\eta_t)+2\zeta^2\f\sum_{j,k=1}^na^{jk}z_{x_j}z_{x_k}+2\mu
\zeta^2\f(1
+2\l\f)z\sum_{j,k=1}^na^{jk}z_{x_j}\psi_{x_k}\\
\ns &\ds\qq+4\zeta\f
z\sum_{j,k=1}^na^{jk}z_{x_j}\zeta_{x_k}
+2\zeta^2\f fz+\zeta^2\f Z^2\Big]dxdt.
\end{array}
$$
Therefore, for any $\e>0$, one has
\begin{equation}\label{nh21}
\begin{array}{ll}
\ds 2\mathbb{E}\int_{Q_0}\th^2\zeta^2\f\sum_{j,k=1}^na^{jk}z_{x_j}z_{x_k}dxdt+\mathbb{E}\int_{Q_0}\th^2\zeta^2\f Z^2dxdt\\
\ns \ds\le
\e\mathbb{E}\int_{Q_0}\th^2\zeta^2\f|\n
z|^2dxdt+ \frac{\cC}{\e}\mathbb{E}\int_{Q_0}
\th^2\(\frac{1}{\l^2\mu^2}f^2+\l^2\mu^2\f^3z^2\)dxdt.
\end{array}
\end{equation}
Since the matrix $(a^{jk})_{1\leq i,j\leq n}$ is
uniformly positive definite, we conclude from
(\ref{nh21}) that
\begin{equation}\label{nh22}
\mathbb{E}\int_0^T\int_{G_1}\th^2\f|\n
z|^2dxdt\le \cC\mathbb{E}\int_{Q_0}
\th^2\(\frac{1}{\l^2\mu^2}f^2+\l^2\mu^2\f^3z^2\)dxdt.
\end{equation}

Combining (\ref{nh20}) and (\ref{nh22}), we
obtain \eqref{h5.2}. This completes the proof of
Theorem \ref{c1t4}.\endpf

\subsubsection{Proof of the observability estimate for backward stochastic
parabolic equations}

We are now in a position to prove Theorem \ref{b heat ob1}.

\ms

{\it Proof of Theorem \ref{b heat ob1}}\,:
Applying Theorem \ref{c1t4} to the equation
\eqref{dual heat 1.1}, we deduce that, for all $\mu\ge
\mu_0$ and $\l\ge \l_0 (\mu)$,
\begin{equation}\label{zth7}
\begin{array}{ll}
\ds
\l^3\mu^4\mathbb{E}\int_Q\th^2\f^3z^2dxdt+\l\mu^2\mathbb{E}\int_Q\th^2\f|\n
z|^2dxdt\\
\ns \ds\le
\cC\Big\{\l^3\mu^4\mathbb{E}\!\int_{Q_0}\!\th^2\f^3
z^2dxdt\!+\!\mathbb{E}\!\int_Q\!
\th^2\Big[\!\sum_{j=1}^n\big(a_{1j}
z\big)_{x_j}\!-\!a_2 z\!-\! a_3 Z\Big]^2dxdt
\!+\!\l^2\mu^2\mathbb{E}\int_Q\!\th^2\f^2Z^2dxdt\Big\}\\
\ns \ds \le
\cC\[\l^3\mu^4\mathbb{E}\!\int_{Q_0}\!\th^2\f^3
z^2dxdt\!+\!\mathbb{E}\!\int_Q\! \th^2\big(|\n
z|^2\!+\!\l^2\mu^2\f^2z^2+Z^2\big)dxdt
+\l^2\mu^2\mathbb{E}\int_Q\th^2\f^2Z^2dxdt\].
\end{array}
\end{equation}
Choosing $\mu=\mu_0$ and $\l=\cC$,
from (\ref{zth7}), we obtain that
\begin{equation}\label{zt7}
\mathbb{E}\int_Q\th^2\f^3 z^2dxdt\le \cC \(\mathbb{E}\int_{Q_0}\th^2\f^3
z^2dxdt+ \mathbb{E}\int_Q\th^2\f^2Z^2dxdt\).
\end{equation}
Recalling \eqref{alphad}, it follows from
\eqref{zt7} that
\begin{equation}\label{zhoppo07}
\begin{array}{ll}\ds
\mathbb{E}\int_{T/4}^{3T/4}\int_G z^2dxdt\\
\ns \ds\le \cC \frac{\ds\max_{(t,x)\in
Q}\(\th^2(t,x)\f^3(t,x) +\th^2(t,x)\f^2(t,x)
\)}{\ds\min_{x\in G}\(\th^2(T/4,x)\f^3(T/2,x)
\)}\(\mathbb{E}\int_{Q_0}
z^2dxdt + \mathbb{E}\int_QZ^2dxdt\)\\
\ns \ds\le \cC \(\mathbb{E}\int_{Q_0} z^2dxdt+
\mathbb{E}\int_QZ^2dxdt\).
\end{array}
\end{equation}

By \eqref{20160619e1} in Proposition
\ref{20160619c1}, it follows that
\begin{equation}\label{lk---j12}
\mathbb{E}\int_G z^2(0)dx\le \cC\mathbb{E}\int_G z^2(t)dx,\qq\forall\;t\in
[0,T].
 \end{equation}
Combining  (\ref{zhoppo07}) and
(\ref{lk---j12}), we conclude that, the solution
$(z,Z)$ to the equation \eqref{dual heat 1.1}
satisfies \eqref{b heat obser1}. This completes
the proof of Theorem \ref{b heat ob1}.
\endpf


\section{Pontryagin-type maximum
principle for controlled stochastic evolution equations in infinite dimensions}\label{ch-Pon}


This section is addressed to studying
the first order necessary optimality condition, i.e., Pontryagin-type maximum principle, for optimal control problems for nonlinear stochastic evolution equations in infinite dimensions, in which
both drift and diffusion terms can contain the
control variables, and the control domain  is
allowed to be nonconvex. The results in this part are taken from
\cite{LZ1,LZ2}.

\subsection{Formulation of the problem}

In this section, $U$ and $\cU[0,T] $ are the same as that in Section \ref{s3}. To simplify the presentation, we assume that $H$ is a separable, real  Hilbert space.

Let us impose the following condition.

\ms

\no{\bf (B1)} {\it For $\f=a,b$, suppose that
$\f(\cd,\cd,\cd):[0,T]\times H\times U\to H$
satisfies : i) For any $(x,u)\in H\times U$, the
functions $\f(\cd,x,u):[0,T]\to H$ is Lebesgue
measurable; ii) For any $(t,x)\in [0,T]\times
H$, the functions $\f(t,x,\cd):U\to H$ is continuous;
and iii) For all $(t,x_1,x_2,u)\in [0,T]\times H\times
H\times U$,
\begin{equation}\label{ab0}
\left\{
\begin{array}{ll}\ds
|\f(t,x_1,u) - \f(t,x_2,u)|_H  \leq
C_L|x_1-x_2|_H,\\
\ns\ds |\f(t,0,u)|_H +|b(t,0,u)|_{H} \leq C_L.
\end{array}
\right.
\end{equation}}

Consider the following controlled (forward)
stochastic evolution equation:
\begin{equation}\label{ch-10-fsystem1}
\left\{
\begin{array}{lll}\ds
dx = \big[Ax +a(t,x,u)\big]dt + b(t,x,u)dW(t) &\mbox{ in }(0,T],\\
\ns\ds x(0)=x_0,
\end{array}
\right.
\end{equation}
where $u\in \cU[0,T]$ and $x_0\in
L^{p_0}_{\cF_0}(\Omega;H)$ for a given $p_0\geq
2$.  Under the assumption (B1), one can show that the equation
\eqref{ch-10-fsystem1} is well-posed in the
sense of mild solution.

Also, we need the following condition:

\ms

\no{\bf (B2)} {\it Suppose that
$g(\cd,\cd,\cd):[0,T]\times H\times U\to \dbR$
and $h(\cd):H\to \dbR$ are two functions
satisfying: i) For any $(x,u)\in H\times U$, the
function $g(\cd,x,u):[0,T]\to \dbR$ is Lebesgue
measurable; ii) For any $(t,x)\in [0,T]\times
H$, the function $g(t,x,\cd):U\to \dbR$ is
continuous; and iii) For all $(t,x_1,x_2,u)\in [0,T]\times
H\times H\times U$,
\begin{equation}\label{ch-10-gh} \left\{
\begin{array}{ll}\ds
|g(t,x_1,u) - g(t,x_2,u)|_{H} +|h(x_1) -
h(x_2)|_H
 \leq C_L|x_1-x_2|_H,\\
\ns\ds |g(t,0,u)|_H +|h(0)|_H \leq C_L.
\end{array}
\right.
\end{equation}}

\ms

Define a cost functional $\cJ(\cdot)$ (for the
controlled system \eqref{ch-10-fsystem1}) as
follows:
\begin{equation}\label{jk1xu}
\cJ(u(\cdot))\triangleq \dbE\Big[\int_0^T
g(t,x(t),u(t))dt + h(x(T))\Big],\q\forall\,
u(\cdot)\in \cU[0,T],
\end{equation}
where $x(\cd)$ is the corresponding solution to
\eqref{ch-10-fsystem1}.

Let us consider the following optimal control
problem for the system \eqref{ch-10-fsystem1}:

\ms

\no {\bf Problem (OP)} {\it Find a $\bar
u(\cdot)\in \cU[0,T]$ such that
\begin{equation}\label{jk2}
\ds\cJ (\bar u(\cdot)) = \inf_{u(\cdot)\in
\cU[0,T]} \cJ (u(\cdot)).
\end{equation}
Any $\bar u(\cdot)$ satisfying (\ref{jk2}) is
called an {\it optimal control}. The
corresponding state process $\bar x(\cdot)$ is
called an {\it optimal state (process)}, and
$(\bar x(\cdot),\bar u(\cdot))$ is called an
{\it optimal pair}.}

\ms

The main aim of this section is to derive the first order necessary optimality condition, i.e., Pontryagin-type maximum principle, for the above Problem (OP).

As we shall see later, the main difficulty to deal with the case of non-convex control domain $U$ is that one
needs to study the following $\cL(H)$-valued
backward stochastic evolution equation\footnote{Throughout this section, for any
operator-valued process (\resp random variable)
$R$, we denote by $R^*$ its pointwise dual
operator-valued process (\resp random variable).
For example, if $R\in L^{r_1}_\dbF(0,T;
L^{r_2}(\Omega; \cL(H)))$, then $R^*\in
L^{r_1}_\dbF(0,T; L^{r_2}(\Omega; \cL(H)))$, and
$|R|_{L^{r_1}_\dbF(0,T; L^{r_2}(\Omega;
\cL(H)))}=|R^*|_{L^{r_1}_\dbF(0,T;
L^{r_2}(\Omega; \cL(H)))}$.}:
\begin{equation}\label{op-bsystem3}
\left\{
\begin{array}{ll}
\ds dP  =  - (A^*  + J^* )P dt  -  P(A + J )dt
-K^*PKdt \\
\ns\ds \hspace{2.4cm} - (K^* Q +  Q K)dt
  +   Fdt  +  Q dW(t) \qq\mbox{ in } [0,T),\\
\ns\ds P(T) = P_T.
\end{array}
\right.
\end{equation}
Here and henceforth, $F\in
L^1_\dbF(0,T;L^2(\Omega;\cL(H)))$, $P_T\in
L^2_{\cF_T}(\Omega;\cL(H))$, $J\in L^4_\dbF(0,T;
L^\infty(\Omega; $ $\cL(H)))$ and $K\in
L^4_\dbF(0,T; L^\infty(\Omega; \cL(H)))$. For
the special case when $H=\dbR^n$, it is easy to
see that \eqref{op-bsystem3} is an
$\dbR^{n\times n}$ (matrix)-backward stochastic
differential equation, and therefore, the
desired well-posedness follows from that of an
$\dbR^{n^2}$(vector)-valued backward stochastic
differential equation. However, one has to face
a real challenge in the study of
\eqref{op-bsystem3} when $\dim H=\infty$,
without further assumption on the data $F$ and
$P_T$. Indeed, in the infinite dimensional
setting, although $\cL(H)$ is still a Banach
space, it is neither reflexive (needless to say
to be a Hilbert space) nor separable even if $H$
itself is separable (See Problem 99 in
\cite{Halmos}). As far as we know, in the
previous literatures there exists no such a
stochastic integration/evolution equation theory
in general Banach spaces that can be employed to
treat the well-posedness of \eqref{op-bsystem3}.
For example, the existing result on stochastic
integration/evolution equation in UMD Banach
spaces (e.g. \cite{NVW2}) does not fit the
present case because, if a Banach space is UMD,
then it is reflexive.

To overcome the above-mentioned difficulty, we employ the stochastic transposition method developed in \cite{LZ}. More precisely, we introduce a
concept of relaxed transposition solution to the
equation \eqref{op-bsystem3}, and develop a way to study the corresponding well-posedness.
Our method can be further modified to treat the second order necessary conditions for stochastic optimal controls and the feedback control design for linear quadratic stochastic optimal control problems in infinite dimensions but all of these topics are beyond the scope of this short course.


\subsection{Pontryagin-type maximum principle for
convex control domain}


In this subsection, we give a
necessary condition for optimal controls of Problem (OP)
for the
case of special control domain $U$, i.e.,
$U$ is a convex subset of another separable
Hilbert space $\wt H$, and the metric of
$U$ is introduced by the norm of $\wt H$
(i.e.,
$\mathbf{d}(u_1,u_2)=|u_1-u_2|_{\wt H}$).

First, we need to study the following $H$-valued
backward stochastic evolution equation:
\begin{eqnarray}\label{ch-10-bsystem1}
\left\{
\begin{array}{lll}
\ds dy(t) = -  A^* y(t) dt + f(t,y(t),Y(t))dt + Y(t) dW(t) &\mbox{ in }[0,T),\\
\ns\ds y(T) = y_T.
\end{array}
\right.
\end{eqnarray}
Here $y_T \in L_{\cF_T}^{2}(\Omega;H)$),
$f(\cd,\cd,\cd):[0,T]\times H\times H \to H$
satisfies
\begin{equation}\label{Lm1.1.1}
\left\{
\begin{array}{ll}\ds
f(\cd,0,0)\in
L^1_{\dbF}(0,T;
L^2(\Omega;H)),\\\ns\ds
|f(t,x_1,y_1)-f(t,x_2,y_2)|_H\leq
C_L\big(|x_1-x_2|_H+|y_1-y_2|_H \big),\\
\ns\ds\hspace{3.5cm} \ae (t,\omega)\in
[0,T]\times\Omega,\;\;
\forall\;x_1,x_2,y_1,y_2\in H.
\end{array}
\right.
\end{equation}

To define the solution to \eqref{Lm1.1.1}, we introduce the following (forward)
stochastic evolution equation:
\begin{eqnarray}\label{ch-4-fsystem2}
\left\{
\begin{array}{lll}\ds
dz = (A^*z + v_1)ds +  v_2 dW(s) &\mbox{ in }(t,T],\\
\ns\ds z(t)=\eta,
\end{array}
\right.
\end{eqnarray}
where $t\in[0,T]$, $v_1\in
L^1_{\dbF}(t,T;L^{2}(\Omega;H))$, $v_2\in
L^{2}_{\dbF}(t,T;H)$, $\eta\in
L^{2}_{\cF_t}(\Omega;H)$. The equation
\eqref{ch-4-fsystem2} admits a unique mild
solution $z\in C_\dbF([t,T];L^2(\Omega;H))$, and
\begin{equation}\label{12.12-eq1}
\ba{ll} \ds|z|_{C_\dbF([t,T];L^2(\Omega;H))}
\leq \cC\big(|\eta|_{L^{2}_{\cF_t}(\Omega;H)} +
|v_1|_{L^1_{\dbF}(t,T;L^{2}(\Omega;H))} +
|v_2|_{L^{2}_{\dbF}(t,T;H)} \big). \ea
\end{equation}

We now introduce the following notion.

\begin{definition}\label{ch-8-definition1}
We call $ (y(\cdot), Y(\cdot)) \in
D_{\dbF}([0,T];L^{2}(\Omega;H)) \times
L^2_{\dbF}(0,T; H) $ a transposition
solution to \eqref{ch-10-bsystem1}
if for any $t\in [0,T]$, $v_1(\cdot)\in
L^1_{\dbF}(t,T;L^2(\Omega;H))$, $v_2(\cdot)\in
L^2_{\dbF}(t,T;H)$, $\eta\in
L^2_{\cF_t}(\Omega;H)$ and the corresponding
solution $z\in C_{\dbF}([t,T];L^2(\Omega;H))$ to
\eqref{ch-4-fsystem2}, it holds that
\begin{equation}\label{eq def solzz}
\begin{array}{ll}\ds
\dbE \big\langle z(T),y_T\big\rangle_{H}
-\dbE\int_t^T \big\langle z(s),F(s,y(s),Y(s) )\big\rangle_Hds\\
\ns\ds = \dbE \big\langle\eta,y(t)\big\rangle_H
+ \dbE\int_t^T \big\langle
v_1(s),y(s)\big\rangle_H ds + \dbE\int_t^T
\big\langle v_2(s),Y(s)\big\rangle_{H} ds.
\end{array}
\end{equation}
\end{definition}

We have the following well-posedness result for
 \eqref{ch-10-bsystem1} (See \cite{LZ1} for its proof).

\begin{theorem}\label{ch-8-the1}
For any $y_T \in L^2_{\cF_T}(\Omega;H)$ and
$f(\cdot,0,0)\in L^1_{\dbF}(0,T;
L^2(\Omega;H))$, the equation
\eqref{ch-10-bsystem1} admits a unique
transposition solution $(y(\cdot), Y(\cdot)) \in
D_{\dbF}([0,T];L^{2}(\Omega; H)) \times
L^{2}_{\dbF}(0,T;H)$. Furthermore,
\begin{equation}\label{s4eq1z}
\begin{array}{ll}\ds
 |(y(\cdot), Y(\cdot))|_{
D_{\dbF}([0,T];L^p(\Omega;H)) \times
L^2_{\dbF}(0,T;H)} \leq \cC\big(|y_T|_{
L^p_{\cF_T}(\Omega;H)}+
 |f(\cdot,0,0)|_{ L^1_{\dbF}(0,T;L^2(\Omega;H))}\big).
\end{array}
\end{equation}
\end{theorem}

We introduce the following further assumptions
for $a(\cd,\cd,\cd)$, $b(\cd,\cd,\cd)$,
$g(\cd,\cd,\cd)$ and $h(\cd)$.

\ms

\no{\bf (B3)} {\it The functions $a(t,x,u)$ and
$b(t,x,u)$, and the functional $g(t,x,u)$ and
$h(x)$ are $C^1$ with respect to $x$ and $u$.
Moreover, for any $(t,x,u)\in [0,T]\times H\times
U$,
\begin{equation}\label{ab}
\left\{
\begin{array}{ll}\ds
 |\!|a_x(t,x,u)|\!|_{\cL(H)}+|\!|b_x(t,x,u)|\!|_{\cL(H)} + |g_x(t,x,u)|_H+|h_x(x)|_{H}\leq C_L,\\
\ns\ds |\!|a_u(t,x,u)|\!|_{\cL(\wt H,H)}+
|\!|b_u(t,x,u)|\!|_{\cL(\wt H,H)} +
|g_u(t,x,u)|_{\wt H} \leq C_L.
\end{array}
\right.
\end{equation}}

\medskip

We have the following result.

\begin{theorem}\label{th max}
Assume that $x_0\in L^2_{\cF_0}(\Omega;H)$. Let
the assumptions (B1), (B2) and (B3) hold, and
let $(\bar x(\cd),\bar u(\cd))$ be an optimal
pair of Problem (OP). Let $(y(\cdot),Y(\cdot))$
be the transposition solution to the equation
\eqref{ch-10-bsystem1} with $y_T$ and
$f(\cd,\cd,\cd)$ given by
\begin{equation}\label{zv1}
\left\{
\begin{array}{ll}
\ds y_T =
-h_x\big(\bar x(T)\big),\\
\ns \ds f(t,y_1,y_2)=-a_x(t,\bar x(t),\bar
u(t))^*y_1 - b_x\big(t,\bar x(t),\bar
u(t)\big)^*y_2 + g_x\big(t,\bar x(t),\bar
u(t)\big).
\end{array}
\right.
\end{equation}
Then,
\begin{equation}\label{maxth ine1}
\begin{array}{ll}\ds
 \big\langle a_u(t,\bar x(t),\bar u(t))^* y(t)
\!+\! b_u(t,\bar x(t),\bar u(t))^*Y(t)\! -\!
g_u(t,\bar u(t),\bar x(t)), u \!-\!\bar u(t)
\big\rangle_{\wt H}
 \leq 0, \\
 \ns\ds \hspace{6cm} \ae (t,\omega)\in [0,T]\times \Omega,\q
\forall\; u \in U.
\end{array}
\end{equation}
\end{theorem}


To prove Theorem \ref{th max}, we need the
following result.
\begin{lemma}\label{lemma4}
If $F(\cdot)\in
L_{\dbF}^2(0,T;\wt H)$ and $\bar u(\cdot)\in
\cU[0,T]$ such that
\begin{equation}\label{lemma4 ine1}
 \dbE \int_0^T \big\langle F(t,\cd), u(t,\cd) -
\bar u(t,\cd) \big\rangle_{\wt H} dt \leq 0,
\end{equation}
holds for any $u(\cdot)\in \cU[0,T]$ satisfying
$u(\cd)-\bar u(\cd)\in
L^2_\dbF(0,T;L^2(\Omega;\wt H))$, then, for any
point $u\in U$, the following pointwise
inequality holds:
\begin{equation}\label{lemma4 ine2}
 \big\langle F(t,\omega), u- \bar
u(t,\omega)\big\rangle_{\wt H} \leq 0, \,\ \ae
(t,\omega)\in [0,T]\times\Omega.
\end{equation}
\end{lemma}

\medskip

{\it Proof}\,: We use the contradiction
argument. Suppose that the inequality
\eqref{lemma4 ine2} did not hold. Then, there
would exist a $u_0\in U$ and an $\e>0$ such that
$$
\a_\e\=\int_\Omega\int_0^T
\chi_{\L_\e}(t,\omega)dtd\dbP
>0,
$$
where $ \L_\e \triangleq \Big\{ (t,\omega)\in
[0,T]\times \Omega\;:\; \Re\big\langle
F(t,\omega), u_0 - \bar
u(t,\omega)\big\rangle_{\wt H} \geq \e \Big\}$,
and $\chi_{\L_\e}$ is the characteristic
function of $\L_\e$. For any $m\in \dbN$, define
$ \L_{\e,m}\=\L_\e\cap \big\{(t,\omega)\in
[0,T]\times\Omega\;\big|\;|\bar
u(t,\omega)|_{H_1}\le m\big\}$. It is clear that
$\ds\lim_{m\to\infty}\L_{\e,m}=\L_\e$. Hence,
there is an $m_\e\in\dbN$ such that
$$
\int_\Omega\int_0^T
\chi_{\L_{\e,m}}(t,\omega)dtd\dbP
>\frac{\a_\e}{2}>0, \qq \forall\;m\ge m_\e.
$$
Since $\big\langle F(\cdot), u_0 - \bar
u(\cdot)\big\rangle_{\wt H}$ is $\mathbf{F}$-adapted,
so is the process $\chi_{\L_{\e,m}}(\cdot)$.
Define
$$
\hat u_{\e,m}(t,\omega) = u_0
\chi_{\L_{\e,m}}(t,\omega)+ \bar
u(t,\omega)\chi_{\L_{\e,m}^c}(t,\omega),\q
(t,\omega)\in [0,T]\times \Omega.
$$
Noting that $|\bar u(\cd)|_{\wt H}\le m$ on
$\L_{\e,m}$, we see that $\hat
u_{\e,m}(\cdot)\in \cU[0,T]$ and satisfies $\hat
u_{\e,m}(\cd)-\bar u(\cd)\in L^2_\dbF(0,T;\wt
H)$. Hence, for any $m\ge m_\e$, we obtain that
$$
\begin{array}{ll}\ds
\dbE\int_0^T \big\langle F(t), \hat
u_{\e,m}(t) -\bar u(t) \big\rangle_{\wt H} dt
 = \int_\Omega\int_0^T
\chi_{\L_{\e,m}}(t,\omega) \big\langle
F(t,\omega), u_0 -
\bar u(t,\omega) \big\rangle_{\wt H}\, dtd\dbP  \\
\ns \ds \geq \e\int_\Omega\int_0^T
\chi_{\L_{\e,m}}(t,\omega)dtd\dbP  \geq
\frac{\e\a_\e}{2}
>0,
\end{array}
$$
which contradicts \eqref{lemma4 ine1}. This
completes the proof of Lemma \ref{lemma4}.\endpf

\ms

We are now in a position to prove Theorem \ref{th max}.

\ms

{\it Proof of Theorem \ref{th max}}\,: We use
the convex perturbation technique and divide the
proof into several steps.

\ms

{\bf Step 1}. For the optimal pair $(\bar
x(\cdot),\bar u(\cdot))$, we fix arbitrarily a
control $u(\cdot)\in \cU[0,T]$ satisfying $
u(\cd)-\bar u(\cd)\in
L^2_\dbF(0,T;L^2(\Omega;\wt H))$. Since $U$ is
convex, we see that
$$
u^\e(\cdot) = \bar u(\cdot) + \e [u(\cdot) -
\bar u(\cdot)] = (1-\e)\bar u(\cdot) + \e
u(\cdot) \in \cU[0,T], \q\forall\;\e \in [0,1].
$$
Denote by $x^\e(\cdot)$ the state process of
\eqref{ch-10-fsystem1} corresponding to the
control $u^\e(\cdot)$. It is easy to show that
\begin{equation}\label{th max eq0.0xz}
|x^\e|_{C_\dbF(0,T;L^2(\Omega;H))}\leq \cC
\big(1+|x_0|_{L^2_{\cF_0}(\Omega;H)}\big),\q
\forall\; \e \in [0,1].
\end{equation}
Write $\ds x_1^\e(\cd) =
\frac{1}{\e}\big[x^\e(\cd)-\bar x(\cd)\big]$ and
$\d u(\cd) = u(\cd) - \bar u(\cd)$. Since $(\bar
x(\cd),\bar u(\cd))$ satisfies
\eqref{ch-10-fsystem1}, it is easy to see that
$x_1^\e(\cdot)$ satisfies the following
stochastic differential equation:
\begin{equation}\label{fsystem3x}
\left\{
\begin{array}{lll}\ds
dx_1^\e = \big(Ax_1^\e +  a_1^\e x^\e_1 +
a_2^\e\d u  \big)dt + \big( b_1^\e x^\e_1 +
b_2^\e\d u \big)dW(t) &\mbox{ in
}(0,T],\\
\ns\ds x_1^\e(0)=0,
\end{array}
\right.
\end{equation}
where for $\f=a,b$,
\begin{equation}\label{tatb}
\ds   \f_1^\e (t)  = \int_0^1 \f_x(t,\bar x(t) +
\si\e x_1^\e(t), u^\e(t))d\si, \q \f_2^\e (t) =
\int_0^1 \f_u(t,\bar x(t), \bar u(t)+\si\e\d
u(t))d\si.
\end{equation}
Consider the following stochastic differential
equation:
\begin{equation}\label{fsystem3.1}
\left\{
\begin{array}{lll}\ds
dx_2 = \big[Ax_2 + a_1(t)x_2 +  a_2(t)\d u
\big]dt + \big[ b_1(t) x_2 + b_2(t)\d u \big]
dW(t) &\mbox{ in
}(0,T],\\
\ns\ds x_2(0)=0,
\end{array}
\right.
\end{equation}
where for $\f=a,b$,
\begin{equation}\label{barab}
\f_1(t) = a_x(t,\bar x(t),\bar u(t)),\q \f_2(t)
= a_u(t,\bar x(t),\bar u(t)).
\end{equation}

\ms

{\bf Step 2}. In this step, we shall show that
\begin{equation}\label{th max eq0.2}
\lim_{\e\to 0+} \big|x_1^\e -
x_2\big|_{L_\dbF^\infty(0,T;L^2(\Omega;H))}=0.
\end{equation}

First, using Burkholder-Davis-Gundy's inequality (See Theorem \ref{BDG})
and by the assumption (B1), we find that
\begin{equation}\label{th max eq0}
\begin{array}{ll}\ds
\mE|x_1^\e(t)|^2_H \3n&\ds  = \mE\Big| \int_0^t
S(t-s) a_1^\e(s) x_1^\e(s) ds + \int_0^t S(t-s)
a_2^\e(s)\d u(s) ds  \\
\ns&\ds \q +
\int_0^t S(t-s) b_1^\e(s) x_1^\e(s) dw(s) +  \int_0^t S(t-s)b_2^\e(s)\d u(s) dw(s)\Big|_H^2\\
\ns&\ds  \leq \cC \mE\[\Big| \int_0^t S(t-s)
a_1^\e(s) x_1^\e(s) ds \Big|_H^2
+  \Big|\int_0^t S(t-s) b_1^\e(s) x_1^\e(s) dw(s) \Big|_H^2 \\
\ns&\ds \q  + \Big| \int_0^t S(t-s)a_2^\e(s)\d u(s)ds\Big|_H^2 + \Big| \int_0^t S(t-s)b_2^\e(s)\d u(s)dw(s)\Big|_H^2\] \\
\ns&\ds   \leq \cC\[ \int_0^t
\mE|x_1^\e(s)|_H^2ds + \int_0^T\mE |\d
u(s)|_{H_1}^2 dt\].
\end{array}
\end{equation}
It follows from \eqref{th max eq0} and
Gronwall's inequality that
\begin{equation}\label{th max eq0.1}
\begin{array}{ll}\ds
\mE|x_1^\e(t)|^2_H \leq  \cC|\bar u -
u|^2_{L^2_\dbF(0,T;H_1)},\q\forall\;t\in [0,T].
\end{array}
\end{equation}
By a similar computation, we see that
\begin{equation}\label{th max eq2}
\begin{array}{ll}\ds
\mE|x_2(t)|^2_H \leq  \cC |\bar u -
u|^2_{L^2_\dbF(0,T;H_1)},\q\forall\;t\in [0,T].
\end{array}
\end{equation}

On the other hand, put $x_3^\e = x_1^\e - x_2$.
Then, $x_3^\e$ solves the following equation:
\begin{equation}\label{fsystem3x2}
\left\{
\begin{array}{lll}\ds
dx_3^\e = \big[Ax_3^\e +  a_1^\e(t) x^\e_3 +
\big(a_1^\e(t)- a_1(t)\big)x_2 +
\big( a_2^\e(t)-a_2(t)\big)\d u \big]dt \\
\ns\ds \hspace{1.3cm} + \big[  b_1^\e(t) x^\e_3
+ \big(b^\e_1(t) - b_1(t)\big) x_2 + \big(
b_2^\e(t)-b_2(t)\big)\d u \big] dW(t) &\mbox{ in
}(0,T],\\
\ns\ds x_3^\e(0)=0.
\end{array}
\right.
\end{equation}
It follows from \eqref{th max
eq2}--\eqref{fsystem3x2} that
\begin{equation}
\begin{array}{ll}\ds
\mE|x_3^\e(t)|^2_{H} \\
\ns\ds = \mE\Big| \int_0^t S(t-s)
a_1^\e(s)x_3^\e(s)ds
+ \int_0^t S(t-s) b_1^\e(s)x_3^\e(s) dW(s) \\
\ns\ds\q  + \int_0^t S(t - s)\big[ a^\e_1\!(s)
-\!a_1\!(s)\big] x_2\!(s) ds \!+ \int_0^t
S(t\!-\!s)\big[ b_1^\e\!(s) \!-\! b_1(s)\big]
x_2 (s) dW(s) \nonumber
\\
\ns\ds\q + \int_0^t S(t - s)\big[ a^\e_2 (s)-a_2
(s)\big]\d u(s) ds
+ \int_0^t S(t - s)\big[ b_2^\e (s)  -  b_2 (s)\big]\d u(s)dW(s)\Big|_H^2\\
\ns\ds  \leq \cC\[\mE\!\int_0^t\!
|x_3^\e(s)|_H^2 ds\! +\!
|x_2(\cd)|^2_{L^\infty_\dbF(0,T;L^2(\Omega;H))}
\int_0^T\!
\mE\big(|\!|a^\e_1(s)\!-a_1(s)|\!|_{\cL(H)}^2
\!+ |\!|b^\e_1(s) \!- b_1(s)|\!|_{\cL(H)}^2
\big) dt
\\
\ns\ds \q +  |u - \bar
u|^2_{L^2_\dbF(0,T;L^2(\Omega;H_1))} \int_0^T
\mE\big(|\!|a_2^\e (s) - a_2
(s)|\!|_{\cL(H_1,H)}^2 + |\!|b_2^\e (s) - b_2
(s)|\!|_{\cL(H_1,H)}^2 \big)dt
\]\\
\ns\ds \leq  \cC(1 \!+\! |u\! -\! \bar
u|^2_{L^2_\dbF(0,T;L^2(\Omega;H_1))})\Big\{\mE\!
\int_0^t\! |x_3^\e(s)|_H^2 ds \!+ \!\int_0^T\!\!
\mE\[ |\!|a^\e_1(s)\! - \!a_1(s)|\!|_{\cL(H)}^2
\!+\! |\!|b^\e_1(s) \!-\!
b_1(s)|\!|_{\cL(H)}^2 \\
\ns\ds \q + |\!|a_2^\e(s) -
a_2(s)|\!|_{\cL(H_1,H)}^2 + |\!|b_2^\e(s) -
b_2(s)|\!|_{\cL(H_1,H)}^2
\]dt\Big\}.
\end{array}
\end{equation}
This, together with Gronwall's inequality,
implies that
\begin{equation}\label{th max eq3}
\begin{array}{ll}\ds
\mE|x_3^\e(t)|^2_{H} \3n&\ds \leq \cC
e^{\cC|u-\bar
u|_{L^2_\dbF(0,T;L^2(\Omega;H_1))}}\int_0^T \mE\[ |\!|a^\e_1(s)-a_1(s)|\!|_{\cL(H)}^2  + |\!|b^\e_1(s)-b_1(s)|\!|_{\cL(H)}^2\\
\ns&\ds \q + |\!|a_2^\e(s) -
a_2(s)|\!|_{\cL(H_1,H)}^2 + |\!|b_2^\e(s) -
b_2(s)|\!|_{\cL(H_1,H)}^2 \]ds,\q \forall\, t\in
[0,T].
\end{array}
\end{equation}

Note that \eqref{th max eq0.1} implies
$x^\e(\cd)\to \bar x(\cd)$ (in $H$) in
probability, as $\e\to0$. Hence, by
\eqref{tatb}, \eqref{barab} and the continuity
of $a_x(t,\cd,\cd)$, $b_x(t,\cd,\cd)$,
$a_u(t,\cd,\cd)$ and $b_u(t,\cd,\cd)$, we deduce
that
$$
\ba{ll}\ds \lim_{\e\to 0}\int_0^T\mE\[
|\!|a^\e_1(s)-a_1(s)|\!|_{\cL(H)}^2 +
|\!|b^\e_1(s)-b_1(s)|\!|_{\cL(H)}^2
\\\ns\ds\qq\qq+ |\!|a_2^\e(s) - a_2(s)|\!|_{\cL(H_1,H)}^2 +|\!|b_2^\e(s) -
b_2(s)|\!|_{\cL(H_1,H)}^2 \]ds=0. \ea
$$
This, combined with (\ref{th max eq3}), gives
\eqref{th max eq0.2}.

\ms

{\bf Step 3}. Since $(\bar x(\cdot),\bar
u(\cdot))$ is an optimal pair of Problem (OP),
from \eqref{th max eq0.2}, we find that
\begin{equation}\label{var 1}
\begin{array}{ll}\ds
 0\leq \lim_{\e\to 0}\frac{\cJ(u^\e(\cdot)) - \cJ(\bar u(\cdot))}{\e} \\
\ns\ds\;\;\;= \Big\{\dbE\!\int_0^T\!\!
\Big[\big\langle g_1(t),x_2(t)\big\rangle_H
\!+\! \big\langle g_2(t),\d u(t)
\big\rangle_{H_1}\Big] dt + \dbE\big\langle
h_x(\bar x(T)),x_2(T)\big\rangle_H\Big\},
\end{array}
\end{equation}
where
$$
g_1(t) = g_x(t,\bar x(t),\bar u(t)),\q g_2(t) =
g_u(t,\bar x(t),\bar u(t)).
$$

Now, it follows from It\^os formula that
\begin{equation}\label{max eq1}
\begin{array}{ll}\ds
-\dbE \big\langle h_x(\bar
x(T)),x_2(T)\big\rangle_H - \dbE\int_0^T
\big\langle g_1(t),x_2(t)\big\rangle_H dt \\
\ns\ds= \dbE \int_0^T \Big[\big\langle a_2(t)\d
u(t), y(t)\big\rangle_H + \big\langle b_2(t)\d
u(t), Y(t)\big\rangle_H\Big]dt.
\end{array}
\end{equation}
Combining \eqref{var 1} and \eqref{max eq1}, we
find
\begin{eqnarray}\label{max ine2}
 \dbE\int_0^T\big\langle a_2(t)^* y(t) +
b_2(t)^*Y(t) - g_2(t), u(t)-\bar u(t)
\big\rangle_{\wt H}dt\leq 0
\end{eqnarray}
holds for any $u(\cdot)\in \cU[0,T]$ satisfying
$u(\cd)-\bar u(\cd)\in
L^2_\dbF(0,T;L^2(\Omega;\wt H))$. Hence, by
means of Lemma \ref{lemma4}, we conclude that
\begin{equation}\label{max ine3}
 \big\langle a_2(t)^* y(t) + b_2(t)^*Y(t) -
g_2(t), u - \bar u(t) \big\rangle_{\wt H} \leq
0, \qq\ae [0,T]\times \Omega,\;\forall \;u \in U.
\end{equation}
This completes the proof of Theorem \ref{th
max}.\endpf


\subsection{Pontryagin-type maximum principle for the
general
case}


In this subsection, we give a
necessary condition for optimal controls of Problem (OP)
for the general
case.

\subsubsection{Relaxed transposition solution to
operator-valued backward stochastic evolution equations}


To define the solution to \eqref{op-bsystem3} in
the transposition sense, we need to introduce
the following two (forward) stochastic evolution
equations:
\begin{equation}\label{op-fsystem1}
\left\{
\begin{array}{ll}
\ds dx_1 = (A+J)x_1ds + u_1ds + Kx_1 dW(s) + v_1 dW(s) &\mbox{ in } (t,T],\\
\ns\ds x_1(t)=\xi_1
\end{array}
\right.
\end{equation}
and
\begin{equation}\label{op-fsystem2}
\left\{
\begin{array}{ll}
\ds dx_2 = (A+J)x_2ds + u_2ds + Kx_2 dW(s) + v_2 dW(s) &\mbox{ in } (t,T],\\
\ns\ds x_2(t)=\xi_2.
\end{array}
\right.
\end{equation}
Here
$$
\xi_1,\xi_2 \!\in\! L^{4}_{\cF_t}(\Omega;H),\,
u_1,u_2\in L^2_\dbF(t,T;L^{4}(\Omega;H)),\,
v_1,v_2\in L^2_\dbF(t,T;L^{4}(\Omega;H)).
$$

Also, we need to introduce the solution space
for \eqref{op-bsystem3}.  Write
\begin{equation}\label{jshi1}
\begin{array}{ll}\ds
\q C_{\dbF,w}([0,T];L^{2}(\Omega;\cL(H))\\
\ns\ds\= \Big\{P(\cd,\cd)\;\Big|\; P(\cd,\cd)\in
L^\infty_\dbF(0,T;L^2(\Omega;\cL(H))) \mbox{ and
for every }
t\in[0,T]\hb{ and }\xi\in L^{4}_{\cF_t}(\Omega;H),\\
\ns\ds\q P(\cd,\cd)\xi\in
C_{\dbF}([t,T];L^{\frac{4}{3}}(\Omega;H)) \mbox{
and }
|P(\cd,\cd)\xi|_{C_{\dbF}([t,T];L^{\frac{4}{3}}(\Omega;H))}
\leq \cC|\xi|_{L^{4}_{\cF_t}(\Omega;H)} \Big\}
\end{array}
\end{equation}
and
\begin{equation}\label{jshi2}
\begin{array}{ll}\ds \cQ[0,T] \=
\Big\{\big(Q^{(\cd)},\widehat
Q^{(\cd)}\big)\;\Big|\;\mbox{For any } t\in [0,T],
\mbox{ both }Q^{(t)}\mbox{ and }\widehat
Q^{(t)}\mbox{ are bounded
linear operators}\\
\ns\ds\hspace{1.7cm}\mbox{ from
}L^{4}_{\cF_t}(\Omega;H)\times
L^2_\dbF(t,T;L^{4}(\Omega;H))\times
L^2_\dbF(t,T;L^{4}(\Omega;H)) \mbox{ to }
L^{2}_\dbF(t,T;L^{\frac{4}{3}}(\Omega;H))\\
\ns\ds \hspace{1.7cm} \mbox{ and
}Q^{(t)}(0,0,\cd)^*=\widehat
Q^{(t)}(0,0,\cd)\Big\}.
\end{array}
\end{equation}

We now employ the stochastic transposition method, and define the relaxed transposition solution
to \eqref{op-bsystem3} as follows:

\begin{definition}\label{op-definition2x}
We call
$
\big(P(\cd),Q^{(\cd)},\widehat Q^{(\cd)}\big)\in
C_{\dbF,w}([0,T];L^{2}(\Omega; \cL(H)))\times
\cQ[0,T]
$
a relaxed transposition solution to the equation
\eqref{op-bsystem3} if for any $t\in [\tau,T]$,
$\xi_1,\xi_2\in L^{4}_{\cF_t}(\Omega;H)$,
$u_1(\cd)$, $u_2(\cd)\in
L^2_{\dbF}(t,T;L^{4}(\Omega;H))$ and $v_1(\cd),
v_2(\cd)\in L^2_{\dbF}(t,T; L^{4}(\Omega;H))$,
it holds that
\begin{equation}\label{6.18eq1}
\begin{array}{ll}
\ds \mE\big\langle P_T x_1(T), x_2(T)
\big\rangle_{H} - \mE \int_t^T \big\langle
F(s) x_1(s), x_2(s) \big\rangle_{H}ds\\
\ns\ds =\mE\big\langle P(t) \xi_1,\xi_2
\big\rangle_{H} + \mE \int_t^T \big\langle
P(s)u_1(s), x_2(s)\big\rangle_{H}ds + \mE
\int_t^T \big\langle P(s)x_1(s),
u_2(s)\big\rangle_{H}ds \\
\ns\ds \q  + \mE \int_t^T \big\langle P(s)K
(s)x_1 (s), v_2 (s)\big\rangle_{H}ds + \mE
\int_t^T \big\langle  P(s)v_1 (s), K (s)x_2 (s) + v_2(s)\big\rangle_{H}ds\\
\ns\ds \q  + \mE \int_t^T  \big\langle v_1(s),
\widehat
Q^{(t)}(\xi_2,u_2,v_2)(s)\big\rangle_{H}ds + \mE
\int_t^T \big\langle Q^{(t)}(\xi_1,u_1,v_1)(s),
v_2(s) \big\rangle_{H}ds,
\end{array}
\end{equation}
Here, $x_1(\cd)$ and $x_2(\cd)$ solve
\eqref{op-fsystem1} and \eqref{op-fsystem2},
respectively.
\end{definition}

We have the following well-posedness result for
the equation (\ref{op-bsystem3}) (See \cite{LZ1} for its proof).

\begin{theorem}\label{OP-th2}
Suppose that  $L^p_{\cF_T}(\Omega)$ ($1\leq
p < \infty$) is a separable Banach space. Then the
equation \eqref{op-bsystem3} admits one and only one
relaxed transposition solution
$\big(P(\cd),Q^{(\cd)},\widehat Q^{(\cd)}\big)\in
C_{\dbF,w}([0,T];$ $L^{2}(\Omega; \cL(H)))\times
\cQ[0,T]$.
Furthermore,
 $$
\begin{array}{ll}\ds
\q |P|_{C_{\dbF,w}([0,T];L^{2}(\Omega; \cL(H)))}
+ \big|\big(Q^{(\cd)},\widehat
Q^{(\cd)}\big)\big|_{\cQ[0,T]}\leq \cC\big(
|F|_{L^1_\dbF(0,T;\;L^2(\Omega;\cL(H)))} +
|P_T|_{L^2_{\cF_T}(\Omega;\;\cL(H))}\big).
\end{array}
 $$
\end{theorem}

Next, we give a regularity result for the
relaxed transposition solution. For this
purpose, we first give two preliminary results (See \cite{LZ2} for their proofs).
\begin{lemma}\label{lemma5}
For each $t\in[0,T]$, if $u_2=v_2=0$ in the
equation \eqref{op-fsystem2}, then there exists
an operator $U(\cd,t)\in
\cL\big(L^{4}_{\cF_t}(\Omega;H);
C_\dbF([t,T];L^{4}(\Omega;H))\big)$ such that
the solution to \eqref{op-fsystem2} can be
represented as $x_2(\cd) = U(\cd,t)\xi_2$.
\end{lemma}

Let $\{\D_n\}_{n=1}^\infty$ be a sequence of
partitions of $[0,T]$, that is,

$$
\D_n\= \Big\{t_i^n\;\Big|\;i=0,1,\cdots,n, \hb{
and }0=t_0^n < t_1^n < \cds < t_{n}^n =T\Big\}
$$
such that $\D_n\subset \D_{n+1}$ and $\d(\D_n)\=
\max_{0\leq i\leq n-1} (t_{i+1}^n - t_{i}^n)\to
0$ as $n\to\infty$. We introduce the following
subspaces of $L^2_\dbF(0,T;L^{4}(\Om;H))$:
\begin{equation}\label{cH}
\cH_n=\Big\{\sum_{i=0}^{n-1}
\chi_{[t_i^n,t_{i+1}^n)}(\cd)U(\cd,t_i^n)h_i\;\Big|\;
h_i\in L^{4}_{\cF_{t_i^n}}(\Om;H)\Big\}.
\end{equation}
Here $U(\cd,\cd)$ is the operator introduced in
Lemma \ref{lemma5}. We have the following
result.
\begin{lemma}\label{5.13-prop1}
The set $\bigcup_{n=1}^\infty \cH_n$ is dense in
$L^2_\dbF(0,T;L^{4}(\Om;H))$.
\end{lemma}

The regularity result for solutions to
\eqref{op-bsystem3} can be stated as follows (See \cite{LZ2} for its proof).

\begin{lemma}\label{10.1th}
Suppose that the assumptions in Theorem
\ref{OP-th2} hold and let
$(P(\cd),Q^{(\cd)},\widehat Q^{(\cd)})$ be the
relaxed transposition solution to the equation
\eqref{op-bsystem3}. Then,  there exist an
$n\in\dbN$ and two pointwise defined linear
operators $Q^n$ and $\widehat Q^n$, both of
which are from $\cH_n$ to
$L^{2}_\dbF(0,T;L^{\frac{4}{3}}(\Om;H))$, such
that, for any $\xi_1,\xi_2\in
L^{4}_{\cF_0}(\Om;H)$, $u_1(\cd), u_2(\cd)\in
L^{4}_\dbF(\Om;L^2(0,T;H))$ and
$v_1(\cd),v_2(\cd)\in \cH_{n}$, it holds that
\begin{equation}\label{10.9eq2}
\begin{array}{ll}
\ds \mE \int_{0}^T \big\langle v_1(s), \widehat
Q^{(0)}(\xi_2,u_2,v_2)(s) \big\rangle_{H}ds +
\mE \int_{0}^T \big\langle
 Q^{(0)}(\xi_1,u_1,v_1) (s),
v_2(s) \big\rangle_{H}ds \\
\ns\ds =\mE \int_{0}^T
 \[\big\langle
\big(Q^n v_1\big)(s), x_2 (s)
\big\rangle_{H}+\big\langle x_1 (s),
\big(\widehat Q^n v_2\big)(s)
\big\rangle_{H}\]ds,
\end{array}
\end{equation}
where, $x_1(\cd)$ and $x_2(\cd)$ solve
accordingly \eqref{op-fsystem1} and
\eqref{op-fsystem2} with $t=0$. Further, there
is a positive constant $\cC(n)$,  depending on $n$,
such that
\begin{equation}\label{4.14-eq3}
\big|Q^n
v_1\big|_{L^{2}_\dbF(0,T;L^{\frac{4}{3}}(\Om;H))}
+ \big|\widehat Q^n
v_2\big|_{L^{2}_\dbF(0,T;L^{\frac{4}{3}}(\Om;H))}\leq
\cC(n)\big(|\tilde
v_1|_{L^2_\dbF(0,T;L^{4}(\Om;H))}+|\tilde
v_2|_{L^2_\dbF(0,T;L^{4}(\Om;H))}\big),
\end{equation}
where
$$
\tilde v_1 = \sum_{i=0}^{n-1}
\chi_{[t^n_i,t^n_{i+1})}(\cd) h_i \q\mbox{ for
}\; v_1 = \sum_{i=0}^{n-1}
\chi_{[t^n_i,t^n_{i+1})}(\cd)U(\cd,t_i) h_i
$$
and
$$
\tilde v_2 = \sum_{j=0}^{n-1}
\chi_{[t^n_j,t^n_{j+1})}(\cd) h_j  \q\mbox{ for
}\; v_2 = \sum_{j=0}^{n-1}
\chi_{[t^n_j,t^n_{j+1})}(\cd)U(\cd,t_j)h_j.
$$
\end{lemma}


\subsubsection{Statement of the Pontryagin-type
maximum principle}


We assume the following further conditions for
the optimal control problem (OP).

\ms

\medskip

\no{\bf (B4)} {\it The function  $a(t,x,u)$ and
$b(t,x,u)$, and the functional $g(t,x,u)$ and
$h(x)$ are $C^2$ with respect to $x$, such that
for  $\psi(t,x,u)= g(t,x,u),h(x)$, it holds that
 $\f_x(t,x,u)$, $\psi_x(t,x,u)$,
$\f_{xx}(t,x,u)$ and $\psi_{xx}(t,x,u)$  are
continuous with respect to $u$. Moreover, for all
$(t,x,u)\in [0,T]\times H\times U$,
\begin{equation}\label{ab1}
\left\{
\begin{array}{ll}\ds
|\!|a_x(t,x,u)|\!|_{\cL(H)}+
|\!|b_x(t,x,u)|\!|_{\cL(H)}+|\psi_x(t,x,u) |_H
\leq C_L,\\
\ns\ds |\!|a_{xx}(t,x,u)|\!|_{\cL(H\times
H,\;H)} + |\!|b_{xx}(t,x,u)|\!|_{\cL(H\times
H,\;H)} +|\!|\psi_{xx}(t,x,u)|\!|_{\cL(H)}
 \leq C_L.
\end{array}
\right.
\end{equation}}

\ms

Let
\begin{equation}\label{H}
\begin{array}{ll}\ds
\dbH(t,x,u,k_1,k_2) \= \big\langle k_1,a(t,x,u)  \big\rangle_H + \big\langle k_2, b(t,x,u)  \big\rangle_{H} - g(t,x,u),\\
\ns\ds \hspace{4cm} (t,x,u,k_1,k_2)\in
[0,T]\times H \times U\times H\times H.
\end{array}
\end{equation}

We have the following result.
\begin{theorem}\label{maximum p2}
Suppose that $L^p_{\cF_T}(\Omega)$ ($1\leq
p < \infty$) is a separable Banach space, $U$ is
a separable metric space, and $x_0\in
L^8_{\cF_0}(\Omega;H)$. Let the assumptions
(B1), (B2) and (B4) hold, and let $(\bar
x(\cd),\bar u(\cd))$ be an optimal pair of
Problem (OP). Let $\big(y(\cdot),Y(\cdot)\big)$
be the transposition solution to
\eqref{ch-10-bsystem1} with $y_T$ and
$f(\cd,\cd,\cd)$ given by \eqref{zv1}.
Assume that  $(P(\cd),Q^{(\cd)},\widehat
Q^{(\cd)})$ is the relaxed transposition
solution to the equation \eqref{op-bsystem3} in
which $P_T$, $J(\cd)$, $K(\cd)$ and $F(\cd)$ are
given by\vspace{-0.1cm}
\begin{equation}\label{MP2-eq9}
\left\{
\begin{array}{ll} \ds P_T =
-h_{xx}\big(\bar x(T)\big),\q J(t) = a_x(t,\bar
x(t),\bar
u(t)), \\
\ns \ds K(t) =b_x(t,\bar x(t),\bar u(t)), \q
F(t)= -\dbH_{xx}\big(t,\bar x(t),\bar
u(t),y(t),Y(t)\big).
\end{array}
\right.
\end{equation}
Then,   for $\ae$ $(t,\omega)\in [0,T]\times
\Omega$ and for all $u \in U$,\vspace{-0.1cm}
\begin{equation}\label{MP2-eq1}
\begin{array}{ll}\ds
 \dbH\big(t,\bar x(t),\bar u(t),y(t),Y(t)\big) - \Re\dbH\big(t,\bar x(t),u,y(t),Y(t)\big) \\
\ns\ds \q - \frac{1}{2}\big\langle  P(t)\big[
b\big(t,\bar x(t),\bar u(t)\big)-b\big(t,\bar
x(t),u\big) \big], b\big(t,\bar x(t),\bar
u(t)\big)-b\big(t,\bar x(t),u\big)
\big\rangle_{H} \\
\ns\ds\geq 0.
\end{array}
\end{equation}
\end{theorem}
%


\subsubsection{Proof of the Pontryagin-type
stochastic maximum principle}


We are now in a position to prove Theorem \ref{maximum p2}.

\ms

{\it Proof of Theorem \ref{maximum p2}}\,: We
divide the proof into two steps.

\ms

{\bf Step 1}. For each $\e>0$, let $E_\e\subset
[0,T]$ be a measurable set with measure $\e$.
Put\vspace{-0.1cm}
\begin{equation}\label{s7ue}
u^\e(t) = \left\{
\begin{array}{ll}
\ds \bar u(t), & t\in [0,T]\setminus E_\e,\\
\ns\ds u(t), & t\in E_\e.
\end{array}
\right.
\end{equation}
where  $u(\cdot)$ is an arbitrary given element
in $\cU[0,T]$.

We introduce some notations which will be used
in what follows. For $\f=a,b,g$, we let
\begin{equation}\label{s7tatb1}
\left\{
\begin{array}{ll}
\ds \f_1(t) = \f_x(t,\bar x(t),\bar u(t)), \q
\f_{11}(t) = \f_{xx}(t,\bar x(t),\bar u(t)),\\
\ns\ds \tilde \f_1^\e (t)  = \int_0^1 \f_x\big(t,\bar x(t) + \si (x^\e (t) - \bar x(t)), u^\e (t)\big)d\si, \\
\ns\ds \tilde \f_{11}^\e (t)  = 2\int_0^1
(1-\si) a_{xx}\big(t,\bar x(t) + \si (x^\e(t) -
\bar x(t)), u^\e (t)\big)d\si,
\end{array}
\right.
\end{equation}
and
\begin{equation}\label{s7tatb2}
\left\{
\begin{array}{ll}
\ds \d \f(t)  = \f(t,\bar x(t), u(t)) -
\f(t,\bar x(t),\bar u(t)),
\\
\ns\ds \d \f_1(t) = \f_x(t,\bar x(t), u(t)) -
\f_x(t,\bar x(t),\bar u(t)),  \\
\ns\ds \d \f_{11}(t) = \f_{xx}(t,\bar x(t),
u(t)) - \f_{xx}(t,\bar x(t),\bar u(t)).
\end{array}
\right.
\end{equation}

Let $x^\e(\cdot)$ be the state process of the
system \eqref{ch-10-fsystem1} corresponding to
the control $u^\e(\cdot)$. Then, $x^\e(\cdot)$
solves\vspace{-0.321cm}
\begin{eqnarray}\label{s7fsystem2}
\left\{
\begin{array}{lll}\ds
d x^\e = \big[Ax^\e +a(t,x^\e,u^\e)\big]dt + b(t,x^\e,u^\e)dW(t) &\mbox{ in }(0,T],\\
\ns\ds x^\e(0)=x_0.
\end{array}
\right.
\end{eqnarray}
It is easy to prove that\vspace{-0.21cm}
\begin{equation}\label{th max eq0.0zx}
|x^\e|_{C_\dbF([0,T];L^8(\Omega;H))}\leq C
\big(1+|x_0|_{L^8_{\cF_0}(\Omega;H)}\big),\q
\forall\; \e
>0.
\end{equation}
 Let
$x_1^\e(\cd) = x^\e(\cd)-\bar x(\cd)$. Then, by
(\ref{th max eq0.0zx}) and noting that the
optimal pair $(\bar x(\cd),\bar u(\cd))$ solves
the equation \eqref{ch-10-fsystem1}, we see that
$x_1^\e(\cdot)$ satisfies the following
stochastic evolution equation:
\begin{equation}\label{fsystem3zx}
\left\{
\begin{array}{lll}\ds
dx_1^\e = \big[Ax_1^\e + \tilde a_1^\e(t) x^\e_1
+ \chi_{E_\e} (t)\d a(t) \big]dt + \big[ \tilde
b_1^\e(t) x^\e_1 + \chi_{E_\e} (t)\d b(t)
\big]dW(t) &\mbox{ in
}(0,T],\\
\ns\ds x_1^\e(0)=0.
\end{array}
\right.
\end{equation}

Consider the following two stochastic
differential equations:\vspace{-0.1cm}
\begin{equation}\label{s7fsystem3.1}
\left\{
\begin{array}{lll}\ds
dx_2^\e = \big[Ax_2^\e + a_1(t)x_2^\e \big]dt +
\big[ b_1(t) x_2^\e + \chi_{E_\e} (t)\d b(t)
\big] dW(t) &\mbox{ in
}(0,T],\\
\ns\ds x_2^\e(0)=0
\end{array}
\right.
\end{equation}
and\vspace{-0.1cm}\footnote{Recall that, for any
$C^2$-function $f(\cd)$ defined on a Banach
space $X$ and $x_0\in X$, $f_{xx}(x_0)\in \cL
(X,X;X)$. This means that, for any
$x_1,x_2\in X$, $f_{xx}(x_0)(x_1,x_2)\in X$.
Hence, by (\ref{s7tatb1}),
$a_{11}(t)\big(x_2^\e,x_2^\e\big)$ (in
(\ref{s7fsystem3.2})) stands for $a_{xx}(t,\bar
x(t),\bar u(t))\big(x_2^\e(t),x_2^\e(t)\big)$.
One has a similar meaning for
$b_{11}(t)\big(x_2^\e,x_2^\e\big)$ and so on.}
\begin{equation}\label{s7fsystem3.2}
\left\{
\begin{array}{lll}\ds
dx_3^\e = \[Ax_3^\e + a_1(t)x_3^\e +
\chi_{E_\e}(t)\d a(t) +
\frac{1}{2}a_{11}(t)\big(x_2^\e,x_2^\e\big) \]dt \\
\ns\ds \hspace{1cm} + \[ b_1(t) x_3^\e +
\chi_{E_\e} (t)\d b_1(t)x_2^\e +
\frac{1}{2}b_{11}(t)\big(x_2^\e,x_2^\e\big)\]
dW(t) &\mbox{ in
}(0,T],\\
\ns\ds x_3^\e(0)=0.
\end{array}
\right.
\end{equation}
Similar to Steps 1-2 in the proof of Theorem
\ref{th max}, we can show that\vspace{-0.1cm}
\begin{equation}\label{s7th max eq0.1}
|x_1^\e(\cd)|^8_{C_\dbF([0,T];L^8(\Om;H))} \leq
\cC(x_0) \e^4,
\end{equation}
\begin{equation}\label{s7th max eq2}
|x_2^\e(\cd)|^8_{C_\dbF([0,T];L^8(\Om;H))} \leq
\cC(x_0)\e^4,
\end{equation}
\begin{equation}\label{s7th max eq2.1}
\max_{t\in[0,T]}\mE|x_3^\e(t)|^4_H \leq \cC(x_0)
\e^4,
\end{equation}
\begin{equation}\label{s7th max eq3}
|x_4^\e(\cd)|_{C_\dbF([0,T];L^2(\Om;H))} \leq
\cC(x_0) \e,
\end{equation}
\begin{equation}\label{sqqeq0.2}
|x_5^\e(\cd)|_{C_\dbF([0,T];L^2(\Om;H))}^2
=o(\e^2),\qq\hb{as } t\to0,
\end{equation}
\begin{equation}\label{2z10}
\begin{array}{ll}
\ds\Big| \int_0^1(1-\si) \big( g_{xx}\big(t,
\bar x(t)+ \si x_1^\e(t), u^\e(t)\big) -
g_{xx}\big(t,\bar
x(t),u^\e(t)\big)\big)d\si\Big|
_{\cL(H)}\\\ns\ds\le \cC\(\int_0^1 \big|
g_{xx}\big(t,\bar x(t) + \si x_1^\e(t), \bar u
(t)\big)- g_{xx}\big(t,\bar x(t), \bar
u(t)\big)\Big| _{\cL(H)}d\si+\chi_{E_\e}(t)\),
\end{array}
\end{equation}
\begin{equation}\label{220160221z10}
\begin{array}{ll}
\ds\Big| \int_0^1(1-\si)\big( h_{xx}\big( \bar
x(T)+\si x_1^\e(T) \big) -h_{xx}\big(\bar
x(T)\big) \big)d\si\Big| _{\cL(H)}\\\ns\ds\le
\cC\(\int_0^1 \big| h_{xx}(\bar x(T) + \si
x_1^\e(T))- h_{xx}(\bar x(T))\Big|
_{\cL(H)}d\si+\chi_{E_\e}(t)\),
\end{array}
\end{equation}
and
\begin{equation}\label{s7th max eq0.2}
|x_1^\e - x_2^\e -
x_3^\e|_{L^\infty_\dbF(0,T;L^2(\Omega;H))}=o(\e),\qq\hb{as
} \e\to0.
\end{equation}

\ms

{\bf Step 2}.  We need to compute the value of
$\cJ(u^\e(\cd)) - \cJ(\bar
u(\cd))$.\vspace{-0.1cm}
\begin{equation}\label{s7eq3}
\begin{array}{ll}\ds
 \cJ(u^\e(\cd)) -  \cJ(\bar u(\cd))\\
\ns\ds = \mE\int_0^T \big[g(t,x^\e(t),u^\e(t)) -
g(t,\bar x(t),\bar u(t))\big]
dt + \mE h\big(x^\e(T)\big) - \mE h\big(\bar x(T)\big)\\
\ns\ds  = \mE \int_0^T \Big\{
\chi_{E_\e}(t)\d g(t) + \big\langle g_x(t,\bar
x(t),
u^\e(t)),x_1^\e(t) \big\rangle_H  \\
\ns\ds \q + \int_0^1 \big\langle(1-\si)
g_{xx}\big(t, \bar x(t) + \si
x_1^\e(t),u^\e(t)\big)x_1^\e(t),x_1^\e(t)\big\rangle_H
d\si \Big\}dt \\
\ns\ds  \q +\mE\big\langle h_x(\bar x(T)),
x_1^\e(T)\big\rangle_H + \mE\int_0^1
\big\langle (1-\si) h_{xx}\big( \bar x(T)+ \si
x_1^\e(T)\big)x_1^\e(T), x_1^\e(T)\big\rangle_H
d\si.
\end{array}
\end{equation}
This, together with the definition of
$x_i^\e(\cd)$ ($i=1,2,3,4,5$), yields that
\begin{equation}\label{s7eq4}
\!\!\begin{array}{ll}\ds
 \cJ(u^\e(\cd)) -  \cJ(\bar u(\cd))\\
\ns\ds = \mE\int_0^T \Big\{ \chi_{E_\e}(t)\d
g(t) + \big\langle \d g_1(t),x_1^\e(t)
\big\rangle_H \chi_{E_\e}(t) + \big\langle
g_1(t), x_2^\e(t) + x_3^\e(t) \big\rangle_H +
\big\langle g_1(t),
x_5^\e(t) \big\rangle_H \\
\ns\ds\q \!\! + \int_0^1\big\langle (1-\si)
\big[ g_{xx}\big(t, \bar x(t) + \si x_1^\e(t),
u^\e(t)\big) - g_{xx}\big(t,\bar
x(t),u^\e(t)\big)\big]x_1^\e(t), x_1^\e(t)
\big\rangle_H d \si \\
\ns\ds \q\!\! + \frac{1}{2}\big\langle \d
g_{11}(t)x_1^\e(t), x_1^\e(t) \big\rangle_H
\chi_{E_\e}(t)\! +\! \frac{1}{2}\big\langle
g_{11}(t)x_2^\e(t),x_2^\e(t) \big\rangle_H \!+\!
\frac{1}{2} \big\langle g_{11}(t)x_4^\e(t),
x_1^\e(t) \!+\! x_2^\e(t) \big\rangle_H \Big\}dt \\
\ns\ds \q\!\! + \mE\big\langle h_x\big(\bar
x(T)\big), x_2^\e(t) \!+\! x_3^\e(t)
\big\rangle_H \!+\! \mE \big\langle h_x\big(\bar
x(T)\big), x_5^\e(t) \big\rangle_H \!+\!
\frac{1}{2}\mE\big\langle h_{xx}\big(\bar
x(T)\big) x_2^\e(t),x_2^\e(t)
\big\rangle_H \\
\ns\ds \q\!\! + \frac{1}{2}\mE \big\langle
h_{xx}\big( \bar x(T) \big)x_4^\e(T), x_1^\e(T)
+
x_2^\e(T)  \big\rangle_H \\
\ns\ds \q \!\!+ \mE\int_0^1 \big\langle
 (1-\si)\big[ h_{xx}\big( \bar x(T) + \si x_1^\e(T) \big) - h_{xx}\big(\bar x(T)\big)
\big]x_1^\e(T),x_1^\e(T) \big\rangle_H d \si.
\end{array}
\end{equation}
For $\ae t\in[0,T]$, we find that
\begin{equation}\label{2z10-990}
\begin{array}{ll}
\ds\big|\!\big|\int_0^1 (1-\si) \big[
g_{xx}\big(t, \bar x(t) + \si x_1^\e(t),
u^\e(t)\big) - g_{xx}\big(t,\bar
x(t),u^\e(t)\big)\big]d\si\big|\!\big|_{\cL(H\times
H,\;H)}\\\ns\ds =\big|\!\big| \int_0^1 (1-\si)
\[g_{xx}\big(t,\bar x(t) + \si x_1^\e(t), \bar u (t)\big)-
g_{xx}\big(t,\bar x(t), \bar u(t)\big)\]d\si
\\\ns\ds
 \q+\int_0^1 (1-\si) \chi_{E_\e}(t)g_{xx}\big(t,\bar
x(t) + \si x_1^\e(t),  u
(t)\big)d\si+\chi_{E_\e}(t)g_{xx}\big(t,\bar
x(t), u(t)\big) \big|\!\big|_{\cL(H\times
H,\;H)} d\si
\\\ns\ds \le \cC\[\int_0^1
\big|\!\big|g_{xx}\big(t,\bar x(t) + \si
x_1^\e(t), \bar u (t)\big)- g_{xx}\big(t,\bar
x(t), \bar u(t)\big)\big|\!\big|_{\cL(H\times
H,\;H)}d\si +\chi_{E_\e}(t)\].
\end{array}
\end{equation}
By (\ref{s7eq4}), noting \eqref{s7th max eq0.1},
\eqref{s7th max eq2}, \eqref{s7th max eq2.1},
\eqref{s7th max eq3}, \eqref{sqqeq0.2} and
(\ref{2z10}), and using the continuity of both
$h_{xx}(x)$ and $g_{xx}(x)$ with respect to $x$,
we end up with
\begin{equation}\label{s7eq5}
\begin{array}{ll}
\ds
\cJ(u^\e(\cd)) -  \cJ(\bar u(\cd))\\
\ns\ds = \mE \int_0^T\[ \big\langle
g_1(t),x_2^\e(t) + x_3^\e(t) \big\rangle_H +
\frac{1}{2}\big\langle
g_{11}(t)x_2^\e(t),x_2^\e(t) \big\rangle_H +
\chi_{E_\e}(t)\d g(t)
\] dt \\
\ns\ds\q + \mE \big\langle h_x\big(\bar
x(T)\big), x_2^\e(T)+x_3^\e(T) \big\rangle_H +
\frac{1}{2}\mE\big\langle h_{xx}\big(\bar
x(T)\big)x_2^\e(t),x_2^\e(t) \big\rangle_H +
o(\e).
\end{array}
\end{equation}

In the sequel, we shall get rid of $x_2^\e(\cd)$
and $x_3^\e(\cd)$ in \eqref{s7eq5} by solutions
to the equations \eqref{ch-10-bsystem1} and
\eqref{op-bsystem3}. By the definition of the
transposition solution to the equation
\eqref{ch-10-bsystem1} (with $y_T$ and
$f(\cd,\cd,\cd)$ given by \eqref{zv1}), we
obtain that
\begin{equation}\label{5.26-eq6}
-\mE\big\langle h_x(\bar
x(T))),x_2^\e(T)\big\rangle_H - \mE \int_0^T
\big\langle g_1(t),x_2^\e(t)\big\rangle_H dt =
\mE \int_0^T\big\langle Y(t), \d
b(t)\big\rangle_H\chi_{E_\e}(t) dt
\end{equation}
and
\begin{equation}\label{5.26-eq7}
\begin{array}{ll}
\ds-\mE\big\langle h_x(\bar
x(T))),x_3^\e(T)\big\rangle_H - \mE \int_0^T
\big\langle g_1(t),x_3^\e(t)\big\rangle_H dt
\\
\ns\ds = \mE \int_0^T \Big\{  \frac{1}{2}\[
\big\langle y(t),a_{11}(t)\big(x_2^\e(t),
x_2^\e(t)\big) \big\rangle_H + \big\langle Y(t),
b_{11}(t)\big(
x_2^\e(t), x_2^\e(t)\big) \big\rangle_H \] \\
\ns\ds \hspace{1.8cm} + \chi_{E_\e}(t)
\[ \big\langle y(t),\d a(t) \big\rangle_H +
\big\langle Y,\d b_1(t)x_2^\e(t) \big\rangle_H
 \]\Big\}dt.
\end{array}
\end{equation}
According to \eqref{s7eq5}--\eqref{5.26-eq7}, we
conclude that
\begin{equation}\label{5.26-eq8}
\begin{array}{ll} \ds
\cJ(u^\e(\cd)) -  \cJ(\bar u(\cd))\\
\ns\ds = \frac{1}{2}\mE\int_0^T\[ \big\langle
g_{11}(t)x_2^\e(t), x_2^\e(t)\big\rangle_H -
\big\langle y(t),a_{11}(t)\big(x_2^\e(t),
x_2^\e(t)\big) \big\rangle_H \\
\ns\ds \q - \big\langle Y,
b_{11}(t)\big(x_2^\e(t),
x_2^\e(t)\big)\big\rangle_H
\]dt + \mE\int_0^T \chi_{E_\e}(t)\[ \d g(t) - \big\langle
y(t),\d a(t)\big\rangle_H \\
\ns\ds \q-\big\langle Y(t),\d b(t) \big\rangle_H
\]dt + \frac{1}{2}\mE \big\langle h_{xx}\big(\bar
x(T)\big)x_2^\e(T), x_2^\e(T) \big\rangle_H +
o(\e),\qq\hb{as }\e\to0.
\end{array}
\end{equation}
By the definition of the relaxed transposition
solution to the equation \eqref{op-bsystem3}
(with $P_T$, $J(\cd)$, $K(\cd)$ and $F(\cd)$
given by \eqref{MP2-eq9}), we obtain that
\begin{equation}\label{5.26-eq9}
\begin{array}{ll}\ds
-\mE\big\langle h_{xx}\big(\bar x(T)\big)
x_2^\e(T), x_2^\e(T) \big\rangle_H + \mE\int_0^T
\big\langle \dbH_{xx}\big(t,\bar
x(t),\bar u(t),y(t),Y(t)\big) x_2^\e(t), x_2^\e(t) \big\rangle_H dt\\
\ns\ds = \dbE\int_0^T \chi_{E_\e}(t)\big\langle
b_1(t)x_2^\e(t), P(t)^*\d b(t)\big\rangle_{H} dt
+ \dbE\int_0^T \chi_{E_\e}(t)\big\langle P(t)\d
b(t),
b_1(t)x_2^\e(t)\big\rangle_{H} dt\\
\ns\ds \q  + \dbE\int_0^T
\chi_{E_\e}(t)\big\langle P(t)\d b(t), \d
b(t)\big\rangle_{H} dt + \dbE\int_0^T
\chi_{E_\e}(t)\big\langle
\d b(t),\widehat Q^{(0)}(0,0,\chi_{E_\e}\d b)(t)\big\rangle_{H} dt \\
\ns\ds \q + \dbE\int_0^T
\chi_{E_\e}(t)\big\langle Q^{(0)}(0,0,\d
b)(t),\d b(t)\big\rangle_{H} dt.
\end{array}
\end{equation}

Now, we estimate the terms in the right hand
side of \eqref{5.26-eq9}. By \eqref{s7th max
eq2}, we have
\begin{equation}\label{s7eq9.1}
\begin{array}{ll}\ds
\Big|\dbE\int_0^T \chi_{E_\e}(t)\big\langle
b_1(t)x_2^\e(t), P(t)^*\d b(t)\big\rangle_{H}
dt+\dbE\int_0^T \chi_{E_\e}(t)\big\langle P(t)\d
b(t), b_1(t)x_2^\e(t)\big\rangle_{H} dt\Big|=
o(\e).
\end{array}
\end{equation}

In what follows, for any $\tau\in [0,T)$, we
choose $E_{\e}=[\tau,\tau+\e]\subset [0,T]$.

By Lemma \ref{5.13-prop1}, we can find a
sequence $\{\beta_n\}_{n=1}^\infty$ such that
$\b_n\in\cH_n$ (Recall \eqref{cH} for the
definition of $\cH_n$) and $
\lim_{n\to\infty}\beta_n = \d b $ in $
L^2_\dbF(0,T;L^4(\Om;H))$. Hence, for some
positive constant $C(x_0)$ (depending on $x_0$),
\begin{equation}\label{10.9qqq3}
|\beta_n|_{L^2_\dbF(0,T;L^4(\Om;H))}\le
C(x_0)<\infty,\qq\forall\;n\in\dbN,
\end{equation}
and there is a subsequence
$\{n_k\}_{k=1}^\infty\subset \{
n\}_{n=1}^\infty$ such that
\begin{equation}\label{s7eq9.2-11}
\lim_{k\to\infty} |\b_ {n_k}(t)-\d
b(t)|_{L^4_{\cF_t}(\Om;H)} = 0\q\mbox{ for }\ae
t\in [0,T].
\end{equation}

Denote by $Q^{n_k}$ and $\widehat Q^{n_k}$ the
corresponding pointwise defined linear
operators from $\cH_{n_k}$ to $L^2_\dbF(0,T;$
$L^{\frac{4}{3}}(\Om;H))$, given in Lemma
\ref{10.1th}.

Consider the following equation:
\begin{equation}\label{s7fsystem3.1x}
\left\{
\begin{array}{lll}\ds
dx_{2,n_k}^{\e} = \big[Ax_{2,n_k}^{\e} +
a_1(t)x_{2,n_k}^{\e} \big]dt + \big[ b_1(t)
x_{2,n_k}^{\e} + \chi_{E_{\e}} (t)\b_ {n_k}(t)
\big] dW(t) &\mbox{ in
}(0,T],\\
\ns\ds x_{2,n_k}^{\e}(0)=0.
\end{array}
\right.\vspace{-0.2cm}
\end{equation}
We have \vspace{-0.2cm}
\begin{equation}\label{s7th max eq1x}
\begin{array}{ll}\ds
\mE|x_{2,n_k}^{\e}(t)|^4_H \\
\ns\ds = \mE\Big| \int_0^t
S(t-s)a_1(s)x_{2,n_k}^{\e}(s) ds + \int_0^t
S(t-s)b_1(s)x_{2,n_k}^{\e}(s)
dW(s) \\
\ns\ds \qq +
\int_0^t S(t-s)\chi_{E_{\e}}(s)\b_ {n_k}(s) dW(s)\Big|_H^4\\
\ns\ds  \leq C \bigg[\mE\Big| \int_0^t
S(t-s)a_1(s)x_{2,n_k}^{\e}(s) ds \Big|_H^4
+ \mE\Big|\int_0^t S(t-s)b_1(s)x_{2,n_k}^{\e}(s) dW(s) \Big|_H^4 \\
\ns\ds \q  + \mE\Big| \int_0^t S(t-s)\chi_{E_{\e}}(s)\b_ {n_k}(s) dW(s)\Big|_H^4\bigg] \\
\ns\ds   \leq C\[\int_0^t
\mE|x_{2,n_k}^{\e}(s)|_H^4 ds  + \e
\int_{E_{\e}}\mE|\beta_{n_k}(s)|_{H}^4ds
\].
\end{array}
\end{equation}
By \eqref{10.9qqq3} and thanks to Gronwall's
inequality, \eqref{s7th max eq1x} leads to
\begin{equation}\label{s7th max eq2x}
|x_{2,n_k}^{\e}(\cd)|^4_{L^\infty_\dbF(0,T;L^4(\Om;H))}
\leq C(x_0,k)\e^2.
\end{equation}
Here and henceforth, $C(x_0,k)$ is a generic
constant (depending on $x_0$,
 $k$, $T$, $A$ and $C_L$), which may
be different from line to line. For any fixed $
k\in\dbN$, since $Q^{n_k}\b_ {n_k}\in
L^2_{\dbF}(0,T;L^{\frac{4}{3}}(\Om;H))$, by
\eqref{s7th max eq2x}, we find that
\begin{equation}\label{s7eq9.3}
\begin{array}{ll}\ds
\Big|\dbE\int_0^T \chi_{E_{\e}}(t)\big\langle
\big(Q^{n_k} \b_
{n_k}\big)(t),x_{2,n_k}^{\e}(t)\big\rangle_{H}
dt \Big|\\
\ns\ds \leq
|x_{2,n_k}^{\e}(\cd)|_{L^\infty_\dbF(0,T;L^4(\Om;H))}
\int_{E_{\e}}\big|\big(Q^{n_k} \b_ {n_k}\big)(t)\big|_{L^{\frac{4}{3}}_{\cF_t}(\Om;H)}dt\\
\ns\ds \leq C(x_0,k)\sqrt{{\e}}
\int_{E_{\e}}\big|\big(Q^{n_k} \b_
{n_k}\big)(t)\big|_{L^{\frac{4}{3}}_{\cF_t}(\Om;H)}dt=
o({\e}), \qq\hbox{as }\e\to0.
\end{array}
\end{equation}
Similarly,
\begin{equation}\label{s7eq9.3x}
\Big|\dbE\int_0^T \chi_{E_{\e}}(t)\big\langle
x_{2,n_k}^{\e}(t),\big(\widehat Q^{n_k}\b_
{n_k}\big)(t)\big\rangle_{H} dt \Big| = o({\e}),
\qq\hbox{as }\e\to0.
\end{equation}

From \eqref{10.9eq2} in Theorem \ref{10.1th},
and noting that both $Q^{n_k}$ and $\widehat
Q^{n_k}$ are pointwise defined, we arrive at
the following equality:
\begin{equation}\label{wws1}
\begin{array}{ll}
\ds \mE \int_{0}^T \big\langle
\chi_{E_{\e}}(t)\b_ {n_k}(t), \widehat
Q^{(0)}(0,0,\chi_{E_{\e}} \b_{n_k})(t)
\big\rangle_{H}dt +  \mE \int_{0}^T \big\langle
Q^{(0)}(0,0,\chi_{E_{\e}} \b_ {n_k}) (t), \chi_{E_{\e}}\b_ {n_k}(t) \big\rangle_{H}dt \\
\ns\ds =\mE \int_{0}^T
\chi_{E_{\e}}\[\big\langle \big(Q^{n_k} \b_
{n_k}\big)(t), x_{2,n_k}^{\e} (t)
\big\rangle_{H}+\big\langle x_{2,n_k}^{\e} (t),
\big(\widehat Q^{n_k} \b_ {n_k}\big)(t)
\big\rangle_{H}\]dt.
\end{array}
\end{equation}
Hence,
\begin{equation}\label{wws2}
\begin{array}{ll}\ds
\mE \int_{0}^T \big\langle \chi_{E_{\e}}(t)\d
b(t), \widehat Q^{(0)}(0,0,\chi_{E_{\e}}\d b)(t)
\big\rangle_{H}dt + \mE \int_{0}^T \big\langle
 Q^{(0)}(0,0,\chi_{E_{\e}}\d
b) (t), \chi_{E_{\e}}(t)\d
b(t) \big\rangle_{H}dt \\
\ns\ds \q\!\!-\mE \int_{0}^T \chi_{E_{\e}}(t)
\[\big\langle \big(Q^{n_k} \b_ {n_k}\big)(t),
x_{2,n_k}^{\e}(t) \big\rangle_{H}+\big\langle
x_{2,n_k}^{\e} (t), \big(\widehat Q^{n_k} \b_
{n_k}\big)(t)
\big\rangle_{H}\]dt\\
\ns\ds = \mE \int_{0}^T \big\langle
\chi_{E_{\e}}(t)\d b(t), \widehat
Q^{(0)}(0,0,\chi_{E_{\e}} \d b)(t)
\big\rangle_{H}dt + \mE \int_{0}^T \big\langle
 Q^{(0)}(0,0,\chi_{E_{\e}} \d
b) (t), \chi_{E_{\e}}(t)\d
b(t) \big\rangle_{H}dt \\
\ns\ds \q\!\! -\mE\!\!\int_{0}^T\!\!\big\langle
\chi_{E_{\e}}\!(t)\b_ {n_k}\!(t), \widehat
Q^{(0)}(0,0,\chi_{E_{\e}}\b_{n_k})(t)\big\rangle_{H}dt
\! - \!\mE\!\!\int_{0}^T\!\!\big\langle
Q^{(0)}(0,0,\chi_{E_{\e}} \b_{n_k} ) (t),
\chi_{E_{\e}}\!(t)\b_{n_k}\!(t)
\big\rangle_{H}dt.
\end{array}
\end{equation}
It is easy to see that
\begin{equation}\label{s7eq9.2xx}
\!\!\!\begin{array}{ll}\ds\Big|\mE \int_{0}^T
\big\langle \chi_{E_{\e}}(t)\d b(t), \widehat
Q^{(0)}(0,0,\chi_{E_{\e}} \d b)(t)
\big\rangle_{H}dt  - \mE \int_{0}^T \big\langle
\chi_{E_{\e}}(t)\b_ {n_k}(t), \widehat
Q^{(0)}(0,0,\chi_{E_{\e}} \b_ {n_k})(t)
\big\rangle_{H}dt \Big|\\
\ns\ds \leq \!\Big|\mE \!\int_{0}^T\!
\big\langle \chi_{E_{\e}}(t)\d b(t), \widehat
Q^{(0)}(0,0,\chi_{E_{\e}} \d b)(t)
\big\rangle_{H}dt - \mE\! \int_{0}^T\!
\big\langle \chi_{E_{\e}}(t)\d b(t), \widehat
Q^{(0)}(0,0,\chi_{E_{\e}} \b_ {n_k})(t)
\big\rangle_{H}dt \Big|\\
\ns\ds \q\!\!\!\!\! + \Big|\mE\!
\int_{0}^T\!\!\! \big\langle \chi_{E_{\e}}(t)\d
b(t), \widehat Q^{(0)}(0, 0, \chi_{E_{\e}} \b_
{n_k})(t) \big\rangle_{H}dt\! - \!\mE\!
\int_{0}^T\!\! \big\langle \chi_{E_{\e}}(t)\b_
{n_k}(t), \widehat Q^{(0)}(0,0,\chi_{E_{\e}} \b_
{n_k})(t) \big\rangle_{H}dt \Big|.
\end{array}
\end{equation}
From \eqref{s7eq9.2-11} and the density of the
Lebesgue points, we find that for $\ae\tau\in
[0,T)$, it holds that
\begin{equation}\label{s7eq9.3-1}
\begin{array}{ll}\ds
\lim_{k\to\infty}\lim_{\e\to
0}\frac{1}{{\e}}\Big|\mE \int_{0}^T \big\langle
\chi_{E_{\e}}(t)\d b(t),\widehat
Q^{(0)}(0,0,\chi_{E_{\e}} \d b)(t)
\big\rangle_{H}dt
\\
\ns\ds \qq\qq\q - \mE \int_{0}^T \big\langle
\chi_{E_{\e}}(t)\d b(t),\widehat Q^{(0)}(0,0,\chi_{E_{\e}} \b_ {n_k})(t) \big\rangle_{H}dt \Big|\\
\ns\ds \leq \lim_{k\to\infty}\lim_{\e\to
0}\frac{1}{{\e}}\[\int_0^T \chi_{E_{\e}}(t)
\(\mE|\d b(t)|^4_{H}\)^{\frac{1}{2}} dt
\Big]^{\frac{1}{2}}|\widehat
Q^{(0)}(0,0,\chi_{E_{\e}} (\d b-\b_ {n_k}))
|_{L^2_\dbF(0,T;L^{\frac{4}{3}}(\Om;H))}
\\
\ns\ds \leq C\lim_{k\to\infty}\lim_{\e\to
0}\frac{1}{{\e}}
\[\int_0^T \chi_{E_{\e}}(t) \(\mE|\d b(t) |^4_{H}\)^{\frac{1}{2}} dt
\Big]^{\frac{1}{2}}\big|\chi_{E_{\e}} (\d b-\b_ {n_k})\big|_{L^2_\dbF(0,T;L^4(\Om;H))}\\
\ns\ds \leq C\lim_{k\to\infty}\lim_{\e\to
0}\frac{|\d
b(\tau)|_{L^4_{\cF_\tau}(\Om;H)}}{\sqrt{{\e}}}\[\int_0^T
\chi_{E_{\e}}(t) \(\mE|\d b(t) - \b_
{n_k}(t)|^4_{H}\)^{\frac{1}{2}} dt
\Big]^{\frac{1}{2}}\\
\ns\ds = C\lim_{k\to\infty}\lim_{\e\to 0}|\d
b(\tau)|_{L^4_{\cF_\tau}(\Om;H)}\[\frac{1}{{\e}}\int_\tau^{\tau+{\e}}
|\d b(t) - \b_ {n_k}(t)|_{L^4_{\cF_t}(\Om;H)}^2
dt
\Big]^{\frac{1}{2}}\\
\ns\ds = C\lim_{k\to\infty}|\d
b(\tau)|_{L^4_{\cF_\tau}(\Om;H)}|\d b(\tau) -
\b_ {n_k}(\tau)|_{L^4_{\cF_\tau}(\Om;H)}
\\\ns\ds= 0.
\end{array}
\end{equation}
Similarly,
\begin{equation}\label{s7eq9.3-1x}
\begin{array}{ll}\ds
\lim_{k\to\infty}\lim_{\e\to
0}\frac{1}{{\e}}\Big|\mE \int_{0}^T \big\langle
\chi_{E_{\e}}(t)\d b(t), \widehat
Q^{(0)}(0,0,\chi_{E_{\e}} \b_ {n_k})(t)
\big\rangle_{H}dt \\
\ns\ds \qq\qq\q - \mE  \int_{0}^T \big\langle
\chi_{E_{\e}}(t)\b_ {n_k}(t), \widehat
Q^{(0)}(0,0,\chi_{E_{\e}} \b_ {n_k})(t)
\big\rangle_{H}dt \Big|\\
\ns\ds \leq \lim_{k\to\infty}\lim_{\e\to
0}\frac{1}{{\e}}\big|\widehat
Q^{(0)}(0,0,\chi_{E_{\e}} \b_
{n_k})\big|_{L^2_\dbF(0,T;L^{\frac{4}{3}}(\Om;H))}
\[\int_0^T \chi_{E_{\e}}(t) \(\mE|\d b(t) - \b_ {n_k}(t)|^4_{H}\)^{\frac{1}{2}} dt
\Big]^{\frac{1}{2}}\\
\ns\ds \leq C\lim_{k\to\infty}\lim_{\e\to
0}\frac{1}{{\e}}\big|\chi_{E_{\e}} \b_
{n_k}\big|_{L^2_\dbF(0,T;L^4(\Om;H))}
\[\int_0^T \chi_{E_{\e}}(t) \(\mE|\d b(t) - \b_ {n_k}(t)|^4_{H}\)^{\frac{1}{2}} dt
\Big]^{\frac{1}{2}}\\
\ns\ds \leq C\lim_{k\to\infty}\lim_{\e\to
0}\frac{1}{{\e}}\Big\{\big|\chi_{E_{\e}} \d
b\big|_{L^2_\dbF(0,T;L^4(\Om;H))}
\[\int_0^T \chi_{E_{\e}}(t) \(\mE|\d b(t) - \b_ {n_k}(t)|^4_{H}\)^{\frac{1}{2}} dt
\Big]^{\frac{1}{2}}\\\ns\ds\qq\qq\qq\qq\q+\int_0^T \chi_{E_{\e}}(t) \(\mE|\d b(t) - \b_ {n_k}(t)|^4_{H}\)^{\frac{1}{2}} dt\Big\}\\
\ns\ds \leq C\lim_{k\to\infty}\lim_{\e\to
0}\Big\{\frac{|\d
b(\tau)|_{L^4_{\cF_\tau}(\Om;H)}}{\sqrt{{\e}}}\[\int_0^T
\chi_{E_{\e}}(t) \(\mE|\d b(t) - \b_
{n_k}(t)|^4_{H}\)^{\frac{1}{2}} dt
\Big]^{\frac{1}{2}}\\\ns\ds\qq\qq\qq\q+\frac{1}{{\e}}\int_0^T
\chi_{E_{\e}}(t) \(\mE|\d b(t) - \b_
{n_k}(t)|^4_{H}\)^{\frac{1}{2}} dt
\Big\}\\
\ns\ds = C\lim_{k\to\infty}\lim_{\e\to
0}\Big\{|\d
b(\tau)|_{L^4_{\cF_\tau}(\Om;H)}\[\frac{1}{{\e}}\int_\tau^{\tau+{\e}}
|\d b(t) - \b_ {n_k}(t)|_{L^4_{\cF_t}(\Om;H)}^2
dt
\Big]^{\frac{1}{2}}\\\ns\ds\qq\qq\qq\q+\frac{1}{{\e}}\int_\tau^{\tau+{\e}}
|\d b(t) - \b_ {n_k}(t)|_{L^4_{\cF_t}(\Om;H)}^2
dt
\Big\}\\
\ns\ds = C\lim_{k\to\infty}\big[|\d
b(\tau)|_{L^4_{\cF_\tau}(\Om;H)}|\d b(\tau) -
\b_ {n_k}(\tau)|_{L^4_{\cF_\tau}(\Om;H)}+|\d
b(\tau) - \b_
{n_k}(\tau)|_{L^4_{\cF_\tau}(\Om;H)}^2\big]
\\\ns\ds= 0.
\end{array}
\end{equation}
From \eqref{s7eq9.2xx}--\eqref{s7eq9.3-1x}, we
find that
\begin{equation}\label{s7eq9.3-1xx}
\begin{array}{ll}\ds
\lim_{k\to\infty} \lim_{\e\to
0}\frac{1}{{\e}}\Big|\mE \int_{0}^T \big\langle
\chi_{E_{\e}}(t)\d b(t), \widehat
Q^{(0)}(0,0,\chi_{E_{\e}} \d b)(t)
\big\rangle_{H}dt \\
\ns\ds \qq - \mE  \int_{0}^T \big\langle
\chi_{E_{\e}}(t)\b_ {n_k}(t), \widehat
Q^{(0)}(0,0,\chi_{E_{\e}} \b_ {n_k})(t)
\big\rangle_{H}dt \Big| = 0.
\end{array}
\end{equation}
By a similar argument, we obtain that
\begin{equation}\label{s7eq9.3-1xxx}
\begin{array}{ll}\ds
\lim_{k\to\infty} \lim_{\e\to
0}\frac{1}{{\e}}\Big|\mE \int_{0}^T \big\langle
Q^{(0)}(0,0,\chi_{E_{\e}} \d
b)(t),\chi_{E_{\e}}(t)\d b(t)
\big\rangle_{H}dt \\
\ns\ds \qq - \mE  \int_{0}^T \big\langle
Q^{(0)}(0,0,\chi_{E_{\e}} \b_
{n_k})(t),\chi_{E_{\e}}(t)\b_ {n_k}(t)
\big\rangle_{H}dt \Big| = 0.
\end{array}
\end{equation}

From \eqref{s7eq9.3}--\eqref{wws2} and
\eqref{s7eq9.3-1xx}--\eqref{s7eq9.3-1xxx}, we
obtain that
\begin{equation}\label{s7eq9.2-111}
 \ba{ll}\ds
\Big|\dbE\int_0^T \chi_{E_{\e}}(t)\big\langle \d
b(t),\widehat Q^{(0)}(0,0,\chi_{E_{\e}}\d
b)(t)\big\rangle_{H} dt + \dbE\int_0^T
\chi_{E_{\e}}(t)\big\langle Q^{(0)}(0,0,\d
b)(t),\d b(t)\big\rangle_{H} dt\Big|
\\\ns\ds=o({\e}), \qq\hbox{as }\e\to0.
 \ea
\end{equation}

Combining \eqref{5.26-eq8}, \eqref{5.26-eq9},
\eqref{s7eq9.1} and \eqref{s7eq9.2-111}, we end
up with
 $$
 \ba{ll}\ds
\cJ(u^{\e}(\cd)) -  \cJ(\bar u(\cd))\\\ns\ds =
\mE\int_0^T \[
 \d g(t) - \big\langle
y(t),\d a(t)\big\rangle_H -\big\langle Y(t),\d
b(t) \big\rangle_H - \frac{1}{2}\big\langle
P(t)\d b(t), \d b(t) \big\rangle_H
 \]\chi_{E_{\e}}(t)dt + o({\e}).
 \ea
 $$
Since $\bar u(\cd)$ is the optimal control,
$\cJ(u^{\e}(\cd)) - \cJ(\bar u(\cd))\geq 0$.
Thus,
\begin{equation}\label{s7eq11}
\mE\int_0^T \chi_{E_{\e}}(t)\[
  \big\langle
y(t),\d a(t)\big\rangle_H +\big\langle Y(t),\d
b(t) \big\rangle_H -\d g(t)+
\frac{1}{2}\big\langle P(t)\d b(t), \d b(t)
\big\rangle_H
 \]dt \leq o({\e}),
\end{equation}
as $\e\to0$.

Finally, by \eqref{s7eq11}, we obtain
\eqref{MP2-eq1}. This completes the proof of
Theorem \ref{maximum p2}.
\endpf

\ms

It is worth to mention that, the stochastic transposition method has some other applications, say it can be used to establish the equivalence
between the existence of optimal feedback
operator for infinite dimensional stochastic linear quadratic control
problems with random coefficients and the
solvability of the corresponding
operator-valued, backward stochastic Riccati
equations (See \cite{LZ5} for more details).




\end{document}